\documentclass[11pt]{article}

\pdfoutput=1

\usepackage{endnotes}
\let\footnote=\endnote
\let\enotesize=\normalsize
\def\notesname{Endnotes}%
\def\makeenmark{\hbox to1.275em{\theenmark.\enskip\hss}}
\def\enoteformat{\rightskip0pt\leftskip0pt\parindent=1.275em
  \leavevmode\llap{\makeenmark}}

\usepackage{graphicx}
\usepackage{ifthen}
\usepackage{amsmath}
\usepackage{amssymb}
\usepackage{amsthm}
\usepackage{fullpage}
\usepackage{setspace}
\usepackage{charter}
\newtheorem{theorem}{Theorem}
\newtheorem{lemma}{Lemma}
\newtheorem{proposition}{Proposition}
\usepackage{comment}
\usepackage{subcaption}
\usepackage{accents}
\usepackage{enumerate}
\usepackage{hyperref}
\usepackage{fix-cm}

\newcommand{\na}{N}
\newcommand{\nb}{M}
\newcommand{\ca}{c}
\newcommand{\cb}{C}
\newcommand{\cm}{\triangle C}
\newcommand{\ka}{k}
\newcommand{\rv}{\alpha}
\newcommand{\rva}[1][]{%
\ifthenelse{\equal{#1}{}}{\alpha_y}{\alpha_{y,#1}}%
}
\newcommand{\rvabar}[1][]{%
\ifthenelse{\equal{#1}{}}{\bar{\alpha}_y}{\bar{\alpha}_{y,#1}}%
}
\newcommand{\rvb}[1][]{%
\ifthenelse{\equal{#1}{}}{\alpha_x}{\alpha_{x,#1}}%
}
\newcommand{\rvbbar}[1][]{%
\ifthenelse{\equal{#1}{}}{\bar{\alpha}_x}{\bar{\alpha}_{x,#1}}%
}

\newcommand{\qbpj}[1][]{%
\ifthenelse{\equal{#1}{}}{x_i}{x_{i,#1}}%
}
\newcommand{\qapi}{y_{j}}
\newcommand{\qfpi}[1][]{%
\ifthenelse{\equal{#1}{}}{f_j}{f_{j,#1}}%
}
\newcommand{\qbp}{\mathbf{x}}
\newcommand{\qfp}{\mathbf{f}}
\newcommand{\qap}{\mathbf{y}}
\newcommand{\qbpn}{\mathbf{x}_{-i}}
\newcommand{\qfpn}{\mathbf{f}_{-j}}
\newcommand{\qapn}{\mathbf{y}_{-j}}
\newcommand{\qbpjj}{x_{i'}}

\newcommand{\qapii}{y_{j'}}

\newcommand{\qapiopt}{y_j}
\newcommand{\qapopt}{\mathbf{y}}

\newcommand{\qapiiopt}{y_{j'}}

\newcommand{\pfpis}{\phi_{j}^{(s)}}

\newcommand{\qf}{f}
\newcommand{\qb}{x}
\newcommand{\qa}{y}

\newcommand{\pfpif}{\phi_{j}}
\newcommand{\pfpjf}{\psi_{i}}

\newcommand{\pf}{p_f}

\newcommand{\Qf}{F}
\newcommand{\Qb}{X}
\newcommand{\Qa}{Y}

\newcommand{\ubar}[1]{\underaccent{\bar}{#1}}

\newboolean{showcomments}
\setboolean{showcomments}{true}
\newcommand{\todo}[1]{\ifthenelse{\boolean{showcomments}}
{\textcolor{red}{(todo: \textit{#1}) }} {}}
\newcommand{\addref}[0]{\ifthenelse{\boolean{showcomments}}
{\textcolor{red}{(add ref) }} {}}
\newcommand{\addcite}[0]{\ifthenelse{\boolean{showcomments}}
{\textcolor{red}{(add cite) }} {}}
\newcommand{\addcites}[0]{\ifthenelse{\boolean{showcomments}}
{\textcolor{red}{(add cites) }} {}}
\newcommand{\desmond}[1]{\ifthenelse{\boolean{showcomments}}
{\textcolor{red}{(Desmond says: \textit{#1}) }}{}}
\newcommand{\adam}[1]{\ifthenelse{\boolean{showcomments}}
{\textcolor{red}{(Adam says: \textit{#1}) }}{}}
\newcommand{\anish}[1]{\ifthenelse{\boolean{showcomments}}
{\textcolor{red}{(Anish says: \textit{#1}) }}{}}

\begin{document}



\title{\bf On the Inefficiency of Forward Markets \\ in Leader-Follower Competition}

\author{{Desmond Cai}\footnote{Dept. of Electrical Engineering, California Institute of Technology, Pasadena, CA 91125, Email:{wccai@caltech.edu}}, \qquad 
{Anish Agarwal}\footnote{Dept. of Computer Science, California Institute of Technology, Pasadena, CA 91125, Email:{aagarwal@caltech.edu}}, \qquad
{Adam Wierman}\footnote{Dept. of Computing and Mathematical Sciences Department, California Institute of Technology, Pasadena, CA 91125, Email:{adamw@caltech.edu}}
} 

\maketitle

\begin{abstract}%
Motivated by electricity markets, this paper studies the impact of forward contracting in situations where firms have capacity constraints and heterogeneous production lead times. We consider a model with two types of firms -- leaders and followers -- that choose production at two different times.  Followers choose productions in the second stage but can sell forward contracts in the first stage. Our main result is an explicit characterization of the equilibrium outcomes.  Classic results on forward contracting suggest that it can mitigate market power in simple settings; however the results in this paper show that the impact of forward markets in this setting is delicate -- forward contracting can enhance or mitigate market power.  In particular, our results show that leader-follower interactions created by heterogeneous production lead times may cause forward markets to be inefficient, even when there are a large number of followers.  In fact, symmetric equilibria do not necessarily exist due to differences in market power among the leaders and followers.
\end{abstract}%

\onehalfspacing

\section{Introduction}
\label{sec:intro}

Forward contracting plays a crucial role in a variety of markets, ranging from finance to cloud computing to commodities, e.g., gas and electricity.  Typically, forward contracting is viewed as a way to increase the efficiency of a marketplace. One way this happens is that forward contracts allow firms to hedge risks, e.g., risk from price fluctuations. In fact, the study of the efficiency created through hedging initiated the academic literature on forward contracting, e.g.,~\cite{Yamey1951,Roy1952,Telser1955,AndersonDanthine1980}. However, as the literature grew, other important benefits of forward markets emerged.

One of the most important of these additional benefits is the role that forward contracting plays in mitigating market power.  The seminal paper on this topic is~\cite{allaz1993cournot}, which studies a two-stage model of forward contracting in a setting where firms have perfect foresight (thus eliminating the gains possible via hedging risk).  This work showed, for the first time, that it is possible to mitigate market power using forward positions. Intuitively, this happens because the presence of forward contracting creates a situation where any individual firm has a strategic incentive to sell forward, which creates a prisoner's dilemma where, in equilibrium, all firms produce more than in the situation without forward contracts. Therefore, the market becomes more competitive.

Following the initial work of~\cite{allaz1993cournot}, the interaction of forward contracting and market power has been studied in depth.  The phenomenon has been investigated empirically in a thread of work that includes, e.g., ~\cite{Joskow2008,Green1999,GreenNewbery1992,TwomeyGreenNeuhoffNewbery2005,Wolfram1999,Wolak2000,BushnellMansurSaravia2008}. The robustness of the phenomenon to model assumptions has been studied in depth as well, e.g.,~\cite{gans1998contracts,newbery1998competition,green1999electricity,su2005equilibrium,liski2006forward,le2006forward,murphy2010impact}. The general consensus from this literature is that forward contracts often mitigate market power but that cases where the opposite occurs do exist.  At this point, characterizing when forward markets mitigate and when they enhance market power is still an open and active line of work.

In this paper, we contribute to the study of the interaction between forward contracting and market power by characterizing the impact of forward contracting in situations where firms have capacity constraints and heterogeneous production lead times.  Heterogeneous lead times create a leader-follower competition, where firms with long lead times (leaders) must decide production quantities well in advance and firms with short lead times (followers) can wait to decide production quantities while still participating in the early market via forward contracts.

Our study is motivated by the operation of electricity markets, where forward markets play a crucial role, capacity constraints are often binding, and generators have very different lead times.  The study of forward contracting in the context of electricity markets has received considerable attention, e.g.,~\cite{kamat2004two,yao2007two,oren2005generation,joskow2007,murphy2010impact}. This is motivated by the fact that as much as 95\% of electricity is traded through forward contracts, any time from minutes to months ahead of delivery.  These studies focus on whether forward contracting mitigates market power through inducing capacity investment or reducing network congestion. Market power has been a significant issue in electricity markets since their deregulation, e.g.,~\cite{BorensteinBushnell1999,jpmorgan}. The physical constraints of Kirchoff's laws and the non-storable nature of electricity have the potential to create hidden monopolies. Supply is also further constrained by generators' ramping limitations. Different generation resources could have significantly different ramping rates, ranging from as low as 1-7 MW/minute for oil and coal, to 50 MW/min for gas, and to more than 100 MW/min for hydro, e.g.,~\cite{WongAlbrechtAllanBillintonChenFongHaddadLiMukerjiPattonothers1999,Meehan2011,Shapiro2011}.


The same issues described above in the context of electricity markets are also prevalent in other marketplaces. For example, when new entrants to an industry such as gas or telecommunications must decide whether to invest in capacities, e.g., see ~\cite{de1997stochastic,demiguel2009stochastic}. The capacity expansion process is time consuming so the new entrants (leaders) must decide in advance the quantities they will supply to the market. The incumbents (followers) already have capacity and need only decide how much goods/services to provide. In both industries, forward contracts form a significant portion of the total output.

\subsection{Contributions of this paper}

This paper initiates the study of forward contracting in leader-follower competition.  While models exist, and are well studied \emph{independently}, for both forward contracting and leader-follower games; the \emph{combination} of the two has not been investigated previously.

In particular, this paper introduces a new model for studying the role of forward contracting in leader-follower Cournot competition. We consider a setting in which there are two types of firms -- \emph{leaders} and \emph{followers} -- that choose production levels at different stages subject to capacity constraints.  Leaders choose production levels before followers. However, followers are allowed to sell forward contracts when leaders are choosing their productions.

Our discussions above highlight that, due to the prisoner's dilemma effect, one may expect that allowing followers to trade forward contracts would increase their productions and, because forward contracts mitigate market power. The consequence of this would be that total production would also increase. In this work, we show that this intuition is not always true, and that the impact of forward contracting is ambiguous. The market power mitigation property of forward contracting might, in fact, create opportunities for leaders to exploit followers' capacity constraints to manipulate the market.

Specifically, the main results of this paper (Theorems~\ref{prop:MarketEquilibrium} and~\ref{prop:StackelbergEquilibrium} in Appendix~\ref{app:market}) provide detailed characterizations of the equilibrium productions with and without forward contracting, which give a complete picture of when forward contracting mitigates and when it enhances market power. As observed by~\cite{murphy2010impact}, capacity constraints may cause profit functions to be non-convex. Therefore, standard techniques used to show existence and uniqueness no longer apply. Nevertheless, we provide closed-form expressions of equilibria as a function of the parameters, including the number of leaders, number of followers, their marginal costs, and the capacity of the followers. Our explicit characterizations enable us to infer tradeoffs between the parameters as well as obtain the asymptotic behavior of the system as the numbers of leaders and followers increase. Among other properties, we show that there is an interval of follower productions just below capacity that are never symmetric equilibria, and that if there are too few leaders relative to followers, then there may not exist symmetric equilibria (Lemmas~\ref{lem:struct-market} and~\ref{lem:struct-market-3}). Moreover, we show that the efficiency loss due to forward contracting remains strictly positive even with a large number of followers (Lemma~\ref{lem:struct-compare}).

Our characterizations also show that market equilibria are especially interesting -- they may not exist or may not be unique -- at the transition between interior equilibria and full capacity utilization due to opportunities for market power exploitation. Thus, the result leads to a variety of qualitative insights about the interaction of forward markets and leader-follower competition.

First and foremost, our results highlight that forward contracting may decrease the efficiency of the market. The reason is that forward contracting may create opportunities for leaders to exploit the capacity constraints of the followers. Forward contracts provide incentives for followers to produce more. However, if this causes followers to become capacity constrained, then leaders would be able to profit by withholding their productions disproportionately, and the net effect is a decrease in total production. Informally, the increased competition due to forward contracting is offset by the decreased competition faced by the leaders due to followers being capacity constrained. Therefore, this is a phenomenon where capacity constraints and forward markets can create opportunities for market manipulation.

Second, and perhaps more damagingly, our results highlight that symmetric equilibrium may not exist as a result of market power exploits via capacity constraints. In particular, we show that symmetric equilibria do not exist precisely when followers are operating close to capacity. Our analysis shows that, when any follower operates close to capacity, other followers have a strategic incentive to exploit the fact that this follower is now \emph{less flexible} by reducing their forward positions. However, if all firms were to reduce their forward positions simultaneously, the high prices would create incentives for them to increase their forward positions. Therefore, there is no symmetric equilibria. This insight is related to the observation by~\cite{murphy2010impact} that equilibria may not exist. However, the argument in~\cite{murphy2010impact} was based on showing that profit functions are not convex and no explicit connection with strategic behavior or the circumstances under which equilibria do not exist were provided. On the other hand, we provide explicit conditions under which symmetric equilibria do not exist, and our analyses reveal the strategic interactions that precludes symmetric equilibrium.

\subsection{Related literature}

Our model, being a combination of the classical forward contracting and leader-follower models, has not been studied before. However, our study fits into the extensive (separate) literatures on each of forward contracting and leader-follower competition. In the following, we review the literature on these and explain how our work contributes to each of them.

\subsubsection*{Forward markets.}


\cite{allaz1993cournot} was the first to provide and analyze a model showing that strategic forward contracting mitigates market power. Later studies by~\cite{gans1998contracts,newbery1998competition,green1999electricity,su2005equilibrium,liski2006forward,le2006forward} reaffirmed or invalidated their findings under other assumptions. As these are not directly relevant to our work, we do not discuss their details here (see~\cite{murphy2010impact} for a survey). However, the general conclusion is that forward markets do not always mitigate market power in settings more general than that considered by \cite{allaz1993cournot}.

The domain of electricity markets has seen the most application of the model from~\cite{allaz1993cournot}. This may be attributed to the fact that the bulk of trade in electricity are through forward contracts and market power was a significant issue in wholesale electricity markets after their deregulation. However, capacity constraints is an important feature in electricity markets, and this feature is not present in their model. Therefore, there have been numerous extensions in this direction.~\cite{kamat2004two} and~\cite{yao2007two} added network constraints and price caps. However, due to the complexity of the problem, only numerical solutions were provided.~\cite{oren2005generation,joskow2007} proposed the idea that forward contracts may increase capacity investment. This idea was then investigated analytically by~\cite{murphy2010impact} by adding an endogenous capacity investment stage. The authors made the interesting finding that forward contracts may not mitigate market power when capacities are endogenous.

To our knowledge, our work is the first to study the robustness of the findings by~\cite{allaz1993cournot} in the classical leader-follower setting with capacity-constrained followers. In addition, our work supplements existing results on the impact of capacity constraints on existence of equilibria, by providing explicit characterizations under which symmetric equilibria exists and vice versa. This paper builds on our preliminary work, described in~\cite{Cai2013}, which illustrated that equilibria may not exist in the specific setting where leaders and followers have equal marginal costs. The work in the current paper considers a more general setting where leaders and followers could have different marginal costs, characterizes the asymptotic behavior of the system as the numbers of leaders and followers increase, and most importantly, characterizes the efficiency loss due to forward contracting.

\subsubsection*{Leader-follower competition.}

The first extension of Stackelberg's framework to multiple leaders and followers was provided by~\cite{sherali1983stackelberg,Sherali1984}. This work also gave conditions for existence and uniqueness of equilibria. Subsequently, there has been significant interest in relaxing the assumptions of the model. However, most studies focus on the technical conditions required for existence and uniqueness, neglecting to study the underlying strategic behavior.~\cite{ehrenmann2004manifolds} showed that equilibrium is no longer unique if one removes Sherali's assumption that identical producers make identical decisions.~\cite{de1997stochastic,xu2005mpcc,demiguel2009stochastic} generalized some of Sherali's existence and uniqueness results to the setting with uncertainty. There are also other efforts by~\cite{pang2005quasi,kulkarni2012global} that provide conditions for existence using variational inequality techniques.

We are not aware of any work that adds capacity constraints to Sherali's model. The closest related work was by~\cite{nie2012duopoly} but the authors were investigating price competition (while we focus on quantity competition).~\cite{demiguel2009stochastic} might appear to have included capacity limits in their analyses. However, the authors used the capacity limits as a technical condition for their proof, since it was defined by the point where marginal cost exceeds price. Therefore, their capacity constraints are never binding, and firms in their model do not strategically withhold productions (unlike in our model).

Our work is the first to extend Sherali's model with capacity constraints on followers while allowing them to sell forward contracts. Similar to Sherali's work, we restrict ourselves to symmetric equilibria in the sense that leaders have equal productions and followers have equal forward positions. We characterize all symmetric equilibria and provide insights into strategic behavior. Note that~\cite{ehrenmann2004manifolds} showed that equilibrium is no longer unique if Sherali's symmetry assumptions are relaxed. However, his findings are technically different from ours. His results are attributed to non-smoothness due to the non-negativity constraints on quantities, while our results are attributed to non-smoothness due to the capacity constraints. Therefore, a symmetric equilibria always exists in~\cite{ehrenmann2004manifolds} but may not exist in our model.

\section{Model}
\label{sec:model}

Our goal is to understand whether forward contracting mitigates market power when firms have capacity constraints and heterogeneous production lead times. To this end, we formulate a model that combines key elements from the classical forward contracting model proposed in~\cite{allaz1993cournot} as well as the classical leader-follower model proposed in~\cite{sherali1983stackelberg}.

We assume that there are two types of firms -- \emph{leaders} and \emph{followers} -- that choose production quantities at two different times. Leaders, who have longer lead times than followers, choose production quantities in the first stage, while followers choose production quantities in the second stage. However, followers sell forward contracts in the first stage. Thus, we also refer to the first stage as the \emph{forward market} and the second stage as the \emph{spot market}.

\subsection{Forward contracting}
\label{subsec:model-market}

Our model for forward contracting is based on the classical model from~\cite{allaz1993cournot}. This model is commonly used in many studies of forward markets~\cite{newbery1998competition,green1999electricity,kamat2004two,liski2006forward,le2006forward,yao2007two,murphy2010impact}. In the forward market, firms sign contracts to deliver a certain quantity of good at a price $\pf$. These contracts are binding and observable pre-commitments. Then, in the spot market, firms sell the good at a price $P(q)$ which is a function of the total quantity $q$ of the good sold in both the forward and spot markets. We assume a linear demand model given by
\begin{align*}
P(q) = \rv - \beta q,
\end{align*}
where the constants $\rv, \beta > 0$. This is a common model for demand~\cite{allaz1993cournot,murphy2010impact} and implies that buyers' aggregate utility is quasilinear in money and quadratic in the quantity of the good consumed.

We assume that there is perfect foresight. That is, in the first stage, both leaders and followers know the demand in the second stage. Equilibrium then requires that the forward and spot prices are aligned:
\begin{align*}
\pf = P(q).
\end{align*}
That is, no arbitrage is possible. This assumption was also used in both the classical forward contracting model~\cite{allaz1993cournot} and the classical leader-follower model~\cite{sherali1983stackelberg}. An extension to the case of uncertain demand is definitely relevant and interesting. But our results show that the model with certain demand is rich enough to capture interesting strategic interactions between leaders and followers. The case of uncertain demand is left to future work.

\subsection{Production lead times}
\label{subsec:model-gen}

Our model for leader-follower competition is based on the classical model from~\cite{sherali1983stackelberg}. We assume that there are $\nb$ leaders and $\na$ followers, with marginal costs $\cb$ and $\ca$, respectively, where $\ca \geq \cb > 0$. We also abuse notation and use $\nb$ and $\na$ to denote the set of leaders and followers, respectively. The assumption that $\ca \geq \cb$ is motivated by the expectation that there is typically a cost to flexibility, e.g. in electricity markets more flexible generators typically have higher operating costs than less flexible generators.

Each leader $i\in \nb$ chooses its production quantity $\qbpj$ in the forward market. Each follower $j \in \na$ chooses its production quantity $\qapi$ in the spot market and also sells a forward contract of quantity $\qfpi$ in the forward market. We assume that leaders sell forward contracts in the forward market equal to their committed productions. It is possible to show that allowing leaders to sell forward contracts that differ from their committed productions does not change the analyses.

We assume that each follower has a production capacity $\ka > 0$ but leaders are not capacity constrained. In practice, followers might only be able to adjust productions within a limited range around operating points. Thus, a more sophisticated model would have followers choose set points in the forward market and impose constraints on deviations from those set points. Our model for followers can be interpreted as them having zero set points and being allowed to ramp productions to a maximum of $\ka$. Similarly, our model for leaders can be interpreted as them choosing their operating points in the forward market and not being allowed to deviate from them.

\subsection{Competitive model}
\label{subsec:model-equil}

We adopt the following equilibria concept for the market. Let the vectors $\qbp = \left(x_{1},\ldots,x_{\nb}\right)$, $\qap = \left(y_{1},\ldots,y_{\na}\right)$, and $\qfp = \left(f_{1},\ldots,f_{\na}\right)$ denote the leaders' productions, followers' productions, and followers' forward contracts, respectively. We also use the notation $\qfpn=\left(f_{1},\ldots,f_{j-1},f_{j+1},\ldots,f_{\na}\right)$ to denote the forward contracts of all followers other than $i$. Similarly, we use the notations $\qbpn=\left(x_{1},\ldots,x_{i-1},x_{i+1},\ldots,x_{\nb}\right)$ and $\qapn=\left(y_{1},\ldots,y_{j-1},y_{j+1},\ldots,y_{\na}\right)$.

\emph{Spot market (followers):}  We define the spot market equilibrium as follows. Only followers compete in the spot market. Follower $j$'s profit from the spot market is:
\begin{align*}
\pfpis\left(\qapi;\qapn\right)
&=
P\left(\sum_{i'=1}^{\nb}\qbpjj + \sum_{j'=1}^{\na}\qapii\right)\cdot \left(\qapi - \qfpi\right) - \ca\qapi.
\end{align*}
Given $\qapn$, follower $j$ chooses a production $\qapi$ to maximize its profit subject to its capacity constraint. Thus, a Nash equilibrium of the spot market is a vector $\qap$ such that for all $j$:
\begin{align*}
\pfpis\left(\qapi;\qapn\right) \geq \pfpis\left(\bar{y}_{j};\qapn\right),\;\mbox{for all }\bar{y}_{j}\in\left[0,\ka\right].
\end{align*}
Theorem 5 of~\cite{jing1999spatial} implies that there always exists a unique spot equilibrium given any leader productions and follower forward positions $(\qbp,\qfp)$. We denote this unique equilibrium by $\qapopt\left(\qfp,\qbp\right)=\left(y_{1}\left(\qfp,\qbp\right),\ldots,y_{\na}\left(\qfp,\qbp\right)\right)$. 

\emph{Forward market:}  The forward market equilibrium depends on behaviors of both followers and leaders. Their profits depend on the outcome of the spot market. In particular, follower $j$'s profit is given by:
\begin{align*}
\pfpif\left(\qfpi;\qfpn,\qbp\right)
&=
P\left(\sum_{i'=1}^{\nb}\qbpjj + \sum_{j'=1}^{\na}\qapiiopt(\qfp,\qbp)\right)\cdot \qfpi + \pfpis\left(\qapopt(\qfp,\qbp)\right)
\\
&=
\left(P\left(\sum_{i'=1}^{\nb}\qbpjj + \sum_{j'=1}^{\na}\qapiiopt(\qfp,\qbp)\right)-\ca\right)\cdot \qapiopt(\qfp,\qbp),
\end{align*}
where the second equality follows by substituting for $\pfpis(\qap(\qfp,\qbp))$. Note that follower $j$ anticipates the impact of the actions in the forward market on the spot market. Given $(\qfpn,\qbp)$, follower $j$ chooses its forward contract $\qfpi$ to maximize its profit. This is an unconstrained maximization as followers can take positive or negative positions in the forward market. Next, leader $i$'s profit is given by:
\begin{align*}
\pfpjf\left(\qbpj;\qbpn,\qfp\right)
&=
\left(P\left(\sum_{i'=1}^{\nb}\qbpjj + \sum_{j'=1}^{\na}\qapii(\qfp,\qbp)\right)-\cb\right)\cdot\qbpj.
\end{align*}
Given $(\qbpn,\qfp)$, leader $i$ chooses a production $\qbpj \in\mathbb{R}_+$ to maximize its profit.

Thus, a subgame perfect Nash equilibrium of the forward market is a tuple $\left(\qfp,\qbp\right)$ such that for all $i$:
\begin{align}
\label{eq:nash-equilibrium-leader}
\pfpjf\left(\qbpj;\qbpn,\qfp\right)
\geq
\pfpjf\left(\bar{x}_{i};\qbpn,\qfp\right),
\;
\mbox{for all }
\bar{x}_{i}\in\mathbb{R}_+,
\end{align}
and for all $j$:
\begin{align}
\label{eq:nash-equilibrium-follower}
\pfpif\left(\qfpi;\qfpn,\qbp\right)
\geq
\pfpif\left(\bar{f}_{j};\qfpn,\qbp\right),
\;
\mbox{for all }
\bar{f}_{j}\in\mathbb{R}.
\end{align}
It is this equilibrium that is the focus of this study. To capture the key strategic interactions between leaders and followers, we focus on equilibria in which leaders have symmetric productions and followers have symmetric forward positions. This symmetric case already offers many insights. 

\section{One Leader and Two Followers}
\label{sec:SimpleExample}

The complete characterization of the model is technical. So, we defer that analysis and discussion to Section \ref{sec:StructuralResults} and begin by developing intuition for our results in a special case.  Specifically, we start by considering only $\nb=1$ leader, $\na=2$ followers, and equal marginal costs $\cb = \ca$. This case, though simple, is already rich enough to expose the structure of the general results and to highlight the inefficiencies that arise from forward contracting in leader-follower competition.

The section is organized as follows.  First, in Sections~\ref{sec:SimpleExample:PeakerReaction} and~\ref{sec:SimpleExample:BaseloadReaction}, we study the reactions of the followers to the leader and vice versa, respectively. In particular, we focus on the impact of followers' capacity constraints and leader's commitment power on their responses to the other producers' actions. Then, in Section~\ref{sec:SimpleExample:ForwardMarketEquilibrium}, we study how they impact the equilibria of the market. Finally, in Sections~\ref{sec:SimpleExample:ForwardMarket}, we study how followers' forward contracting impact market outcomes.

Throughout this section, we denote the normalized demand by:
\begin{align*}
    \bar{\rv} := \frac{1}{\beta}(\rv-\cb).
\end{align*}
Recall that $\rv$ is the maximum price that demand is willing to pay and $\cb$ is the minimum price that producers need to receive for them to supply to the market. Therefore, we restrict our analyses to the case where $\bar{\rv} \geq 0$.

\subsection{Follower reaction}
\label{sec:SimpleExample:PeakerReaction}
We begin by studying how followers respond when the leader produces a fixed quantity $\qb\in\mathbb{R}_+$. We focus on symmetric responses, that is, those where followers take equal forward positions. Let $\Qf :\mathbb{R}_+ \rightarrow \mathbb{P}(\mathbb{R})$ denote the symmetric reaction correspondence of the followers, i.e., for each $f\in \Qf(\qb)$,
\begin{align*}
& \phi_1(\qf;\qf,\qb) \geq \phi_1(\bar{f};\qf,\qb), \quad \forall \bar{f}\in\mathbb{R};
\\
\text{and} \quad & \phi_2(\qf;\qf,\qb) \geq \phi_2(\bar{f};\qf,\qb), \quad \forall \bar{f}\in\mathbb{R}.
\end{align*}
Proposition~\ref{prop:SpotEquilibrium2} in the Appendix implies that the followers produce equal quantities $y_1(\qf;\qf,\qb)=y_2(\qf;\qf,\qb)$. Let $\Qa:\mathbb{R}_+ \rightarrow \mathbb{P}(\mathbb{R}_+)$ denote the production correspondence of the followers, i.e., for each $y\in\Qa(\qb)$, there exists $\qf\in\Qf(\qb)$ such that $y_1(\qf;\qf,\qb) = y_2(\qf;\qf,\qb) = y$. Applying Propositions~\ref{prop:SpotEquilibrium2} and~\ref{prop:PeakerReaction} in the Appendix, the reaction and production correspondences are given by:
\begin{align*}
    \begin{array}{llll}
& \Qf(x) = [-\bar{\rv}+\qb+3\ka,\infty), & \Qa(x) = \{\ka\}, & \mathrm{if} \; \qb \leq \bar{\rv}-3\ka,
\\
& \Qf(x) = \varnothing, & \Qa(x) = \varnothing, & \mathrm{if} \; \bar{\rv} - 3\ka < \qb < \bar{\rv}-\frac{5}{5-2\sqrt{2}}\ka,
\\
& \Qf(x) = \left\{ \frac{1}{5}(\bar{\rv}-\qb) \right\}, & \Qa(x) = \left\{ \frac{2}{5}(\bar{\rv}-x) \right\}, & \mathrm{if} \; \bar{\rv} - \frac{5}{5-2\sqrt{2}}\ka \leq \qb \leq \bar{\rv},
\\
& \Qf(x) = (-\infty,-\bar{\rv}+\qb], & \Qa(x) = \{0\}, & \mathrm{if} \; \bar{\rv} \leq \qb.
    \end{array}
\end{align*}

\begin{figure}[t]
\begin{center}
\includegraphics[height=2in]{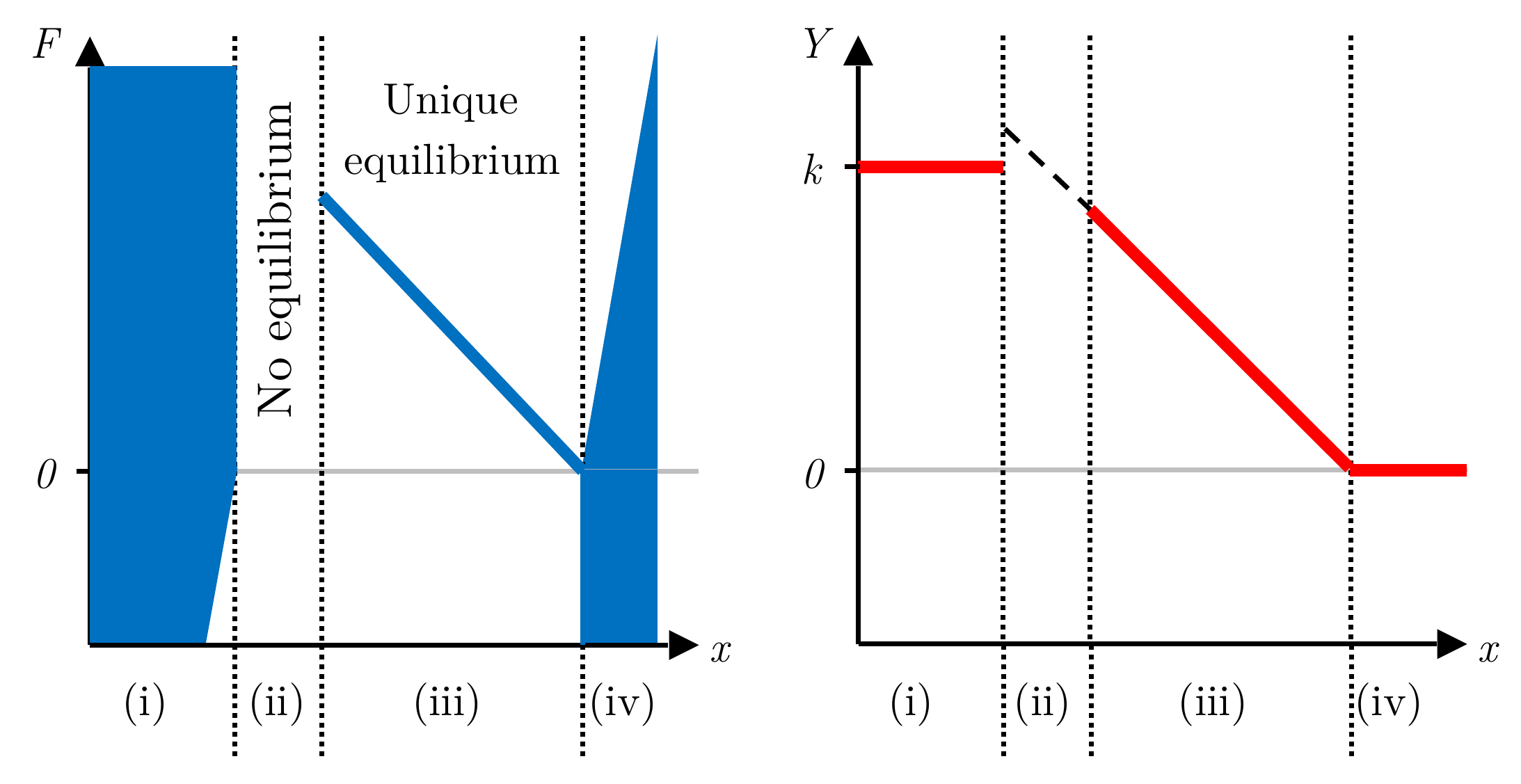}
\caption{Follower reaction correspondence $\Qf$ and production correspondence $\Qa$.}
\label{fig:PeakFwd}
\end{center}
\end{figure}

Figure~\ref{fig:PeakFwd} shows the characteristic shapes of $\Qf$ and $\Qa$. There are four major segments labelled (i)~--~(iv). Note that the follower productions are always $\ka$ in segment (i) and $0$ in segment (iv). In general, one expects followers' reactions to decrease as $\qb$ increases because a higher leader production decreases the demand in the spot market. This behavior indeed holds in segment (iii), which is also the behavior in a conventional forward market in the absence of capacity constraints. However, the capacity constraints lead to complex reactions, as seen in segments (i), (ii), and (iv).

\emph{Segments (i) and (iv): $\qb \leq \bar{\rv}-3\ka$ or $\bar{\rv} \leq \qb$. Multiple equilibria.} These are degenerate scenarios where followers have binding productions, and hence are neutral to a range of different forward positions, as they all lead to the same production outcomes. The structure of the reaction set is also intuitive. Consider segment (i), where followers produce zero quantities. If $\qf'$ is a symmetric reaction, then any $\qf''<\qf'$ is also a symmetric reaction, since decreasing forward positions create incentives to decrease productions, and productions cannot drop below zero. Therefore, the reaction sets are left half-lines. A similar argument applies to segment (iv), but in this case, the reaction sets are right half-lines.

\emph{Segment (ii): $\bar{\rv}-\frac{5}{5-2\sqrt{2}}\ka\leq\qb\leq\bar{\rv}$. No equilibrium.} This is the scenario where followers' capacity constraints create incentives for market manipulation which causes symmetric reactions to disappear. The type of symmetric reactions in segment (iii) are unsustainable here because each follower has incentive to reduce its forward position. For example, when follower $1$ reduces its forward position, it induces follower $2$ to increase its production. However, since follower $2$ can only increase its production up to $\ka$, the total production decreases, the market price increases, and follower $1$'s profit increases. By symmetry, follower $2$ has incentive to manipulate the market in a similar manner. Yet, should both followers reduce their forward positions, there will be excess demand in the market. Therefore, there is no symmetric equilibrium between the followers.

%

\subsection{Leader reaction}
\label{sec:SimpleExample:BaseloadReaction}
Next, we study how the leader responds when both followers take a fixed forward position $\qf\in\mathbb{R}$. Let $\Qb:\mathbb{R}\rightarrow\mathbb{P}(\mathbb{R}_+)$ denote the leader's reaction correspondence, i.e., for each $\qb\in\Qb(\qf)$,
\begin{align*}
\psi_1 \left(\qb;\qf,\qf\right)
\geq
\psi_1\left(\bar{x};\qf,\qf\right),
\quad \forall
\bar{x}\in\mathbb{R}_+.
\end{align*}
Let $\Qa:\mathbb{R}\rightarrow\mathbb{P}(\mathbb{R}_+)$ denote the production correspondence of the followers, i.e., for each $y\in\Qa(\qf)$, there exists $\qb\in\Qb(\qf)$ such that $y_1(\qf;\qf,\qb)=y_2(\qf;\qf,\qb)=y$. The expressions for $\Qb$ and $\Qa$ can be obtained from Propositions~\ref{prop:BaseloadReaction} and~\ref{prop:SpotEquilibrium2} in the Appendix. $\Qb$ and $\Qa$ takes three distinctive shapes depending on the value of $\bar{\rv}$.

\begin{figure}[t]
      \centering
      \begin{subfigure}[b]{0.3\textwidth}
      \includegraphics[width=\textwidth]{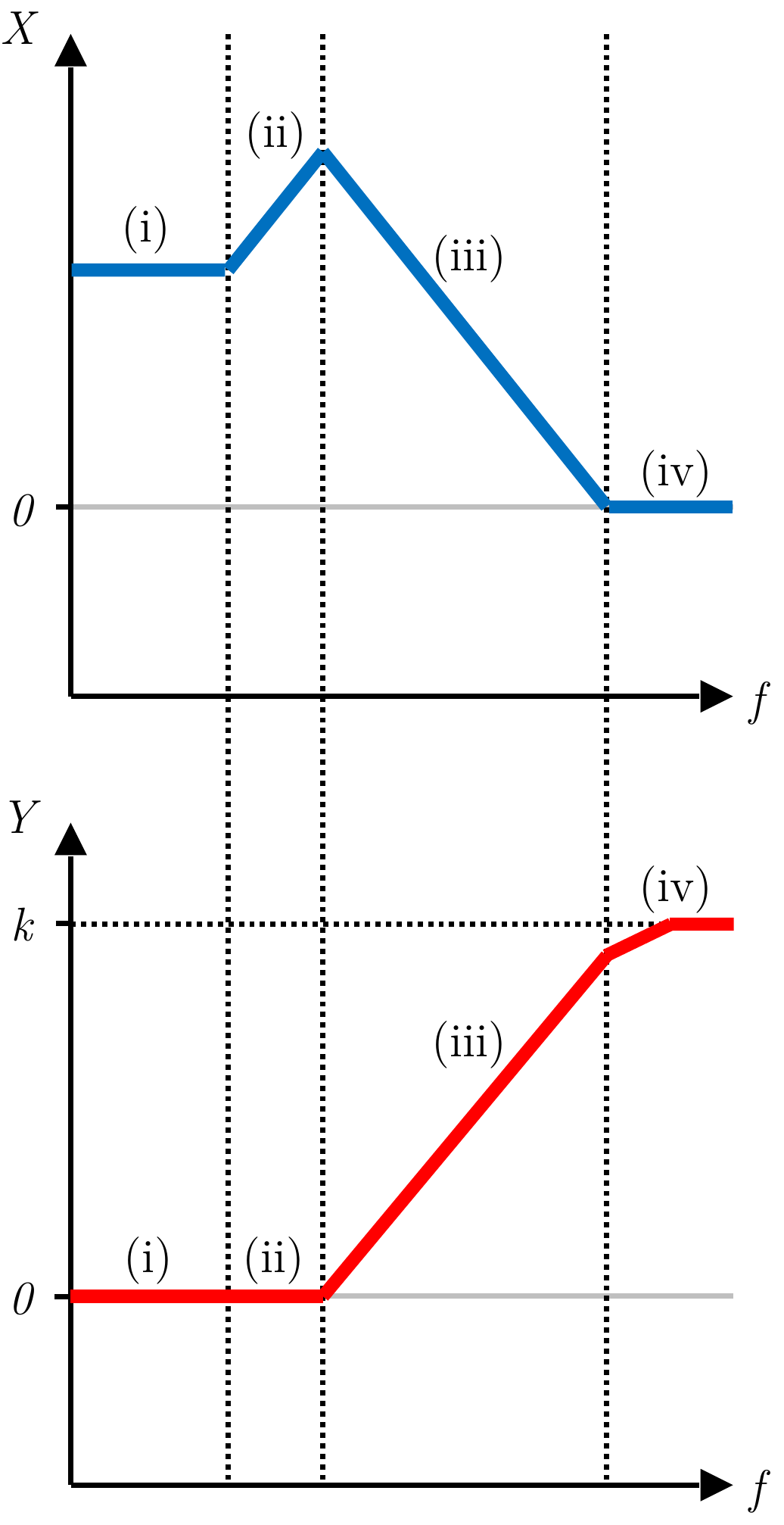}
      \caption{$0 \leq \bar{\rv} \leq 2\ka$}
      \label{fig:BaseFwdA}
      \end{subfigure}
      \begin{subfigure}[b]{0.3\textwidth}
      \includegraphics[width=\textwidth]{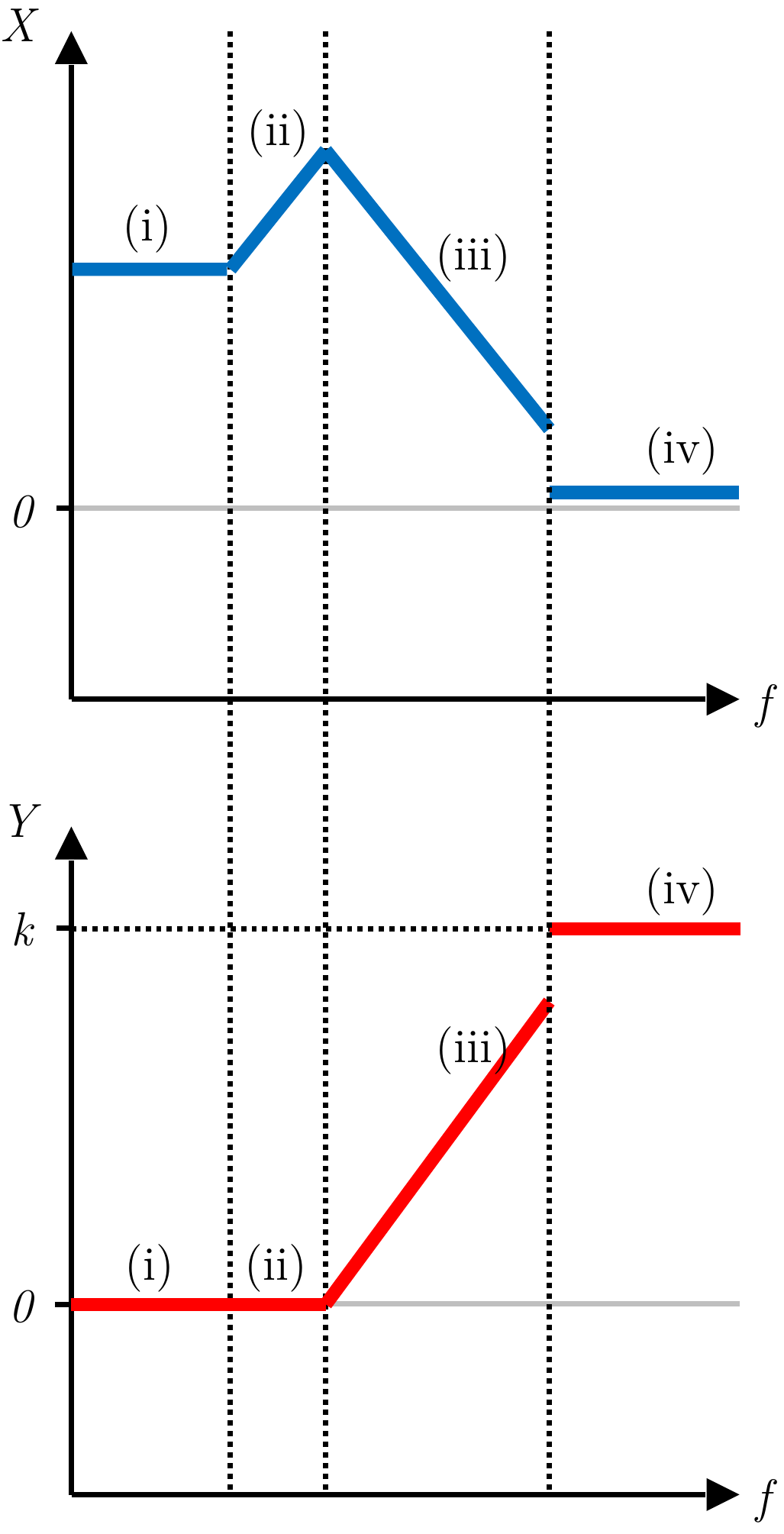}
      \caption{$2\ka < \bar{\rv} < \frac{4}{2-\sqrt{3}}\ka$}
      \label{fig:BaseFwdB}
      \end{subfigure}
      \begin{subfigure}[b]{0.3\textwidth}
      \includegraphics[width=\textwidth]{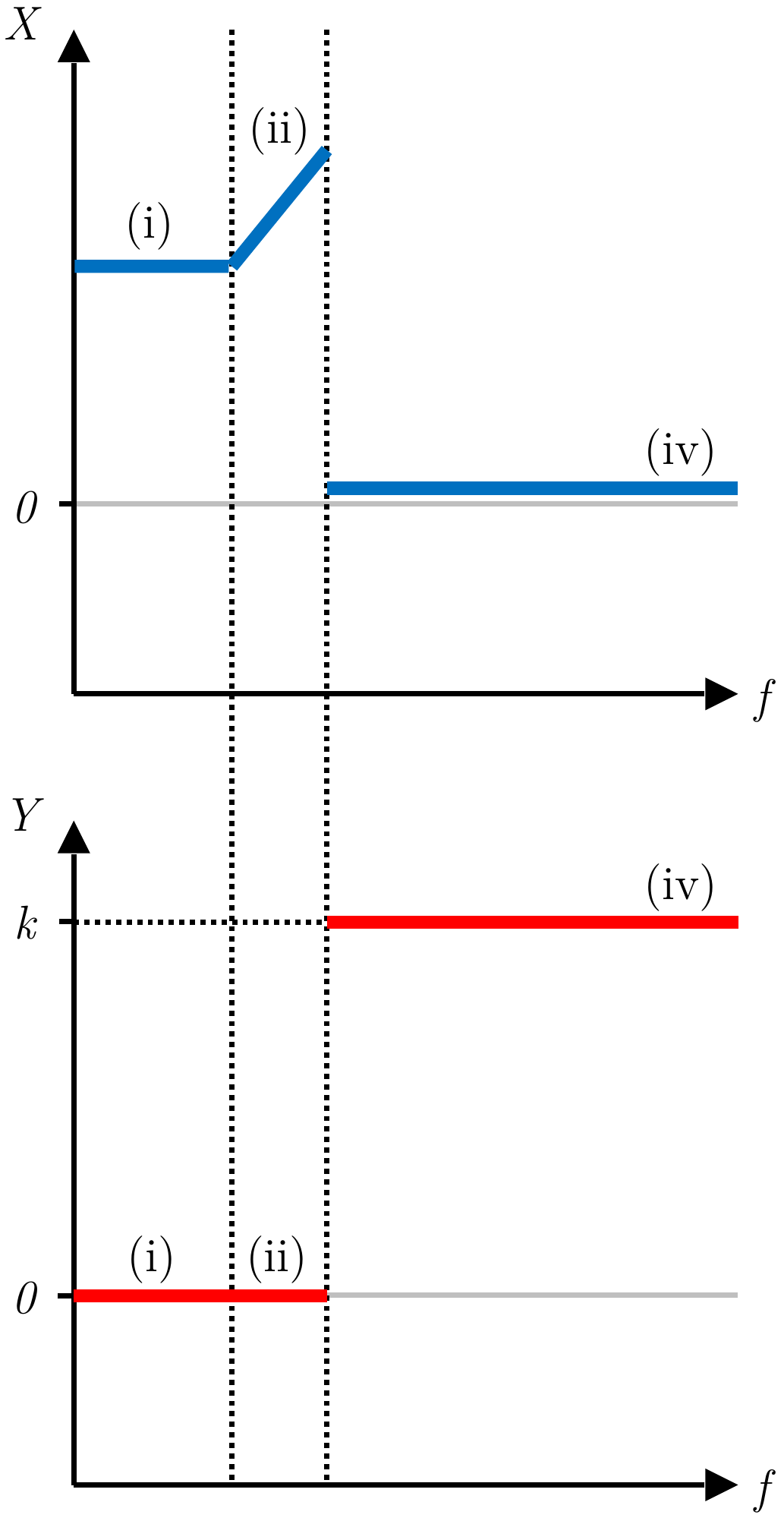}
      \caption{$\frac{4}{2-\sqrt{3}}\ka \leq \bar{\rv}$}
      \label{fig:BaseFwdC}
      \end{subfigure}
      \caption{Leader reaction correspondence $\Qb$ and follower production correspondence $\Qa$.}
      \label{fig:BaseFwd}
\end{figure}

\emph{Low demand: $0 \leq \bar{\rv} \leq 2\ka$.} In this case, the reaction and production correspondences are given by
\begin{align*}
    \begin{array}[]{llll}
        & \Qb(\qf) = \{\frac{1}{2}\bar{\rv}\}, & \Qa(\qf) = \{0\}, & \mathrm{if} \; \qf \leq -\frac{1}{2}\bar{\rv},
        \\
        & \Qb(\qf) = \{\bar{\rv}+\qf\}, & \Qa(\qf) = \{0\}, & \mathrm{if} \; -\frac{1}{2}\bar{\rv} \leq \qf \leq -\frac{1}{4}\bar{\rv},
        \\
        & \Qb(\qf) = \{\frac{1}{2}\bar{\rv}-\qf\}, & \Qa(\qf) = \{\frac{1}{6}\bar{\rv}+\frac{2}{3}\qf\}, & \mathrm{if} \; -\frac{1}{4}\bar{\rv} \leq \qf \leq \frac{1}{2}\bar{\rv},
        \\
        & \Qb(\qf) = \{0\}, & \Qa(\qf) = \{\frac{1}{3}(\bar{\rv}+\qf), & \mathrm{if} \; -\frac{1}{2}\bar{\rv} \leq \qf \leq 3\ka-\bar{\rv},
        \\
        & \Qb(\qf) = \{0\}, & \Qa(\qf) = \{\ka\}, & \mathrm{if} \; 3\ka - \bar{\rv} \leq \qf.
    \end{array}
\end{align*}
Figure~\ref{fig:BaseFwdA} shows the characteristic shapes of $\Qb$ and $\Qa$. There are four major segments labelled (i)~--~(iv). The follower supplies $0$ in segments (i) and (ii) and supplies $\ka$ for a subset of segment (iv). In general, one expects the leader's production to decrease as $\qf$ increases, because larger forward positions lead to larger follower supplies, which decreases the market price. This behavior indeed holds in segment (iii). However, the capacity constraints and leader's commitment power lead to complex reactions in segments (i), (ii), and (iv).

\emph{Segment (i) and (iv): $\qf \leq -\frac{1}{2}\bar{\rv}$ or $3\ka-\bar{\rv}\leq\qf$. Constant production.} These are degenerate scenarios where the leader is insensitive to the followers' forward positions. When $\qf \leq -\frac{1}{2}\bar{\rv}$, it is because followers' always supply zero regardless of their forward positions. When $3\ka-\bar{\rv} \leq \qf$, it is because followers supply large quantities, and drive prices down below the level at which it is profitable for leaders to produce.

\emph{Segment (ii): $-\frac{1}{2}\bar{\rv} \leq \qf \leq -\frac{1}{4}\bar{\rv}$. Increasing production.} In this scenario, the leader uses its commitment power to drive the followers out of the market. As followers increase their forward positions, the leader, instead of decreasing its production (as one typically expects), actually increases its production, as doing so allows it to depress demand below the level at which followers are willing to supply.

\emph{Medium demand: $2\ka < \bar{\rv} < \frac{4}{2-\sqrt{3}}\ka$.} In this scenario, the reaction and production correspondences are given by
\begin{align*}
    \begin{array}[]{llll}
        & \Qb(\qf)=\{\frac{1}{2}\bar{\rv}\}, & \Qa(\qf)=\{0\}, & \mathrm{if}\;\qf \leq -\frac{1}{2}\bar{\rv},
        \\
        & \Qb(\qf)=\{\bar{\rv}+\qf\}, & \Qa(\qf)=\{0\}, & \mathrm{if}\;-\frac{1}{2}\bar{\rv} \leq \qf \leq -\frac{1}{4}\bar{\rv},
        \\
        & \Qb(\qf)=\{-\frac{1}{2}\bar{\rv}\}, & \Qa(\qf)=\{\frac{1}{6}\bar{\rv}+\frac{2}{3}\qf\}, & \mathrm{if}\;-\frac{1}{4}\bar{\rv} \leq \qf \leq -\frac{\sqrt{3}-1}{2}\bar{\rv}+\sqrt{3}\ka,
        \\
        & \Qb(\qf)=\{\frac{1}{2}\bar{\rv}-\qf,\frac{1}{2}\bar{\rv}-\ka\}, & \Qa(\qf)=\{\frac{1}{6}\bar{\rv}+\frac{2}{3}\qf,\ka\}, & \mathrm{if}\;\qf=-\frac{\sqrt{3}-1}{2}\bar{\rv}+\sqrt{3}\ka,
        \\
        & \Qb(\qf)=\{\frac{1}{2}\bar{\rv}-\ka\}, & \Qa(\qf)=\{\ka\}, & \mathrm{if}\;-\frac{\sqrt{3}-1}{2}\bar{\rv}+\sqrt{3}\ka < \qf.
    \end{array}
\end{align*}
Figure~\ref{fig:BaseFwdB} shows the characteristic shapes of $\Qb$ and $\Qa$. There are, again, four major segments labelled (i)~--~(iv). Segments (i), (ii), and (iii), are similar to that in the low demand case when $0 \leq \bar{\rv} \leq 2\ka$. Segment (iv), however, is different in that, while the leader produces zero in this segment when $0 \leq \bar{\rv} \leq 2\ka$, the leader now produces a strictly positive quantity in this segment. This is due to the fact that the leader's profit on each unit is given by
\begin{align*}
    \rv - \beta(y_1+y_2+\qb) - \cb = \beta(\bar{\rv}-y_1-y_2-\qb),
\end{align*}
and hence, when $\bar{\rv} > 2\ka$, the leader is still able to profit from producing when both followers produce $\ka$. Due to this, the leader also has an incentive to exploit followers' capacity constraints, unlike previously when the leader was producing zero. The leader does so by sharply reducing its production at the end of segment (iii). This induces the followers to increase their supply, but since followers can only increase their supply up to $\ka$, the total market production decreases, the market price increases, and the leader's profit increases. Therefore, there is a discontinuity in the leader's reaction curve between segments (iii) and (iv).

\emph{High demand: $\frac{4}{2-\sqrt{3}}\ka \leq \bar{\rv}$.} In this scenario, the reaction and production correspondences are given by
\begin{align*}
    \begin{array}[]{llll}
        & \Qb(\qf)=\{\frac{1}{2}\bar{\rv}\}, & \Qa(\qf)=\{0\}, & \mathrm{if}\;\qf \leq -\frac{1}{2}\bar{\rv},
        \\
        & \Qb(\qf)=\{\bar{\rv}+\qf\}, & \Qa(\qf)=\{0\}, & \mathrm{if}\;-\frac{1}{2}\bar{\rv} \leq \qf < -\frac{1}{2}\bar{\rv}+\sqrt{(\bar{\rv}-\ka)\ka},
        \\
        & \Qb(\qf)=\{\bar{\rv}+\qf,\frac{1}{2}\bar{\rv}-\ka\}, & \Qa(\qf)=\{0,\ka\}, & \mathrm{if}\;\qf=-\frac{1}{2}\bar{\rv}+\sqrt{(\bar{\rv}-\ka)\ka},
        \\
        & \Qb(\qf)=\{\frac{1}{2}\bar{\rv}-\ka\}, & \Qa(\qf)=\{\ka\}, & \mathrm{if}\;-\frac{1}{2}\bar{\rv}+\sqrt{(\bar{\rv}-\ka)\ka} < \qf.
    \end{array}
\end{align*}
Figure~\ref{fig:BaseFwdC} shows the characteristic shapes of $\Qb$ and $\Qa$. There are now only three segments labelled (i), (ii), and (iv). These segments are similar to segments (i), (ii), and (iv), respectively, in the medium demand case where $2\ka < \bar{\rv} < \frac{4}{2-\sqrt{3}}\ka$. The difference is that, now, the leader decreases its production sharply once followers begin to supply to the market. As a consequence, the followers' supply jumps from $0$ to $\ka$.

\subsection{Forward market equilibrium}
\label{sec:SimpleExample:ForwardMarketEquilibrium}

We now study the equilibria of the forward market. Let $Q\subseteq\mathbb{R}\times\mathbb{R}_+$ denote the set of all symmetric equilibria and $\Qa\subseteq\mathbb{R}_+$ denote the set of all follower productions, i.e., $(\qf,\qb)\in Q$ if $(\qf,\qf,\qb)$ is a Nash equilibrium of the forward market, and $y\in\Qa$ if there exists $(\qf,\qb)\in Q$ such that $y_1(\qf;\qf,\qb)=y_2(\qf;\qf,\qb)=y$. From Theorem~\ref{prop:MarketEquilibrium} and Proposition~\ref{prop:SpotEquilibrium2}, the symmetric equilibria and follower productions are given by
\begin{align*}
    \begin{array}[]{llll}
        & Q=\{(\qf,\qb):\;\qf=\frac{1}{8}\bar{\rv},\;\qb=\frac{3}{8}\bar{\rv}\}, & \Qa=\{\frac{1}{4}\bar{\rv}\}, & \mathrm{if} \; 0 \leq \bar{\rv} \leq \frac{8}{4-\sqrt{3}}\ka,
        \\
        & Q=\varnothing, & \Qa=\varnothing, & \mathrm{if} \; \frac{8}{4-\sqrt{3}}\ka < \bar{\rv} < 4\ka,
        \\
        & Q=\{(\qf,\qb):\;\qf\in[-\frac{1}{2}\bar{\rv}+2\ka,\infty),\;\qb=\frac{1}{2}\bar{\rv}-\ka\}, & \Qa=\{\ka\}, & \mathrm{if} \; 4\ka \leq \bar{\rv}.
    \end{array}
\end{align*}
Observe that there are three operating regimes.

\emph{Low demand: $0 \leq \bar{\rv} \leq \frac{8}{4-\sqrt{3}}\ka$.} There is a one symmetric equilibrium. Productions increase as $\bar{\rv}$ increases. This regime is identical to that in the absence of capacity constraints (to see this, substitute $\ka=\infty$).

\emph{Medium demand: $\frac{8}{4-\sqrt{3}}\ka < \bar{\rv} < 4\ka$.} There is no symmetric equilibrium. This phenomena is due to leaders and followers withholding productions and forward contracts respectively. As observed in the separate reaction curves, each individual follower or leader has incentive to exploit the capacity constraints of the followers by reducing its position in the forward market. But should all producers do so, there will be excess demand in the market, and hence no symmetric equilibria are sustainable.

\emph{High demand: $4\ka \leq \bar{\rv}$.} There is a unique equilibrium leader production $\frac{1}{2}\bar{\rv}-\ka$ and infinitely many equilibrium follower forward positions $[-\frac{1}{2}\bar{\rv}+2\ka,\infty)$. The latter is a right half-line because followers are supplying all their capacity and so are indifferent once forward positions exceed a certain value. The leader production increases with demand; although the rate of increase of $\frac{1}{2}$ is slower than in the case when demand is low, where it increased at the rate $\frac{3}{8}$. This distinction is due to the leader facing less competition than before when followers were not capacity-constrained.

    Note that, unlike with the leader's reaction curve, there is no apparent phenomenon where the leader increases its production to drive followers out of the market. This can be attributed to the fact that the leader and followers have equal marginal costs. In Section~\ref{sec:StructuralResults}, we will see that the leader's commitment power does cause its equilibrium production to increase with demand, when we relax the assumption of equal marginal costs.

\subsection{Inefficiency of the forward market}
\label{sec:SimpleExample:ForwardMarket}

To study the efficiency of the forward market, we compare the outcome in our market against that in a Stackelberg competition, where followers do not sell forward contracts. Therefore, the leader continues to commit to its production ahead of the followers.

Note that the symmetric Stackelberg equilibria are simply the symmetric reactions of the leader when followers take neutral forward positions. Therefore, using the notation in Section~\ref{sec:SimpleExample:BaseloadReaction}, we let $\Qb(0)$ denote the set of all symmetric Stackelberg equilibria, i.e., for each $\qb\in\Qb(0)$,
\begin{align*}
    & \psi_1(\qb;0,0) \geq \psi_1(\bar{\qb};0,0), \quad \forall\bar{\qb}\in\mathbb{R}_+,
\end{align*}
and we let $\Qa\subseteq\mathbb{R}_+$ denote the set of all follower productions, i.e., $\qa\in\Qa$ if there exists $\qb\in\Qb(0)$ such that $\qa_1(0;0,\qb)=\qa_2(0;0,\qb)=\qa$. From Theorem~\ref{prop:StackelbergEquilibrium}, the Stackelberg equilibria are given by
\begin{align*}
    \begin{array}[]{llll}
        & \Qb(0) = \left\{ \frac{1}{2}\bar{\rv} \right\}, & \Qa = \left\{\frac{1}{6}\bar{\rv}\right\}, & \mathrm{if} \; 0 \leq \bar{\rv} < \frac{2\sqrt{3}}{\sqrt{3}-1}\ka,
        \\
        & \Qb(0) = \left\{ \frac{1}{2}\bar{\rv},\frac{1}{2}\bar{\rv}-\ka \right\}, & \Qa = \left\{\frac{1}{3-\sqrt{3}}\ka, \ka\right\}, & \mathrm{if} \; \bar{\rv}=\frac{2\sqrt{3}}{\sqrt{3}-1}\ka,
        \\
        & \Qb(0) = \left\{ \frac{1}{2}\bar{\rv}-\ka \right\}, & \Qa = \left\{\ka\right\}, & \mathrm{if} \; \frac{2\sqrt{3}}{\sqrt{3}-1}\ka < \bar{\rv}.
    \end{array}
\end{align*}
Observe that there are two operating regimes. The regime $0\leq\bar{\rv}<\frac{2\sqrt{3}}{\sqrt{3}-1}\ka$ is the regime of low demand. In this regime, the market has a unique equilibrium and both leader and follower productions increase with demand. The regime $\frac{2\sqrt{3}}{\sqrt{3}-1}\ka<\bar{\rv}$ is the regime of high demand. In this regime, followers produce all their capacity. The leader produces $\frac{1}{2}\bar{\rv}-\ka$, which is less than its production $\frac{1}{2}\bar{\rv}$ in the low demand regime, because it faces less competition now since followers have no capacity left to supply.

By comparing the Stackelberg equilibria to the forward market equilibria, we can see that introducing a forward market does not always increase the total market production. In particular, when $4\ka \leq \bar{\rv} < \frac{2\sqrt{3}}{\sqrt{3}-1}\ka$, the total production $\frac{1}{2}\bar{\rv}+\ka$ with the forward market is less than the total production $\frac{5}{6}\bar{\rv}$ in the Stackelberg market. In this scenario, demand is high and followers produce almost all their capacity in the Stackelberg market. Having followers trade forward contracts give them more incentive to produce and they increase their productions to $\ka$. However, this has the side effect of reducing the competition faced by the leader, and giving it an incentive to withhold its production. The net effect is a decrease in total market production. Since all producers have equal marginal costs, a decrease in total market production implies a decrease in social welfare.

\section{General Structural Results}
\label{sec:StructuralResults}

Building on the analysis in the previous section, we now extend the insights obtained from studying the case of $1$ leader, $2$ followers, and equal marginal costs to general numbers of leaders and followers with marginal costs $\cb$ and $\ca$, respectively, that may be different. 

Our main results (Theorems~\ref{prop:MarketEquilibrium} and~\ref{prop:StackelbergEquilibrium} in Appendix~\ref{app:market}) provide complete characterizations of the symmetric equilibrium productions with and without the forward market, which give a complete picture of when forward contracting mitigates and when it enhances market power. Since these results are technical, we highlight the key properties by characterizing the asymptotic behavior as the number of producers increases (Lemmas~\ref{lem:struct-market}~--~\ref{lem:struct-compare}). Among other properties, we show in Lemmas~\ref{lem:struct-market} and~\ref{lem:struct-market-3} that there is an interval of follower productions just below capacity that are never symmetric equilibria, and that if there are too few leaders relative to followers, then there may not exist symmetric equilibria. We also show, in Lemma~\ref{lem:struct-compare}, that the efficiency loss as a function of the number of producers remains strictly positive even with a large number of followers. 

This section is organized as follows. First, in Sections~\ref{sec:struct-follower} and~\ref{sec:struct-leader}, we characterize the structure of the reactions of the followers to the leaders and vice versa, respectively. Then, in Section~\ref{sec:struct-market}, we characterize the structure of the equilibria. Finally, in Section~\ref{sec:struct-inefficiency}, we characterize the efficiency loss of followers' forward contracting.

Throughout this section, we denote by $\rv_x$ and $\rv_y$ the normalized leader and follower demands respectively and by $\cm$ the normalized marginal cost difference between the leaders and followers:
\begin{align*}
    \rv_x &= \frac{1}{\beta}(\rv-\cb),
    \\
    \rv_y &= \frac{1}{\beta}(\rv-\ca),
    \\
    \cm &= \frac{1}{\beta}(\ca-\cb).
\end{align*}
Note that $\rv_y = \rv_x - \cm$. Since $\ca \geq \cb$, it suffices to restrict our analyses to the case where $\rv_x \geq 0$. We focus on symmetric equilibria, by which we mean equilibria where leaders have symmetric productions and followers have symmetric forward positions (which, by Proposition~\ref{prop:SpotEquilibrium2}, implies that the latter have symmetric productions).

Unless otherwise stated, the proofs for all the results in this section are provided in Appendix~\ref{sec:StructuralProofs}.

\subsection{Follower reaction}
\label{sec:struct-follower}

Suppose all leaders produce a quantity $\qb\in\mathbb{R}_+$ and let $\Qf(\qb)\subseteq\mathbb{R}$ denote the set of all symmetric follower reactions, i.e., for each $\qf\in\Qf(\qb)$ and $j\in\na$,
\begin{align*}
    \phi_j(\qf;\qf\mathbf{1},\qb\mathbf{1}) \geq \phi_j(\bar{\qf};\qf\mathbf{1},\qb\mathbf{1}), \quad \forall \bar{\qf}\in\mathbb{R}.
\end{align*}
Proposition~\ref{prop:PeakerReaction} in Appendix~\ref{app:follower} gives the solution for $\Qf(\qb)$. Observe that $\Qf(\qb)$ has a similar shape to the graph in Figure~\ref{fig:PeakFwd}. We focus on the segment where $\Qf(\qb) = \varnothing$ and highlight key properties that attribute this segment to market manipulation when followers are operating just below capacity.

\begin{lemma}
    The following holds.
    \begin{enumerate}
        \item There exists a unique $\bar{\qa} < \ka$, such that there exists $\qf\in\Qf(\qb)$ that satisfies $\qa_j(\qf\mathbf{1},\qb\mathbf{1})=\qa$ if and only if $0 \leq \qa \leq \bar{\qa}$ or $\qa = \ka$. Moreover,
            \begin{align*}
                \bar{\qa} = \left( 1-O\left( \frac{1}{\na} \right) \right)\ka.
            \end{align*}
        \item There exists a unique $\ubar{\xi}\in\mathbb{R}_+$, such that $\qb \leq \frac{1}{\nb}(\rv_y-\ubar{\xi})$ if and only if there exists $\qf\in\Qf(\qb)$ that satisfies $\qa_j(\qf\mathbf{1},\qb\mathbf{1}) = \ka$, and a unique $\bar{\xi}<\ubar{\xi}$, such that $\qb \geq \frac{1}{\nb}(\rv_y-\bar{\xi})$ if and only if there exists $\qf\in\Qf(\qb)$ that satisfies $\qa_j(\qf\mathbf{1},\qb\mathbf{1}) \leq \bar{\qa}$. Moreover, $\Qf(\qb)=\varnothing$ for all $\qb\in\left( \frac{1}{\nb}\left(\rv_y-\ubar{\xi}\right),\frac{1}{\nb}\left(\rv_y-\bar{\xi}\right) \right)$ and $\ubar{\xi}-\bar{\xi} = \ubar{\xi}\cdot O\left( \frac{1}{\na} \right)$.
    \end{enumerate}
\label{lem:struct-follower}
\end{lemma}

Therefore, there exists an open interval of symmetric leader productions inside which there is no symmetric follower reaction. Due to this interval, there is a set of symmetric follower productions just below $\ka$ that are never equilibria. As $\na$ increases, this set shrinks at the rate $\frac{1}{\na}$. In the limit, all symmetric follower productions could be equilibria. These asymptotic behaviors are consistent with the intuition that followers have less ability to manipulate the market as their numbers increase.

\subsection{Leader reaction}
\label{sec:struct-leader}

Suppose all followers take the forward position $\qf\in\mathbb{R}$ and let $\Qb(\qf)\subseteq\mathbb{R}_+$ denote the set of all symmetric leader reactions, i.e., for each $\qb\in\Qb(\qf)$ and $i\in\nb$,
\begin{align*}
    \psi_i(\qb;\qb\mathbf{1},\qf\mathbf{1}) \geq \psi_i(\bar{\qb};\qb\mathbf{1},\qf\mathbf{1}), \quad \forall \bar{\qb}\in\mathbb{R}_+.
\end{align*}
Proposition~\ref{prop:BaseloadReaction} in Appendix~\ref{app:leader} gives the solution for $\Qb(\qf)$. When $\rv_x \leq \na\ka$, we have $\eta_3 \leq \ka-(\rv_x-\na\ka)$, and one can check that $\Qb(\qf)$ has a similar shape to the graph in Figure~\ref{fig:BaseFwdA}. When $\rv_x > \na\ka$, then $\Qb(\qf)$ differs from the graphs in Figures~\ref{fig:BaseFwdB} and~\ref{fig:BaseFwdC} in that segment (iv) may overlap with segments (iii) and (ii), i.e., there may be up to two reactions. Here, we focus on segment (ii) where the leader reaction is strictly increasing, as well as the discontinuous transition between segment (iv) and segments (ii) or (iii). The following result highlights key properties of segment (ii).
\begin{lemma}
    There exists unique $\ubar{\qf},\bar{\qf}\in\mathbb{R}$, such that $\qf\in\left[\ubar{\qf},\bar{\qf}\right]$ if and only if $\frac{1}{\nb}\left( \rv_x-\cm+\qf \right)\in\Qb(\qf)$. Moreover, the following holds:
    \begin{enumerate}
        \item $\qa_j\left( \qf\mathbf{1},\frac{1}{\nb}\left( \rv_x-\cm+\qf \right)\mathbf{1} \right)=0$ for all $\qf\in\left[ \ubar{\qf},\bar{\qf} \right]$.
        \item $\bar{\qf}-\ubar{\qf} = O\left(\frac{\rv_x}{\nb}\right)$.
    \end{enumerate}
    \label{lem:struct-leader}
\end{lemma}

Therefore, there exists a closed interval $\left[\ubar{\qf},\bar{\qf}\right]$ of symmetric follower productions inside which there is a graph of strictly increasing leader reactions. Moreover, the followers' productions are zero. This is due to leaders using their commitment power to drive the followers out of the market. As $\nb$ increases, this interval shrinks at the rate $\frac{\rv_x}{\nb}$.

The next result highlights key properties of the transition between segment (iii) and (iv).
\begin{lemma}
    Suppose $\rv_x > \na\ka$.
    \begin{enumerate}
        \item There exists a unique $\bar{\qa}<\ka$, such that there exists $\qf\in\mathbb{R}$ and $\qb\in\Qb(\qf)$ that satisfies $\qa_j(\qf\mathbf{1},\qb\mathbf{1})=\qa$ if and only if $0\leq\qa\leq\bar{\qa}$ or $\qa=\ka$. Moreover,
            \begin{align*}
                \bar{\qa}
                =
                \begin{cases}
                    \left( 1-O\left( \frac{\rv_x-\na\ka}{\nb\na} \right) \right)\ka, & \mathit{if} \; \rv_x\leq \na\ka\left( 1+\frac{(\nb+1)\sqrt{\na+1}}{(\sqrt{\na+1}-1)^2}+\frac{(\nb-1)\sqrt{\na+1}}{\sqrt{\na+1}-1} \right),
                    \\
                    0, & \mathit{otherwise}.
                \end{cases}
            \end{align*}
        \item There exists a unique $\bar{\qf}\in\mathbb{R}$, such that $\qf\leq\bar{\qf}$ if and only if there exists $\qb\in\Qb(\qf)$ that satisfies $\qa_j(\qf\mathbf{1},\qb\mathbf{1}) \leq \bar{\qa}$, and a unique $\ubar{\qf}\leq\bar{\qf}$, such that $\qf\geq\ubar{\qf}$ if and only if there exists $\qb\in\Qb(\qf)$ that satisfies $\qa_j(\qf\mathbf{1},\qb\mathbf{1}) = \ka$. Moreover, $|\Qb(\qf)|=2$ for all $\qf\in\left[ \ubar{\qf},\bar{\qf} \right]$. Furthermore, if $\rv_x \leq \na\ka\left( 1+\frac{(\nb+1)\sqrt{\na+1}}{(\sqrt{\na+1}-1)^2} \right)$, then
    \begin{align*}
        \bar{\qf}-\ubar{\qf}
        =
        O\left( \frac{\rv_x-\na\ka}{\nb\sqrt{\na}} \right).
    \end{align*}
    \end{enumerate}
    \label{lem:struct-leader-2}
\end{lemma}

Therefore, there exists an open interval of follower productions $(\bar{\qa},\ka)$ that are never supported by any leader reaction. This interval is due to leaders manipulating the market when followers are operating just below capacity. As $\nb$, $\na$, $\rv_x$ increases, This interval shrinks at the rate $\frac{\rv_x-\na\ka}{\nb\na}$. In the limit, all follower productions can be sustained. Moreover, there is also an interval of follower forward positions $\left[\ubar{\qf},\bar{\qf}\right]$ inside which there are two leader reactions that have different follower productions (one equal to $\ka$ and one less than $\bar{\qa}$). This interval shrinks at the rate $\frac{\rv_x-\na\ka}{\nb\sqrt{\na}}$.

Note that, since the follower production is continuous in $\qf$ and $\qb$, the second claim in Lemma~\ref{lem:struct-leader-2} implies that the leader reaction is discontinuous. This was also observed in the case of one leader and two followers. Also, note that followers' capacity constraints have different impacts on the reactions of the followers and that of the leaders. In the case of followers, it led to non-existence of symmetric reactions. In the case of leaders, there always exists a symmetric reaction but there is a discontinuity in the reaction correspondence.

\subsection{Forward market equilibrium}
\label{sec:struct-market}

We now present structural results for the symmetric equilibria of the forward market. Let $Q\subseteq\mathbb{R}\times\mathbb{R}_+$ denote the set of all symmetric equilibria, i.e., for each $(\qf,\qb)\in Q$, $(\qf\mathbf{1},\qb\mathbf{1})$ is a Nash equilibrium. Theorem~\ref{prop:MarketEquilibrium} in Appendix~\ref{app:market}  gives the solution for $Q$.

First, we focus on the case where $\cm = 0$. The structure of the equilibria is almost identical to that in Section~\ref{sec:SimpleExample:ForwardMarketEquilibrium}; the key difference is that, while there is either one or no equilibria in Section~\ref{sec:SimpleExample:ForwardMarketEquilibrium}, there could be up to two equilibria now. This is highlighted in the following result.
\begin{lemma}
    Suppose $\cm = 0$.
    \begin{enumerate}
        \item There exists a unique $\bar{\qa} < \ka$, such that there exists $(\qf,\qb)\in Q$ that satisfies $\qa_j(\qf\mathbf{1},\qb\mathbf{1})=\qa$ if and only if $0 \leq \qa \leq \bar{\qa}$ or $\qa = \ka$. Moreover,
            \begin{align*}
                \bar{\qa} = \left( 1-O\left( \frac{1}{\na} \right) \right)\ka.
            \end{align*}
        \item There exists a unique $\bar{\rv}_x\in\mathbb{R}_+$, such that $\rv_x \leq \bar{\rv}_x$ if and only if there exists $(\qf,\qb)\in Q$ that satisfies $\qa_j(\qf\mathbf{1},\qb\mathbf{1}) \leq \bar{\qa}$, and a unique $\ubar{\rv}_x\in\mathbb{R}_+$, such that $\rv_x \geq \ubar{\rv}_x$ if and only if there exists $(\qf,\qb)\in Q$ that satisfies $\qa_j(\qf\mathbf{1},\qb\mathbf{1}) = \ka$. Moreover, if
            \begin{align*}
                \nb < \na\sqrt{\na+1}-1,
            \end{align*}
        then $\bar{\rv}_x < \ubar{\rv}_x$, $Q=\varnothing$ for all $\rv_x\in(\bar{\rv}_x,\ubar{\rv}_x)$, and $\ubar{\rv}_x - \bar{\rv}_x = \ubar{\rv}_x\cdot O\left( \frac{1}{\na} \right)$. Otherwise, then $\bar{\rv}_x \geq \ubar{\rv}_x$, $|Q|=2$ for all $\rv_x\in[\ubar{\rv}_x,\bar{\rv}_x]$, and $\bar{\rv}_x - \ubar{\rv}_x = \ubar{\rv}_x\cdot O\left( \frac{1}{\na\sqrt{\na}} \right)$.
    \end{enumerate}
    \label{lem:struct-market}
\end{lemma}

Therefore, there exists an open interval of follower productions $(\bar{\qa},\ka)$ that are never symmetric equilibria. As $\na$ increases, this interval shrinks to the empty set at the rate $\frac{1}{\na}$. The latter is independent of the number of leaders $\nb$ or demand $\rv_x$. However, $\nb$ has an impact on whether there might be no symmetric equilibria or multiple symmetric equilibria. In particular, when $\nb < \na\sqrt{\na+1}-1$, there are no symmetric equilibria when $\ubar{\rv}_x < \rv_x < \bar{\rv}_x$. Otherwise, when $\nb\geq\na\sqrt{\na+1}-1$, there are two symmetric equilibria when $\bar{\rv}_x \leq \rv_x \leq \ubar{\rv}_x$. This behavior illustrates a tradeoff between the number of leaders and followers. The more followers in the system, the more leaders are needed for there to exist symmetric equilibria in the market.

Next, we consider the case where $\cm > 0$. In this case, the structure of the equilibria has an additional feature that was not present when $\cm = 0$. In particular, when demand is low, followers might not supply to the market. The next lemma highlights the structure of the transition to strictly positive follower productions.
\begin{lemma}
    Suppose $\cm > 0$. Let $\zeta_1 = (\nb+1)\cm + \min\left( \nb\na\cm,\nb\na\ka+2\nb\sqrt{\na\ka\cm} \right)$. Then there exists $(\qf,\qb)\in Q$, such that $\qa_j(\qf\mathbf{1},\qb\mathbf{1})=0$ if and only if $\rv_x \leq \zeta_1$. Furthermore, if $\rv_x > (\nb+1)\cm$, then $\qa_j(\qf\mathbf{1},\qb\mathbf{1})=0$ if and only if
    \begin{align*}
        (\qf,\qb)\in \left\{ (\qf,\qb)\in\mathbb{R}\times\mathbb{R}_+ \left|\;
            \qb = \frac{1}{\nb}\left( \rv_x-(\cm-\qf) \right) \; \mathrm{and} \;
            0 \leq \qf \leq \bar{\qf}
    \right. \right\} \subseteq Q,
    \end{align*}
    where $\bar{\qf} > 0$ if $\rv_x < \zeta_1$.
    \label{lem:struct-market-2}
\end{lemma}
The proof is omitted as it is a straightforward observation from Theorem~\ref{prop:MarketEquilibrium}. As the market transitions from zero to strictly positive follower productions, there is a regime of demand where there are multiple equilibria, characterized by leaders increasing supply when followers take larger forward positions. This phenomenon is due to leaders using their commitment power to drive followers out of the market (recall Lemma~\ref{lem:struct-leader}). Therefore, although followers are not supplying to the market, their forward positions have an impact on the efficiency of the equilibrium. Moreover, note that the size of the interval of demand values where this phenomenon occurs is $\min\left( \nb\na\cm,\nb\na\ka+2\nb\sqrt{\na\ka\cm} \right) = \Theta\left( \nb\na \right)$.

When demand is high, the structure of the equilibria is similar to that when $\cm = 0$, in that there could be two or zero equilibria, except that the threshold for $\nb$ now depends on $\cm$. Furthermore, even in the limit as $\na$ tends to infinity, certain follower productions are never equilibria..
\begin{lemma}
    Suppose $\cm > 0$.
    \begin{enumerate}
        \item There exists a unique $\bar{\qa} < \ka$, such that there exists $(\qf,\qb)\in Q$ that satisfies $\qa_j(\qf\mathbf{1},\qb\mathbf{1})=\qa$ if and only if $\qa \leq \bar{\qa}$ or $\qa = \ka$. Moreover,
        \begin{align*}
            \bar{\qa} =
            \begin{cases}
                \left( 1-O\left( \frac{1}{\na} \right) \right)\left( \ka-\frac{(\sqrt{\na+1}-1)^2}{\na}\cm \right), & \mathit{if} \; \cm < \frac{\na\ka}{(\sqrt{\na+1}-1)^2},
                \\
                0, & \mathit{otherwise}.
            \end{cases}
        \end{align*}
        \item There exists a unique $\bar{\rv}_x\in\mathbb{R}_+$, such that $\rv_x \leq \bar{\rv}_x$ if and only if there exists $(\qf,\qb)\in Q$ that satisfies $\qa_j(\qf\mathbf{1},\qb\mathbf{1}) \leq \bar{\qa}$, and a unique $\ubar{\rv}_x\in\mathbb{R}_+$, such that $\rv_x \geq \ubar{\rv}_x$ if and only if there exists $(\qf,\qb)\in Q$ that satisfies $\qa_j(\qf\mathbf{1},\qb\mathbf{1}) = \ka$. Moreover, if
            \begin{align}
                \nb <
                \begin{cases}
                    \frac{(\na+1)\ka-\frac{\na^2+1}{\na^2+(\sqrt{\na+1}-1)^2}\left( \na\ka-\left( \sqrt{\na+1}-1 \right)^2\cm \right)}{\na\cm-\ka+\frac{\na+1}{\na^2+(\sqrt{\na+1}-1)^2}\left( \na\ka-\left( \sqrt{\na+1}-1 \right)^2\cm \right)}, & \mathit{if} \; \cm < \frac{\na\ka}{(\sqrt{\na+1}-1)^2},
                    \\
                    \frac{(\na+1)\ka-\cm}{\na\ka+\cm-\ka+2\sqrt{\na\ka\cm}}, & \mathit{otherwise},
                \end{cases}
                \label{eq:struct-market-3-cond}
            \end{align}
            then $\bar{\rv}_x < \ubar{\rv}_x$ and $Q=\varnothing$ for all $\rv_x\in(\bar{\rv}_x,\ubar{\rv}_x)$. Otherwise, $\bar{\rv}_x \geq \ubar{\rv}_x$ and $|Q|=2$ for all $\rv_x\in\left[ \ubar{\rv}_x,\bar{\rv}_x \right]$.
    \end{enumerate}
    \label{lem:struct-market-3}
\end{lemma}

\subsection{Inefficiency of the forward market}
\label{sec:struct-inefficiency}

We compare the outcome against a Stackelberg competition where followers do not sell forward contracts. Note that the symmetric equilibria of a Stackelberg competition are given by the symmetric reactions of the leaders when followers take neutral forward positions, i.e., $\Qb(0)$, where $\Qb$ is defined in Section~\ref{sec:struct-leader}. Theorem~\ref{prop:StackelbergEquilibrium} in Appendix~\ref{app:market} gives the solution for $\Qb(0)$. The structure is similar to the equilibria of the forward market. We highlight the key features in the following three lemmas.

\begin{lemma}
    Suppose $\cm = 0$.
    \begin{enumerate}
        \item There exists a unique $\bar{\qa} < \ka$, such that there exists $\qb\in \Qb(0)$ that satisfies $\qa_j(\mathbf{0},\qb\mathbf{1})=\qa$ if and only if $\qa \leq \bar{\qa}$ or $\qa = \ka$. Moreover,
            \begin{align*}
                \bar{\qa} = \left( 1+O\left( \frac{1}{\sqrt{\na}} \right) \right)\frac{\ka}{2}.
            \end{align*}
        \item There exists a unique $\bar{\rv}_x\in\mathbb{R}_+$, such that $\rv_x \leq \bar{\rv}_x$ if and only if there exists $\qb\in \Qb(0)$ that satisfies $\qa_j(\mathbf{0},\qb\mathbf{1}) \leq \bar{\qa}$, and a unique $\ubar{\rv}_x\leq\bar{\rv}_x$, such that $\rv_x \geq \ubar{\rv}_x$ if and only if there exists $\qb\in \Qb(0)$ that satisfies $\qa_j(\mathbf{0},\qb\mathbf{1}) = \ka$. Moreover, $|\Qb(0)|=2$ for all $\rv_x \in \left[ \ubar{\rv}_x,\bar{\rv}_x \right]$.
    \end{enumerate}
    \label{lem:struct-stack}
\end{lemma}

Again, we see that there is an open interval of follower productions $(\bar{\qa},\ka)$ that are never symmetric equilibria. However, as $\na$ increases, this interval, instead of shrinking as in the case of the forward market, expands at the rate $\frac{1}{\sqrt{\na}}$ to a size of $\frac{\ka}{2}$. That is, as the followers become more competitive, the leaders are better able to exploit the capacity constraints of the followers.

When $\cm > 0$, followers might not supply to the market. The next lemma highlights the structure of this regime.

\begin{lemma}
    Suppose $\cm > 0$. Let $\zeta_1 = (\nb+1)\cm + \min\left( \nb\na\cm,\nb\na\ka+2\nb\sqrt{\na\ka\cm} \right)$. Then there exists $\qb\in \Qb(0)$ such that $\qa_j(\mathbf{0},\qb\mathbf{1})=0$ if and only if $\rv_x \leq \zeta_1$. Furthermore,
    \begin{align*}
        \qb =
        \begin{cases}
            \frac{1}{\nb+1}\rv_x & \mathrm{if} \; 0 \leq \rv_x < (\nb+1)\cm,
            \\
            \frac{1}{\nb}(\rv_x-\cm) & \mathrm{if} \; (\nb+1)\cm \leq \rv_x \leq \zeta_1,
        \end{cases}
    \end{align*}
    \label{lem:struct-stack-2}
\end{lemma}

The proof is omitted as it is a straightforward observation from Theorem~\ref{prop:StackelbergEquilibrium}. The key insight is that this regime exhibits different behavior depending on whether $\rv_x$ is less than or greater than $(\nb+1)\cm$. The leader productions increase at a faster rate when $\rv_x > (\nb+1)\cm$ because leaders use their commitment power to drive followers out of the market.

When demand is high, the structure of the equilibria is similar to the case when $\cm = 0$, except that the range of follower productions that could be equilibria is now smaller. The larger the value of $\cm$, the smaller the range of supportable follower productions.

\begin{lemma}
    Suppose $\cm > 0$.
    \begin{enumerate}
        \item There exists a unique $\bar{\qa} < \ka$, such that there exist $\qb\in \Qb(0)$ that satisfies $\qa_j(\mathbf{0},\qb\mathbf{1})=\qa$ if and only $\qa\leq\bar{\qa}$ or $\qa = \ka$. Moreover,
            \begin{align*}
                \bar{\qa} =
                \begin{cases}
                    \left( 1+O\left( \frac{1}{\sqrt{\na}} \right)\right) \frac{1}{2}\left( \ka-\frac{(\sqrt{\na+1}-1)^2}{\na}\cm \right),  & \mathit{if} \; \cm \leq \frac{\na\ka}{(\sqrt{\na+1}-1)^2},
                    \\
                    0, & \mathit{otherwise}.
                \end{cases}
            \end{align*}
        \item There exists a unique $\bar{\rv}_x\in\mathbb{R}_+$, such that $\rv_x \leq \bar{\rv}_x$ if and only if there exists $\qb\in \Qb(0)$ that satisfies $\qa_j(\mathbf{0},\qb\mathbf{1}) \leq \bar{\qa}$, and a unique $\ubar{\rv}_x\leq\bar{\rv}_x$, such that $\rv_x \geq \ubar{\rv}_x$ if and only if there exists $\qb\in \Qb(0)$ that satisfies $\qa_j(\mathbf{0},\qb\mathbf{1}) = \ka$. Moreover, $|\Qb(0)|=2$ for all $\rv_x\in\left[ \ubar{\rv}_x,\bar{\rv}_x \right]$.
    \end{enumerate}
    \label{lem:struct-stack-3}
\end{lemma}

We now contrast the efficiency of the equilibria in the forward and Stackelberg markets. Given follower and leader productions $\qa$ and $\qb$ respectively, let $\mathsf{SW}(\qa,\qb)$ denote the social welfare:
\begin{align*}
    \mathsf{SW}(\qa,\qb)
    := \int_0^{\nb\qb+\na\qa} P\left( w \right)\,dw - \left( \nb\cb\qb + \na\ca\qa \right).
\end{align*}
The next lemma highlights that adding a forward market to a Stackelberg market could be inefficient.

\begin{lemma}
    Suppose $\cm = 0$. Let $\ubar{\rv}_x := (\nb+\na+1)\ka$ and $\bar{\rv}_x := \frac{(\nb+1)\sqrt{\na+1}}{2(\sqrt{\na+1}-1)}\na\ka$. Then, for all $\rv_x \in [\ubar{\rv}_x,\bar{\rv}_x]$, there exists $(\qf,\qb)\in Q$ and $\qb_S\in\Qb(0)$ such that
    \begin{align*}
        \nb\qb+\na\qa_j(\qf\mathbf{1},\qb\mathbf{1}) &< \nb\qb_S+\na\qa_j(\mathbf{0},\qb_S\mathbf{1}),
        \\
        \mathsf{SW}(\qa_j(\qf\mathbf{1},\qb\mathbf{1}),\qb) &< \mathsf{SW}(\qa_j(\mathbf{0},\qb_S\mathbf{1}),\qb_S).
    \end{align*}
    Moreover,
    \begin{align*}
        \frac{\nb\qb_S+\na\qa_j(\mathbf{0},\qb_S\mathbf{1})}{\nb\qb+\na\qa_j(\qf\mathbf{1},\qb\mathbf{1})} & \leq \frac{(\nb\na+\nb+\na)(\nb+1)}{\nb(\nb+1)(\na+1)+2(\na+1-\sqrt{\na+1})},
        \\
        \frac{\mathsf{SW}(\qa_j(\mathbf{0},\qb_S\mathbf{1}),\qb_S)}{\mathsf{SW}(\qa_j(\qf\mathbf{1},\qb\mathbf{1}),\qb)} & \leq \frac{(\nb+1)^2(\nb\na+\nb+\na)(\nb\na+\nb+\na+2)}{(\na+1) \left( (\nb^2+\nb+2)\sqrt{\na+1}-2\right)\left( (\nb^2+3\nb)\sqrt{\na+1}+2 \right)},
    \end{align*}
    where the inequalities are tight.
    \label{lem:struct-compare}
\end{lemma}

This inefficiency is attributed to equilibria in the forward market where followers produce $\ka$ while there are equilibria in the Stackelberg market where followers produce strictly less than $\ka$. Therefore, this inefficiency is due to leaders exploiting the capacity constraints of the followers in the forward market. To see this, note that this inefficiency does not disappear even with a large number of followers:
\begin{align*}
    \lim_{\na\rightarrow\infty} \frac{\nb\qb_S+\na\qa_j(\mathbf{0},\qb_S\mathbf{1})}{\nb\qb+\na\qa_j(\qf\mathbf{1},\qb\mathbf{1})} &\leq \frac{(\nb+1)^2}{\nb^2+\nb+2},
    \\
    \lim_{\na\rightarrow\infty} \frac{\mathsf{SW}(\qa_j(\mathbf{0},\qb_S\mathbf{1}),\qb_S)}{\mathsf{SW}(\qa_j(\qf\mathbf{1},\qb\mathbf{1}),\qb)} & \leq \frac{(\nb+1)^4}{(\nb^2+\nb+2)(\nb^2+3\nb)}.
\end{align*}
On the other hand, this inefficiency disappears with a large number of leaders:
\begin{align*}
    \lim_{\nb\rightarrow\infty} \frac{\nb\qb_S+\na\qa_j(\mathbf{0},\qb_S\mathbf{1})}{\nb\qb+\na\qa_j(\qf\mathbf{1},\qb\mathbf{1})} &\leq 1,
    \\
    \lim_{\nb\rightarrow\infty} \frac{\mathsf{SW}(\qa_j(\mathbf{0},\qb_S\mathbf{1}),\qb_S)}{\mathsf{SW}(\qa_j(\qf\mathbf{1},\qb\mathbf{1}),\qb)} & \leq 1.
\end{align*}
The statement of Lemma~\ref{lem:struct-compare} does not specify whether there exists forward equilibria that are equally or more efficient than Stackelberg equilibria. However, it is possible to impose further conditions on the system and demand such that the Stackelberg equilibria are always strictly more efficient.

The same approach in the proof of Lemma~\ref{lem:struct-compare} can be used to obtain bounds on the production and efficiency losses when $\cm > 0$. However, the bounds are more complicated and depend on $\cm$ and $\ka$.

\section{Conclusion}

Forward contracts make up a significant share of trade in many markets ranging from finance to cloud computing to commodities, e.g., gas and electricity. The general view is that forward trading improves the efficiency of markets by providing a mechanism for participants to hedge risks and mitigating market power. Since the seminal result by~\cite{allaz1993cournot} that proved forward contracts can mitigate market power, there has been significant interest in understanding the generality of this phenomenon. In this work, we show that leader-follower interactions may cause forward contracting to be inefficient. This is because forward contracting increases followers’ outputs which may create opportunities for leaders to exploit the capacity constraints of followers. Furthermore, we show that the efficiency loss remains strictly positive even with a large number of followers (Lemma~\ref{lem:struct-compare}), and hence this inefficiency may be attributed to oligopoly leaders. Our results contrast with prior work that also showed that forward markets may not mitigate market power as those were due to endogenous investment~\cite{murphy2010impact} or transmission congestion~\cite{kamat2004two}. Our results are important due to the prevalence of forward trading in some industries where there are leader-follower relationships between the firms (e.g. gas and electricity).

Furthermore, due to our closed-form expressions for every symmetric equilibria (Theorems~\ref{prop:MarketEquilibrium} and~\ref{prop:StackelbergEquilibrium}), we are able to characterize the behavior of the system explicitly which provide strategic insights. One key characterisation is that there is an interval of follower productions just below capacity that are never symmetric equilibria (including symmetric leader/follower reactions) (Lemmas~\ref{lem:struct-follower},~\ref{lem:struct-leader-2},~\ref{lem:struct-market}, and~\ref{lem:struct-market-3}). This phenomenon may be attributed to oligopoly followers since this interval shrinks as the number of followers increase. Another key characterisation is that, if there are too few leaders relative to followers, then there may not exist symmetric equilibria (Lemmas~\ref{lem:struct-market} and~\ref{lem:struct-market-3}). This tradeoff shows that, the more competitive the spot market, the easier it is for leaders to exploit followers’ capacity constraints, which is reminiscent of the first-mover advantage in the classical Stackelberg game. Therefore, temporal constraints may create differences in market power between firms.

The strategic interactions that we observed in this work could provide insights into behaviour in other settings. To this point, we have not addressed the impact of the leader’s production inflexibility on the efficiency of the market, as our focus was on the efficiency of forward contracting. Nevertheless, the insights obtained from our results allow us to conjecture the possible impacts. It is well known that Stackelberg competition is less efficient than Cournot. Therefore, the natural inference is that constraining the leader to choose productions in the first stage (versus it selling forward contracts only and choosing productions in the second stage) would decrease the social welfare. However, based on the phenomenon shown in our work, it is plausible that this intuition is not always true. If the leader did not have to choose productions in the first stage, the added competition would cause followers to produce more. But, if this causes followers to be capacity constrained, then there might be opportunities for producers to exploit those constraints, resulting in a loss of efficiency or non-existence of equilibria altogether.

\bibliographystyle{plain}
\bibliography{bibfile}

    \appendix


\section{Closed-Form Solutions for Symmetric Equilibria}
\label{app:closed-form}
In this Appendix, we derive closed-form expressions for the symmetric follower reactions, symmetric leader reactions, the symmetric forward market equilibria, and the symmetric Stackelberg equilibria. These results are used to derive the structural results in Section~\ref{sec:StructuralResults}. We denote by $\rv_x$ and $\rv_y$ the normalized leader and follower demands respectively, and by $\cm$ the normalized marginal cost gap:
\begin{align*}
    \rv_x &= \frac{1}{\beta}(\rv-\cb),
    \\
    \rv_y &= \frac{1}{\beta}(\rv-\ca),
    \\
    \cm &= \frac{1}{\beta}(\ca-\cb).
\end{align*}
We use the following notation. For scalars $z,a,b \in \mathbb{R}$ such that $a \leq b$, let 
\begin{align*}
[z]_a^b := 
\begin{cases}
a, & \mathrm{if} \; z \leq a,
\\
b, & \mathrm{if} \; z \geq b,
\\
z, & \mathrm{otherwise}.
\end{cases}
\end{align*}
We will use the following properties:
\begin{enumerate}[(i)]
\item For any $c\in \mathbb{R}$, $c + \left[z\right]_a^b = \left[z + c\right]_{a + c}^{b + c}$.
\item If $c > 0$, then $c\left[z\right]_a^b = \left[cz\right]_{ca}^{cb}$. 
\item If $c < 0$, then $c\left[z\right]_a^b = \left[cz\right]_{cb}^{ca}$. 
\end{enumerate}

\subsection{Spot Market Analyses}

\begin{proposition}
\label{prop:SpotNash}
Fix a follower $l\in\na$ and suppose $f_j = \qf$ for every $j\neq l$. There is a unique Nash equilibrium $\qap$ in the spot market such that, for each $j \neq l$,
\begin{align}
y_{j} = \left[\frac{1}{\na}\left(\rv_y + \qf - \sum_{i=1}^{M} x_i - y_l\right)\right]_0^\ka.
\label{eq:spot-equilibrium-1}
\end{align}
\end{proposition}

    \proof{Proof.}
The uniqueness of the Nash equilibrium follows from Theorem 5 of~\cite{jing1999spatial}. Each follower $j\in \na$ has a strategy set $[0,\ka]$ which is compact. Its payoff function in the spot market $\phi_j^{(s)}$ is continuous in all arguments and is strictly concave in $y_j$. Thus, from the Karush-Kuhn-Tucker (KKT) conditions, we infer that $\qap\in[0,\ka]^\na$ is a Nash equilibrium of the spot market, if and only if there exists $\boldsymbol\lambda,\boldsymbol\mu\in\mathbb{R}_+^\na$ such that, for each $j\in\na$:
\begin{align}
\nabla_{y_j}\left[\phi_j^{(s)}(y_j;\mathbf{y}_{-j})+\lambda_j y_j + \mu_j (k - y_j)\right] &= 0,
\label{eq:SpotOptimalityFOC}
\\
\lambda_j y_j = \mu_j (k - y_j) &= 0.
\label{eq:SpotOptimalityComplSlack}
\end{align}
Take any $j\neq l$. Expanding the LHS of~\eqref{eq:SpotOptimalityFOC} gives:
\begin{align*}
\nabla_{y_j}\left[\phi_j^{(s)}(y_j;\mathbf{y}_{-j})+\lambda_j y_j + \mu_j (k - y_j)\right]
&= 
\beta\left(\alpha_y + f - \sum_{i=1}^M x_i - y_j - \sum_{j'=1}^{\na} y_{j'} \right) + \lambda_j - \mu_j
\\
&=
\beta\left(\rv_y + f - \sum_{i=1}^{\nb}x_i - y_l - \na y_{j} \right) + \lambda_j - \mu_j.
\end{align*}
Suppose $0 < y_j < \ka$. Then~\eqref{eq:SpotOptimalityComplSlack} imply that $\lambda_j = \mu_j = 0$. From~\eqref{eq:SpotOptimalityFOC}, we obtain
\begin{align}
y_j
=\frac{1}{\na}\left(\rv_y + \qf - \sum_{i=1}^{\nb}x_i - y_l\right).
\label{eq:SpotNash1}
\end{align}
Suppose $y_j = 0$. Then~\eqref{eq:SpotOptimalityComplSlack} imply that $\mu_j = 0$. From~\eqref{eq:SpotOptimalityFOC}, we obtain
\begin{align}
-\left(\rv_y + \qf - \sum_{i=1}^{\nb}x_i - y_l\right)= \lambda_j \geq 0.
\label{eq:SpotNash2}
\end{align}
Suppose $y_j = \ka$. Then~\eqref{eq:SpotOptimalityComplSlack} imply that $\lambda_j = 0$. From~\eqref{eq:SpotOptimalityFOC}, we obtain
\begin{align}
\left(\rv_y + \qf - \sum_{i=1}^{\nb}x_i - y_l - \na\ka\right) = \mu_j \geq 0.
\label{eq:SpotNash3}
\end{align}
Since $0 \leq y_j \leq \ka$, \eqref{eq:SpotNash1}~--~\eqref{eq:SpotNash3} together imply that
\begin{align*}
y_j
&=
\begin{cases}
0, & \mathrm{if} \; \frac{1}{\na}\left(\rv_y + \qf - \sum_{i=1}^{\nb}x_i - y_l\right) \leq 0,
\\
\ka, & \mathrm{if} \; \frac{1}{\na}\left(\rv_y + \qf - \sum_{i=1}^{\nb}x_i - y_l\right) \geq \ka,
\\
\frac{1}{\na}\left(\rv_y + \qf - \sum_{i=1}^{\nb}x_i - y_l\right), & \mathrm{otherwise},
\end{cases}
\end{align*}
which is equivalent to~\eqref{eq:spot-equilibrium-1}.
\endproof

\begin{proposition}
\label{prop:SpotEquilibrium2}
Suppose $f_j = \qf$ for every $j\in N$. There is a unique Nash equilibrium in the spot market, given by
\begin{align}
y_j
&=
\left[\frac{1}{\na+1}\left(\rv_y + \qf - \sum_{i=1}^{\nb}\qbpj\right)\right]_0^k.
\label{eq:spot-equilibrium-2}
\end{align}
\end{proposition}

    \proof{Proof.}
The uniqueness of the Nash equilibrium follows from Theorem 5 of~\cite{jing1999spatial}. Thus, it suffices to show that the given productions form a Nash equilibrium. From the optimality conditions in~\eqref{eq:SpotOptimalityFOC}~--~\eqref{eq:SpotOptimalityComplSlack}, we infer that $y\in\left[0,\ka\right]$ is a symmetric Nash equilibrium in the spot market, if and only if there exists scalars $\lambda,\mu\in\mathbb{R}_+$ such that,
\begin{align*}
\beta\left(\rv_y + f - \sum_{i=1}^{\nb}x_i - (\na+1)y \right) + \lambda - \mu &= 0,
\\
\lambda y = \mu \left(\ka - y\right) &= 0.
\end{align*}
Let
\begin{align*}
\lambda
&=
\left[-\beta\left(\rv_y + \qf - \sum_{i=1}^{\nb}x_i - (\na+1) y\right)\right]_0^\infty,
\\
\mu
&=
\left[\beta\left(\rv_y + \qf - \sum_{i=1}^{\nb}x_i - (\na+1) y\right)\right]_0^\infty.
\end{align*}
It is straightforward to show that $y$ defined in~\eqref{eq:spot-equilibrium-2}, and $\lambda, \mu$ defined above, together satisfy the optimality conditions.
\endproof

\subsection{Follower Reaction Analyses}
\label{app:follower}

\begin{proposition}
\label{prop:PeakerReaction}
Fix the leaders' productions $\qbp\in\mathbb{R}_+^\nb$. Let $\Qf\subseteq\mathbb{R}$ denote the set of symmetric follower reactions, i.e., for each $f\in \Qf$ and $j\in \na$,
\begin{align}
\label{eq:ProofOfProp3QfDefinition}
\pfpif\left(\qf;\qf\mathbf{1},\qbp\right)
\geq
\pfpif\left(\bar{f};\qf\mathbf{1},\qbp\right),
\quad
\forall 
\bar{f}\in\mathbb{R}.
\end{align}
Let $\xi := \rv_y - \sum_{i=1}^{\nb} x_i$. Then,
\begin{align}
\Qf
&=
\begin{cases}
\left(-\infty, -\xi\right], & \mathit{if}\; \xi < 0,
\\
\left\{\frac{\na-1}{\na^2+1}\xi\right\}, & \mathit{if}\; 0 \leq \xi \leq \frac{(\na^2+1)(\na-1)}{\na^2-2\sqrt{\na}+1}\ka,
\\
\varnothing, & \mathit{if}\; \frac{(\na^2+1)(\na-1)}{\na^2-2\sqrt{\na}+1}\ka < \xi < \left(\na+1\right)\ka,
\\
\left[-\xi + (\na+1)\ka,\infty\right), & \mathit{if}\; \left(\na+1\right)\ka \leq \xi.
\end{cases}
\label{eq:follower-reaction}
\end{align}
Moreover, for each $\qf\in\Qf$, 
\begin{align*}
y_j(\qf\mathbf{1},\qbp) = 0 
& \iff \xi \leq 0,
\\
0 < y_j(\qf\mathbf{1},\qbp) < \ka 
& \iff 0 < \xi \leq \frac{(\na^2+1)(\na-1)}{\na^2-2\sqrt{\na}+1}\ka,
\\
y_j(\qf\mathbf{1},\qbp) = \ka
& \iff (\na+1)\ka \leq \xi.
\end{align*}
\end{proposition}

    \proof{Proof.}
The proof proceeds in three steps. In step 1, we reformulate a follower's payoff maximization problem into a problem involving its production quantity only. In step 2, we compute its payoff maximizing production quantity. In step 3, we compute the symmetric follower forward positions that satisfy the condition that every follower is producing at its payoff maximizing quantity. The latter gives the set of symmetric follower reactions.

\emph{Step 1:} Fix a follower $l\in \na$ and suppose $f_j = f$ for every $j\neq l$. Using Proposition~\ref{prop:SpotNash} to substitute for $y_j(f_j;\qfp_{-j},\qbp)$ for every $j\neq l$, we infer that the total production in the spot market is given by
\begin{align*}
\sum_{j=1}^{\na}y_j(\qfpi;\qfp_{-j},\qbp)
&=
y_l(f_l;\qf\mathbf{1},\qbp) + \left(\na-1\right)\left[\frac{1}{\na}\left(\rv_y+\qf-\sum_{i=1}^{\nb}x_i-y_l(f_l;\qf\mathbf{1},\qbp)\right)\right]_0^\ka
\\
&=
y_l(f_l;\qf\mathbf{1},\qbp) + \left[\frac{\na-1}{\na}\left(\rv_y+\qf-\sum_{i=1}^{\nb}x_i-y_l(f_l;\qf\mathbf{1},\qbp)\right)\right]_0^{(\na-1)\ka}
\\
&=
\left[\frac{\na-1}{\na}\left(\rv_y+\qf-\sum_{i=1}^{\nb}x_i\right) + \frac{1}{\na}y_l(f_l;\qf\mathbf{1},\qbp)\right]_{y_l(f_l;\qf\mathbf{1},\qbp)}^{y_l(f_l;\qf\mathbf{1},\qbp)+(\na-1)\ka},
\end{align*}
By substituting the above into follower $l$'s payoff, and using the fact that $y_l(\mathbb{R};\qf\mathbf{1},\qbp) = \left[0,\ka\right]$, we obtain
\begin{align}
\sup_{f_l\in\mathbb{R}} \phi_l(f_l;\qf\mathbf{1},\qbp)
&=
\sup_{f_l\in\mathbb{R}} \left(P\left(\left[\frac{\na-1}{\na}\left(\rv_y+\qf-\sum_{i=1}^{\nb}x_i\right) + \frac{1}{\na}y_l(f_l;\qf\mathbf{1},\qbp)\right]_{y_l(f_l;\qf\mathbf{1},\qbp)}^{y_l(f_l;\qf\mathbf{1},\qbp)+(\na-1)\ka} 
\right.\right.
\nonumber\\
&\qquad\qquad\qquad\left.\left. + \sum_{i=1}^{\nb} x_i\right)-\ca\right)\cdot y_l(f_l;\qf\mathbf{1},\qbp)
\label{eq:follower-profit-aux}
\\
&=\sup_{y\in\left[0,\ka\right]} \hat{\phi}_l(y;\qf,\qbp),
\label{eq:follower-nash-aux}
\end{align}
where
\begin{align*}
\hat{\phi}_l(y;\qf,\qbp)
:=\left(P\left(\left[\frac{\na-1}{\na}\left(\rv_y+\qf-\sum_{i=1}^{\nb}x_i\right) + \frac{1}{\na}y\right]_{y}^{y+(\na-1)\ka} + \sum_{i=1}^{\nb} x_i\right)-\ca\right)\cdot y.
\end{align*}

\emph{Step 2:} We solve for the solution to~\eqref{eq:follower-nash-aux}. Substituting for the demand function yields
\begin{align*}
\hat{\phi}_l(y;\qf,\qbp)
&=
\beta\left(\xi - \left[\frac{\na-1}{\na}\left(\xi + \qf\right) + \frac{1}{\na} y\right]_{y}^{y + (\na-1)\ka}\right) y
\\
&=
\begin{cases}
\beta\left(\xi - y\right)y, & \mathrm{if} \; \eqref{eq:follower-nash-range-1} \; \mathrm{holds},
\\
\begin{cases}
\beta\left(\frac{1}{\na}\xi - \frac{\na-1}{\na}\qf - \frac{1}{\na}y\right) y, & \mathrm{if} \; 0 \leq y < \xi + \qf,
\\
\beta\left(\xi - y\right)y, & \mathrm{if} \; \ka \geq y \geq \xi + \qf,
\end{cases} & \mathrm{if} \; \eqref{eq:follower-nash-range-2} \; \mathrm{holds},
\\
\beta\left(\frac{1}{\na}\xi - \frac{\na-1}{\na}\qf - \frac{1}{\na}y\right) y, & \mathrm{if} \; \eqref{eq:follower-nash-range-3} \; \mathrm{holds},
\\
\begin{cases}
\beta\left(\xi - y - (\na-1)\ka\right) y, & \mathrm{if} \; 0 \leq y \leq \xi + \qf - \na\ka,
\\
\beta\left(\frac{1}{\na}\xi - \frac{\na-1}{\na}\qf - \frac{1}{\na}y\right) y, & \mathrm{if} \; \ka \geq y > \xi + \qf - \na\ka,
\end{cases} & \mathrm{if} \; \eqref{eq:follower-nash-range-4} \; \mathrm{holds},
\\
\beta\left(\xi - y - (\na-1)\ka\right) y, & \mathrm{if} \; \eqref{eq:follower-nash-range-5} \; \mathrm{holds},
\end{cases}
\end{align*}
where the second equality follows from the fact that $y \in \left[0,\ka\right]$ and the five cases~\eqref{eq:follower-nash-range-1}~--~\eqref{eq:follower-nash-range-5} are defined by
\begin{subequations}
\begin{align}
\xi + \qf &\leq 0,
\label{eq:follower-nash-range-1}
\\
0 < \xi + \qf &< \ka,
\label{eq:follower-nash-range-2}
\\
\ka \leq \xi + \qf & \leq \na\ka,
\label{eq:follower-nash-range-3}
\\
\na\ka < \xi + \qf & < (\na+1)\ka,
\label{eq:follower-nash-range-4}
\\
(\na+1)\ka \leq \xi + \qf.
\label{eq:follower-nash-range-5}
\end{align}
\end{subequations}
We analyze each case separately. 

\paragraph{Case (i):} $\xi + \qf \leq 0$. Then $\hat{\phi}_l(y;\qf,\qbp)$ is a smooth function in $y$ over the interval $\left[0,\ka\right]$. The first and second derivatives are given by
\begin{align*}
\frac{\partial}{\partial y}\hat{\phi}_l(y;\qf,\qbp)
&=\beta\left(\xi - 2 y\right),
\\
\frac{\partial^2}{\partial y^2}\hat{\phi}_l(y;\qf,\qbp)
&=-2\beta < 0,
\end{align*}
which implies that $\hat{\phi}_l(y;\qf,\qbp)$ is strictly concave in $y$. Thus, $y$ is a solution to~\eqref{eq:follower-nash-aux} if and only if it satisfies the following first order optimality conditions:
\begin{align}
\frac{\partial^+}{\partial y}\hat{\phi}_l(y;\qf,\qbp) \leq 0, & \quad \mathrm{if} \; 0 \leq y < \ka,
\label{eq:follow-nash-aux-concave-foc-1}\\
\frac{\partial^-}{\partial y}\hat{\phi}_l(y;\qf,\qbp) \geq 0, & \quad \mathrm{if} \; 0 < y \leq \ka.
\label{eq:follow-nash-aux-concave-foc-2}
\end{align}
It is straightforward to show that there is a unique solution  given by
\begin{align}
y
&=
\left[\frac{1}{2}\xi\right]_0^\ka.
\label{eq:follower-optimal-production-i}
\end{align}

\paragraph{Case (ii):} $0 < \xi + \qf < \ka$. Then $\hat{\phi}_l(y;\qf,\qbp)$ is a piecewise smooth function in $y$ over the interval $\left[0,\ka\right]$. The first and second derivatives are given by
\begin{align*}
\frac{\partial}{\partial y}\hat{\phi}_l(y;\qf,\qbp)
&=
\begin{cases}
\beta\left(\frac{1}{\na}\xi-\frac{\na-1}{\na}\qf - \frac{2}{\na}y\right), & \mathrm{if} \; 0 \leq y < \xi + \qf,
\\
\beta\left(\xi - 2y\right), & \mathrm{if} \; \ka \geq y > \xi + \qf,
\end{cases}
\\
\frac{\partial^2}{\partial y^2}\hat{\phi}_l(y;\qf,\qbp)
&=
\begin{cases}
-\frac{2}{\na}\beta, & \mathrm{if} \; 0 \leq y < \xi + \qf,
\\
-2\beta, & \mathrm{if} \; \ka \geq y > \xi + \qf,
\end{cases}
\\
&< 0.
\end{align*}
Moreover, we have
\begin{align*}
\left.\frac{\partial^+}{\partial y}\hat{\phi}_l(y;\qf,\qbp)\right|_{y = \xi + \qf}
&=
\beta\left(- \xi - 2\qf\right)
\\
&=
\frac{1}{\na}\beta\left(- \na\xi - 2\na\qf\right)
\\
&\leq
\frac{1}{\na}\beta\left(- \xi + (\na-1)\qf - 2\na\qf\right)
\\
&=
\frac{1}{\na}\beta\left(\xi - (\na-1)\qf - 2\xi - 2\qf\right)
\\
&=
\left.\frac{\partial^-}{\partial y}\hat{\phi}_l(y;\qf,\qbp)\right|_{y = \xi + \qf},
\end{align*}
where the inequality follows from the fact that $\xi + \qf > 0$. Thus, $\hat{\phi}_l(y;\qf,\qbp)$ is concave in $y$ over $\left[0,\ka\right]$. Therefore, $y$ is a solution to~\eqref{eq:follower-nash-aux} if and only if it satisfies the first order optimality conditions~\eqref{eq:follow-nash-aux-concave-foc-1}~--~\eqref{eq:follow-nash-aux-concave-foc-2}. It is straightforward to show that there is a unique solution given by
\begin{align}
y
=
\begin{cases}
0, & \mathrm{if} \; \xi \leq (\na-1)\qf,
\\
\frac{1}{2}\left(\xi - \left(\na-1\right)\qf\right), & \mathrm{if} \; \xi > \max((\na-1)\qf,-(\na+1)\qf),
\\
\xi + \qf, & \mathrm{if} \; -2\qf \leq \xi \leq -(\na+1)\qf,
\\
\frac{1}{2}\xi, & \mathrm{if} \; \xi < \min(2\ka,-2\qf),
\\
\ka, & \mathrm{if} \; \xi \geq 2\ka.
\end{cases}
\label{eq:follower-optimal-production-ii}
\end{align}

\paragraph{Case (iii):} $\ka \leq \xi + \qf \leq \na\ka$. Then $\hat{\phi}_l(y;\qf,\qbp)$ is a smooth function in $y$ over the interval $\left[0,\ka\right]$. The first and second derivatives are given by
\begin{align*}
\frac{\partial}{\partial y}\hat{\phi}_l(y;\qf,\qbp)
&=\beta\left(\frac{1}{\na}\xi - \frac{\na-1}{\na}\qf - \frac{2}{\na}y\right),
\\
\frac{\partial^2}{\partial y^2}\hat{\phi}_l(y;\qf,\qbp)
&=-\frac{2}{\na}\beta < 0,
\end{align*}
which implies that $\hat{\phi}_l(y;\qf,\qbp)$ is strictly concave in $y$. Thus, $y$ is a solution to~\eqref{eq:follower-nash-aux} if and only if it satisfies the first order optimality conditions~\eqref{eq:follow-nash-aux-concave-foc-1}~--~\eqref{eq:follow-nash-aux-concave-foc-2}. It is straightforward to show that there is a unique solution given by
\begin{align}
y
&=
\left[\frac{1}{2}\left(\xi-(\na-1)\qf\right)\right]_0^\ka.
\label{eq:follower-optimal-production-iii}
\end{align}

\paragraph{Case (iv):} $\na\ka < \xi + \qf < (\na+1)\ka$. Then $\hat{\phi}_l(y;\qf,\qbp)$ is a piecewise smooth function in $y$ over the interval $\left[0,\ka\right]$. The first and second derivatives are given by
\begin{align*}
\frac{\partial}{\partial y}\hat{\phi}_l(y;\qf,\qbp)
&=
\begin{cases}
\beta\left(\xi - 2y - (\na-1)\ka\right), & \mathrm{if} \; 0 \leq y < \xi + \qf - \na\ka,
\\
\beta\left(\frac{1}{\na}\xi - \frac{\na-1}{\na}\qf - \frac{2}{\na}y\right), & \mathrm{if} \; \ka \geq y > \xi + \qf - \na\ka,
\end{cases}
\\
\frac{\partial^2}{\partial y^2}\hat{\phi}_l(y;\qf,\qbp)
&=
\begin{cases}
-2\beta, & \mathrm{if} \; 0 \leq y < \xi + \qf - \na\ka,
\\
-\frac{2}{\na}\beta, & \mathrm{if} \; \ka \geq y > \xi + \qf - \na\ka,
\end{cases}
\\
&< 0.
\end{align*}
It is straightforward to check that $\hat{\phi}_l(y;\qf,\qbp)$ is not concave in $y$ over the interval $\left[0,\ka\right]$. However, $\hat{\phi}_l(y;\qf,\qbp)$ is piecewise concave in $y$. Therefore, solve the following sub-problems:
\begin{align}
\sup_{y\in\left[0,\xi+\qf-(\na-1)\ka\right]} \hat{\phi}_l(y;\qf,\qbp),
\label{eq:follower-nash-iv-prob1}
\end{align}
and
\begin{align}
\sup_{y\in\left[\xi+\qf-(\na-1)\ka,\ka\right]} \hat{\phi}_l(y;\qf,\qbp).
\label{eq:follower-nash-iv-prob2}
\end{align}
The solution of the sub-problem with the larger optimal value is the solution to~\eqref{eq:follower-nash-aux}. Using the first-order optimality conditions, the unique solution to~\eqref{eq:follower-nash-iv-prob1} is given by
\begin{align*}
y
=
\left[\frac{1}{2}\left(\xi - (\na-1)\ka\right)\right]_0^{\xi + \qf - \na\ka} =: z_1,
\end{align*}
and that to~\eqref{eq:follower-nash-iv-prob2} is given by
\begin{align*}
y
=
\left[\frac{1}{2}\left(\xi - (\na-1)\qf\right)\right]_{\xi + \qf -\na\ka}^\ka =: z_2.
\end{align*}
Therefore, the solution(s) to~\eqref{eq:follower-nash-aux} are given by:
\begin{subequations}
\begin{align}
y = z_1, &\quad \mathrm{if} \; \hat{\phi}_l(z_1;\qf,\qbp) > \hat{\phi}_l(z_2;\qf,\qbp),
\\
y = z_2, &\quad \mathrm{if} \; \hat{\phi}_l(z_1;\qf,\qbp) < \hat{\phi}_l(z_2;\qf,\qbp),
\\
y = z_1 \; \mathrm{or} \; z_2, & \quad \mathrm{if} \; \hat{\phi}_l(z_1;\qf,\qbp) = \hat{\phi}_l(z_2;\qf,\qbp).
\end{align}
\label{eq:follower-optimal-production-iv}
\end{subequations}

\paragraph{Case (v):} $(\na+1)\ka \leq \xi + \qf$. Then $\hat{\phi}_l(y;\qf,\qbp)$ is a smooth function in $y$ over the interval $\left[0,\ka\right]$. The first and second derivatives are given by
\begin{align*}
\frac{\partial}{\partial y}\hat{\phi}_l(y;\qf,\qbp)
&=\beta\left(\xi - (\na-1)\ka - 2y\right),
\\
\frac{\partial^2}{\partial y^2}\hat{\phi}_l(y;\qf,\qbp)
&=-2\beta < 0,
\end{align*}
which implies that $\hat{\phi}_l(y;\qf,\qbp)$ is strictly concave in $y$. Therefore, $y$ is a solution to~\eqref{eq:follower-nash-aux} if and only if it satisfies the first order optimality conditions~\eqref{eq:follow-nash-aux-concave-foc-1}~--~\eqref{eq:follow-nash-aux-concave-foc-2}. It is straightforward to show that there is a unique solution given by
\begin{align}
y
&=
\left[\frac{1}{2}\left(\xi-(\na-1)\ka\right)\right]_0^\ka.
\label{eq:follower-optimal-production-v}
\end{align}

\emph{Step 3:} We solve for the symmetric follower forward positions that satisfy the condition that every follower is producing at its payoff maximizing quantity. The latter gives the set of symmetric follower reactions since
\begin{align}
&\phi_l\left(\qf;\qf\mathbf{1},\qbp\right) \geq \phi_l\left(\bar{\qf};\qf\mathbf{1},\qbp\right),\quad \forall\bar{\qf}\in\mathbb{R}
\nonumber\\
\iff &
\hat{\phi}_l\left(y_l(\qf;\qf\mathbf{1},\qbp);\qf,\qbp\right) \geq \hat{\phi}_l\left(y_l(\bar{\qf};\qf\mathbf{1},\qbp);\qf,\qbp\right), \quad \forall \bar{\qf}\in\mathbb{R}
\nonumber\\
\iff &
\hat{\phi}_l\left(y;\qf,\qbp\right) \geq \hat{\phi}_l\left(\bar{y};\qf\mathbf{1},\qbp);\qf,\qbp\right), \quad \forall \bar{y}\in\left[0,\ka\right], 
\; \mathrm{and} \; y = \left[\frac{1}{\na+1}\left(\xi + \qf\right)\right]_0^\ka.
\label{eq:follower-nash-symmetric-production}
\end{align}
The first equivalence is due to~\eqref{eq:follower-profit-aux}. The second equivalence is due to the fact that $y_l(\qf;\qf\mathbf{1},\qbp) = \left[\frac{1}{\na+1}\left(\xi + \qf\right)\right]_0^\ka$ and $y_l\left(\mathbb{R};\qf\mathbf{1},\qbp\right) = \left[0,\ka\right]$. We divide the analyses into five cases depending on the value of $\xi + \qf$. 

\paragraph{Case (i):} $\xi + \qf \leq 0$. Note that $y$ is given by~\eqref{eq:follower-optimal-production-i}. Therefore, the symmetric follower reactions are given by:
\begin{align}
\left[\frac{1}{2}\xi\right]_0^\ka = \left[\frac{1}{\na+1}\left(\xi +\qf\right)\right]_0^\ka \;\mathrm{and} \; \eqref{eq:follower-nash-range-1} \; \mathrm{holds}
&\iff
\left[\frac{1}{2}\xi\right]_0^\ka = 0 \;\mathrm{and} \; \eqref{eq:follower-nash-range-1} \; \mathrm{holds}
\nonumber\\
&\iff
\xi \leq 0 \;\mathrm{and} \; \qf \leq -\xi.
\label{eq:follower-nash-i}
\end{align}

\paragraph{Case (ii):} $0 < \xi + \qf < \ka$. Note that $y$ is given by~\eqref{eq:follower-optimal-production-ii}. Since $0 < \xi + \qf < \ka \implies 0 < \frac{1}{\na+1}\left(\xi + \qf\right) < \ka$, the symmetric follower reactions are given by:
\begin{align}
& \frac{1}{2}\left(\xi - (\na-1)\qf\right) = \frac{1}{\na+1}\left(\xi + \qf\right) \; \mathrm{and} \; \xi > \max((\na-1)\qf,-(\na+1)\qf) \; \mathrm{and} \; \eqref{eq:follower-nash-range-2} \; \mathrm{holds} 
\nonumber\\
& \mathrm{or} \; \xi + \qf = \frac{1}{\na+1}\left(\xi + \qf\right) \; \mathrm{and} \; -2\qf \leq \xi \leq -(\na+1)\qf \; \mathrm{and} \; \eqref{eq:follower-nash-range-2} \; \mathrm{holds} 
\nonumber\\
& \mathrm{or} \; \frac{1}{2}\xi = \frac{1}{\na+1}\left(\xi + \qf\right) \; \mathrm{and} \; \xi < \min(2\ka, -2\qf) \; \mathrm{and} \; \eqref{eq:follower-nash-range-2} \; \mathrm{holds} 
\nonumber\\
\iff & 
\qf = \frac{\na-1}{\na^2+1}\xi \; \mathrm{and} \; \xi > \max((\na-1)\qf,-(\na+1)\qf) \; \mathrm{and} \; \eqref{eq:follower-nash-range-2} \; \mathrm{holds} 
\nonumber\\
& \mathrm{or} \; \qf = -\xi \; \mathrm{and} \; -2\qf \leq \xi \leq -(\na+1)\qf \; \mathrm{and} \; \eqref{eq:follower-nash-range-2} \; \mathrm{holds} 
\nonumber\\
& \mathrm{or} \; \qf = \frac{\na-1}{2}\xi \; \mathrm{and} \; \xi < \min(2\ka, -2\qf) \; \mathrm{and} \; \eqref{eq:follower-nash-range-2} \; \mathrm{holds}
\nonumber\\
\iff &
\qf = \frac{\na-1}{\na^2+1}\xi \; \mathrm{and} \; \xi > \max\left(\frac{(\na-1)^2}{\na^2+1}\xi,-\frac{\na^2-1}{\na^2+1}\xi\right) \; \mathrm{and} \; 0 < \xi < \frac{\na^2+1}{\na(\na+1)}\ka 
\nonumber\\
& \mathrm{or} \; \qf = \frac{\na-1}{2}\xi \; \mathrm{and} \; \xi < \min(2\ka, -(\na-1)\xi) \; \mathrm{and} \; 0 < \frac{\na+1}{2}\xi < \ka
\nonumber\\
\iff & 
\qf = \frac{\na-1}{\na^2+1}\xi \; \mathrm{and} \; 0 < \xi < \frac{\na^2+1}{\na(\na+1)}\ka.
\label{eq:follower-nash-ii}
\end{align}
The second equivalence is due to the fact that $\qf = -\xi \implies \xi + \qf = 0$. The third equivalence is due to the fact that $\xi > 0 \implies \frac{(\na-1)^2}{\na^2+1}\xi \geq -\frac{\na^2-1}{\na^2+1}\xi$ and $\frac{\na+1}{2}\xi > 0 \implies 2\ka > -(\na-1)\xi$. 

\paragraph{Case (iii):} $\ka \leq \xi + \qf \leq \na\ka$. Note that $y$ is given by~\eqref{eq:follower-optimal-production-iii}. Therefore, the symmetric follower reactions are given by:
\begin{align}
&\left[\frac{1}{2}\left(\xi - (\na-1)\qf\right)\right]_0^\ka = \left[\frac{1}{\na+1}\left(\xi + \qf\right)\right]_0^\ka \; \mathrm{and} \; \eqref{eq:follower-nash-range-3} \; \mathrm{holds}
\nonumber\\ 
\iff &
\frac{1}{2}\left(\xi - (\na-1)\qf\right) = \frac{1}{\na+1}\left(\xi + \qf\right) \; \mathrm{and} \; \eqref{eq:follower-nash-range-3} \; \mathrm{holds} \; \mathrm{and} \; 0 < \frac{1}{2}\left(\xi - (\na-1)\qf\right) < \ka
\nonumber\\
\iff &
\qf = \frac{\na-1}{\na^2+1}\xi \; \mathrm{and} \; \eqref{eq:follower-nash-range-3} \; \mathrm{holds} \; \mathrm{and} \; 0 < \frac{1}{2}\left(\xi - (\na-1)\qf\right) < \ka
\nonumber\\
\iff &
\qf = \frac{\na-1}{\na^2+1}\xi \; \mathrm{and} \; \frac{\na^2+1}{\na(\na+1)}\ka \leq \xi \leq \frac{\na^2+1}{\na+1}\ka \; \mathrm{and} \; 0 < \xi < \frac{\na^2+1}{\na-\frac{1}{2}}\ka
\nonumber\\
\iff &
\qf = \frac{\na-1}{\na^2+1}\xi \; \mathrm{and} \; \frac{\na^2+1}{\na(\na+1)}\ka \leq \xi \leq \frac{\na^2+1}{\na+1}\ka.
\label{eq:follower-nash-iii}
\end{align}
The first equivalence is due to the fact that $\eqref{eq:follower-nash-range-3} \implies 0 < \frac{1}{\na+1}\left(\xi + \qf\right) < \ka$. The last equivalence is due to the fact that $\frac{\na^2+1}{\na(\na+1)} < \frac{\na^2+1}{\na+1} < \frac{\na^2+1}{\na-\frac{1}{2}}$.

\paragraph{Case (iv):} $\na\ka < \xi + \qf < (\na+1)\ka$. Note that $y$ is described by~\eqref{eq:follower-optimal-production-iv}. We show that there does not exist a symmetric follower reaction such that $y = z_1$. Suppose otherwise. By Proposition~\ref{prop:SpotNash}, for each $j\neq l$,
\begin{align*}
y_j
&=
\left[\frac{1}{\na}\left(\xi + \qf - \left[\frac{1}{2}\left(\xi -(\na-1)\ka\right)\right]_0^{\xi + \qf - \na\ka}\right)\right]_0^\ka
\\
&=
\left[-\frac{1}{\na}\left[\frac{1}{2}\left(-\xi - 2\qf -(\na-1)\ka\right)\right]_{-\xi-\qf}^{- \na\ka}\right]_0^\ka
\\
&=
\left[\left[\frac{1}{2\na}\left(\xi + 2\qf +  (\na-1)\ka\right)\right]_{\ka}^{\frac{1}{\na}\left(\xi+\qf\right)}\right]_0^\ka
\\
&=
\ka
\end{align*}
However, $\eqref{eq:follower-nash-range-4} \implies \frac{1}{\na+1}\left(\xi + \qf\right) < \ka \implies y < \ka = y_j$, which contradicts with the fact that a symmetric follower reaction implies symmetric productions (by Proposition~\ref{prop:SpotEquilibrium2}).

Therefore, the symmetric follower reactions are given by:
\begin{align}
&\left[\frac{1}{2}\left(\xi - (\na-1)\qf\right)\right]_{\xi + \qf - \na\ka}^{\ka} = \frac{1}{\na+1}\left(\xi + \qf\right) \; \mathrm{and} \; \eqref{eq:follower-nash-range-4} \; \mathrm{holds} \; \mathrm{and} \; \hat{\phi}_l(z_1;\qf,\qbp) \leq \hat{\phi}_l(z_2;\qf,\qbp)
\nonumber\\
\iff &
\frac{1}{2}\left(\xi - (\na-1)\qf\right) = \frac{1}{\na+1}\left(\xi + \qf\right) \; \mathrm{and} \; \eqref{eq:follower-nash-range-4} \; \mathrm{holds} \; \mathrm{and} \; \hat{\phi}_l(z_1;\qf,\qbp) \leq \hat{\phi}_l(z_2;\qf,\qbp)
\nonumber\\
\iff & 
\qf = \frac{\na-1}{\na^2+1}\xi \; \mathrm{and} \; \eqref{eq:follower-nash-range-4} \; \mathrm{holds} 
\; \mathrm{and} \; \hat{\phi}_l(z_1;\qf,\qbp) \leq \hat{\phi}_l(z_2;\qf,\qbp)
\nonumber\\
\iff & 
\qf = \frac{\na-1}{\na^2+1}\xi \; \mathrm{and} \; \frac{\na^2+1}{\na+1}\ka < \xi < \frac{\na^2+1}{\na}\ka 
\; \mathrm{and} \; \hat{\phi}_l(z_1;\qf,\qbp) \leq \hat{\phi}_l(z_2;\qf,\qbp)
\nonumber\\
\iff &
\qf = \frac{\na-1}{\na^2+1}\xi \; \mathrm{and} \; \frac{\na^2+1}{\na+1}\ka < \xi \leq \frac{(\na^2+1)(\na-1)}{\na^2-2\sqrt{\na}+1}\ka.
\label{eq:follower-nash-iv}
\end{align}
The first equivalence follows from the fact that $\xi + \qf - \na\ka = \frac{1}{\na+1}\left(\xi + \qf\right) \implies \xi + \qf = (\na+1)\ka$ and $\ka = \frac{1}{\na+1}\left(\xi + \qf\right) \implies \xi + \qf = (\na+1)\ka$. The last equivalence follows from the following facts. First, note that
\begin{align*}
z_1 
&= \left[\frac{1}{2}\left(\xi - (\na-1)\ka\right)\right]_0^{\frac{\na(\na+1)}{\na^2+1}\xi - \na\ka} 
\\
&= 
\begin{cases}
\frac{\na(\na+1)}{\na^2+1}\xi - \na\ka, & \mathrm{if} \; \frac{\na^2+1}{\na+1}\ka < \xi < \frac{\na^2+1}{\na}\ka \; \mathrm{and} \; \frac{1}{2}\left(\xi - (\na-1)\ka\right) > \frac{\na(\na+1)}{\na^2+1}\xi - \na\ka,
\\
\frac{1}{2}\left(\xi - (\na-1)\ka\right), & \mathrm{if} \; \frac{\na^2+1}{\na+1}\ka < \xi < \frac{\na^2+1}{\na}\ka \; \mathrm{and} \; \frac{1}{2}\left(\xi - (\na-1)\ka\right) \leq \frac{\na(\na+1)}{\na^2+1}\xi - \na\ka,
\end{cases}
\\
&= 
\begin{cases}
\frac{\na(\na+1)}{\na^2+1}\xi - \na\ka, & \mathrm{if} \; \frac{\na^2+1}{\na+1}\ka < \xi < \frac{(\na+1)(\na^2+1)}{\na^2+2\na-1}\ka,
\\
\frac{1}{2}\left(\xi - (\na-1)\ka\right), & \frac{(\na+1)(\na^2+1)}{\na^2+2\na-1}\ka \leq \xi < \frac{\na^2+1}{\na}\ka.
\end{cases}
\end{align*}
where the second equality is due to $\xi > \frac{\na^2+1}{\na+1}\ka > \frac{\na^2-1}{\na+1}\ka = (\na-1)\ka$. Thus, if $\frac{\na^2+1}{\na+1}\ka < \xi < \frac{(\na+1)(\na^2+1)}{\na^2+2\na-1}\ka$, then
\begin{align*}
&\hat{\phi}_l(z_1;\qf,\qbp) \leq \hat{\phi}_l(z_2;\qf,\qbp)
\\
\iff & \left(\xi - (\na-1)\ka - z_1\right)z_1 \leq \frac{1}{\na}\left(\xi - (\na-1)\qf - z_2\right)z_2
\\
\iff & \left(\ka - \qf\right)\left(\xi + \qf - \na\ka\right) \leq \frac{1}{4\na}\left(\xi - (\na-1)\qf\right)^2
\\
\iff & \mathrm{True}.
\end{align*}
On the other hand, if $\frac{(\na+1)(\na^2+1)}{\na^2+2\na-1}\ka \leq \xi < \frac{\na^2+1}{\na}\ka$, then
\begin{align*}
&\hat{\phi}_l(z_1;\qf,\qbp) \leq \hat{\phi}_l(z_2;\qf,\qbp)
\\
\iff & \left(\xi - (\na-1)\ka - z_1\right)z_1 \leq \frac{1}{\na}\left(\xi - (\na-1)\qf - z_2\right)z_2
\\
\iff & \frac{1}{2}\left(\xi - (\na-1)\ka\right)^2 \leq \frac{1}{4\na}\left(\xi - (\na-1)\qf\right)^2
\\
\iff & \xi \leq \frac{(\na^2+1)(\na-1)}{\na^2-2\sqrt{\na}+1}\ka,
\end{align*}
where $\frac{(\na+1)(\na^2+1)}{\na^2+2\na-1}\ka \leq \frac{(\na^2+1)(\na-1)}{\na^2-2\sqrt{\na}+1}\ka < \frac{\na^2+1}{\na}\ka$.

\paragraph{Case (v):} $(\na+1)\ka \leq \xi + \qf$. Note that $y$ is given by~\eqref{eq:follower-optimal-production-v}. Therefore, the symmetric follower reactions are given by:
\begin{align}
&\left[\frac{1}{2}\left(\xi - (\na-1)\ka\right)\right]_0^\ka=\left[\frac{1}{\na+1}\left(\xi + \qf\right)\right]_0^\ka \; \mathrm{and} \; \eqref{eq:follower-nash-range-5} \; \mathrm{holds}
\nonumber\\
\iff &
\left[\frac{1}{2}\left(\xi - (\na-1)\ka\right)\right]_0^\ka = \ka \; \mathrm{and} \; \eqref{eq:follower-nash-range-5} \; \mathrm{holds}
\nonumber\\
\iff &
\qf \geq -\xi + (\na+1)\ka \; \mathrm{and} \; \xi \geq (\na+1)\ka.
\label{eq:follower-nash-v}
\end{align}

Putting together the descriptions in~\eqref{eq:follower-nash-i}~--~\eqref{eq:follower-nash-v} yield~\eqref{eq:follower-reaction}.
\endproof

\subsection{Leader Reaction Analyses}
\label{app:leader}

\begin{proposition}
\label{prop:BaseloadReaction}
Fix the followers' forward positions $\qfp = \qf\mathbf{1} \in \mathbb{R}^\na$. Let $\Qb\subseteq\mathbb{R}_+$ denote the set of symmetric leader reactions, i.e., for each $\qb\in\Qb$ and $i\in\nb$, 
\begin{align}
\pfpjf\left(x;x\mathbf{1},\qf\mathbf{1}\right)
\geq
\pfpjf\left(\bar{x};x\mathbf{1},\qf\mathbf{1}\right),
\quad
\forall 
\bar{x}\in\mathbb{R}_+.
\end{align}
Let:
\begin{align*}
    x_1 &=  \frac{1}{\nb+1}\left[\rv_x\right]_0^\infty,
    \\
    x_2 &= \frac{1}{\nb}\left(\rv_x - \cm + \qf\right),
    \\
    x_3 &= \frac{1}{\nb+1}\left[\rv_x + \na\left(\cm - \qf\right)\right]_0^\infty,
    \\
    x_4 &= \frac{1}{\nb+1}\left[\rv_x - \na\ka\right]_0^\infty.
\end{align*}
Then,
\begin{align*}
    \Qb = \left\{ 
        x\in\mathbb{R}_+\left|
        \begin{array}{l}
            x = x_1 \; \mathit{if} \; \qf - \cm < -\frac{\rv_x}{\nb+1},
            \\
            \mathit{or} \; x = x_2 \; \mathit{if} \; -\frac{\rv_x}{\nb+1} \leq \qf - \cm \leq \min\left(-\frac{\rv_x}{\nb\na+\nb+1}, \eta_4\right),
            \\
            \mathit{or} \; x = x_3 \; \mathit{if} \; -\frac{\rv_x}{\na\nb+\nb+1} < \qf - \cm \leq \max\left(\eta_3, \ka - (\rv_x - \na\ka)\right),
            \\
            \mathit{or} \; x = x_4 \; \mathit{if} \; \qf - \cm \geq 
\begin{cases}
\ka - (\rv_x - \na\ka), & \mathit{if} \; \rv_x < \na\ka,
\\
\eta_2, & \mathit{if} \; \rv_x \geq \left(1+\frac{(\nb+1)\sqrt{\na+1}}{(\sqrt{\na+1}-1)^2}\right)\na\ka,
\\
\eta_1, & \mathit{otherwise}.
\end{cases}
        \end{array}
        \right.
    \right\},
\end{align*}
where
\begin{align*}
\eta_1
& :=
\begin{array}{l}
\ka - \frac{\rv_x - \na\ka}{\na}\left(\frac{2(\sqrt{\na+1}-1)}{\nb+1}\right),
\end{array}
\\
\eta_2
& :=
\begin{array}{l}
- \frac{1}{2}\left(\frac{2(\rvb-\na\ka)}{\nb+1}+\na\ka\right)\left(1-\sqrt{1-\left(\frac{\frac{2(\rvb-\na\ka)}{\nb+1}}{\frac{2(\rvb-\na\ka)}{\nb+1}+\na\ka}\right)^2}\right),
\end{array}
\\
\eta_3
& :=
\begin{array}{l}
\ka - \frac{\rv_x - \na\ka}{\na}\left(\frac{2(\sqrt{\na+1}-1)}{2+(\nb-1)\sqrt{\na+1}}\right),
\end{array}
\\
\eta_4
& :=
\begin{array}{l}
- \frac{1}{2}\left(\frac{2(\rvb-\nb\na\ka)}{\nb+1}+\left(\frac{2\nb}{\nb+1}\right)^2\na\ka\right)\left(1-\sqrt{1-\left(\frac{\frac{2(\rvb-\nb\na\ka)}{\nb+1}}{\frac{2(\rvb-\nb\na\ka)}{\nb+1}+\left(\frac{2\nb}{\nb+1}\right)^2\na\ka}\right)^2}\right).
\end{array}
\end{align*}
Moreover, for each $x\in \Qb$,
\begin{align*}
y_j(\qf\mathbf{1},x\mathbf{1}) = 0 
& \iff x = x_1 \; \mathit{or} \; x_2,
\\
0 < y_j(\qf\mathbf{1},x\mathbf{1}) < \ka 
& \iff x = x_3,
\\
y_j(\qf\mathbf{1},x\mathbf{1}) = \ka 
& \iff x = x_4.
\end{align*}
\end{proposition}

    \proof{Proof.}
The proof proceeds in three steps. In step 1, we solve for a leader's payoff maximizing production quantity given that all other leaders produce equal quantities. In step 2, we solve for the symmetric leader productions that satisfy the condition that every leader is producing at its payoff maximizing quantity. The latter gives the set of symmetric leader reactions. In step 3, we explain how the solutions obtained in step 2 is equivalent to $\Qb$.

\emph{Step 1:} Fix a leader $l\in\nb$ and suppose $x_i = x$ for every $i\neq l$. We solve for the solution to
\begin{align}
\sup_{x_l \in\mathbb{R}_+} \psi_l\left(x_l;x\mathbf{1},\qf\mathbf{1}\right).
\label{eq:leader-nash-aux}
\end{align}
Substituting for the demand function yields
\begin{align*}
\psi_l\left(x_l;x\mathbf{1},\qf\mathbf{1}\right)
&=
\beta\left(\rv_x - x_l - (\nb-1)x - \sum_{j=1}^{\na} y_j(\qf\mathbf{1},x_l,x\mathbf{1})\right)x_l
\end{align*}
where the follower productions are given by
\begin{align*}
&\sum_{j=1}^{\na} y_j(\qf\mathbf{1},x_l,x\mathbf{1})
\\
&=
\na\left[\frac{1}{\na+1}\left(\rv_y + \qf - x_l - (\nb-1)x\right)\right]_0^\ka
\\
&=
\begin{cases}
0, & \mathrm{if} \; \eqref{eq:leader-nash-range-1} \; \mathrm{holds},
\\
\begin{cases}
0, & \mathrm{if} \; \rv_y + \qf - x_l - (\nb-1)x \leq 0,
\\
\frac{\na}{\na+1}\left(\rv_y + \qf - x_l - (\nb-1)x\right), & \mathrm{otherwise},
\end{cases}
 & \mathrm{if} \; \eqref{eq:leader-nash-range-2} \; \mathrm{holds},
\\
\begin{cases}
0, & \mathrm{if} \; \rv_y + \qf - x_l - (\nb-1)x \leq 0,
\\
\ka, & \mathrm{if} \; \rv_y + \qf - x_l - (\nb-1)x \geq (\na+1)\ka,
\\
\frac{\na}{\na+1}\left(\rv_y + \qf - x_l - (\nb-1)x\right), & \mathrm{otherwise},
\end{cases}
 & \mathrm{if} \; \eqref{eq:leader-nash-range-3} \; \mathrm{holds},
\end{cases}
\end{align*}
where the second equality is due to the fact that $x_l, x \geq 0$ and the three cases~\eqref{eq:leader-nash-range-1}~--~\eqref{eq:leader-nash-range-3} are defined by
\begin{subequations}
\begin{align}
\rv_y + \qf - (\nb-1)x & \leq 0, 
\label{eq:leader-nash-range-1}
\\
0 < \rv_y + \qf - (\nb-1)x & \leq (\na+1)\ka,
\label{eq:leader-nash-range-2}
\\
(\na+1)\ka < \rv_y + \qf - (\nb-1)x &. 
\label{eq:leader-nash-range-3}
\end{align}
\end{subequations}
We analyze each case separately.

\paragraph{Case (i):} $\rv_y + \qf  - (\nb-1)x  < 0$. We obtain
\begin{align*}
\psi_l\left(x_l;x\mathbf{1},f\mathbf{1}\right)
=
\beta\left(\rv_x - x_l - (\nb-1)x\right) x_l.
\end{align*}
Therefore, $\psi_l(x_l;x\mathbf{1},f\mathbf{1})$ is a smooth function in $x_l$ over $\mathbb{R}_+$. The first and second derivatives are given by
\begin{align*}
\frac{\partial}{\partial x_l}\psi_l\left(x_l;x\mathbf{1},f\mathbf{1}\right)
&=
\beta\left(\rv_x - (\nb-1)x - 2x_l\right),
\\
\frac{\partial^2}{\partial x_l^2}\psi_l\left(x_l;x\mathbf{1},f\mathbf{1}\right)
&=
-2\beta < 0,
\end{align*}
which implies that $\psi_l\left(x_l;x\mathbf{1},f\mathbf{1}\right)$ is strictly concave in $x_l$. Therefore, $x_l$ is a solution to~\eqref{eq:leader-nash-aux} if and only if it satisfies the following first order optimality conditions:
\begin{align}
\frac{\partial^+}{\partial x_l}\psi_l\left(x_l;x\mathbf{1},f\mathbf{1}\right) & \leq 0, \quad \mathrm{if} \; 0 \leq x_l,
\label{eq:leader-nash-foc-1}
\\
\frac{\partial^-}{\partial x_l}\psi_l\left(x_l;x\mathbf{1},f\mathbf{1}\right) & \geq 0, \quad \mathrm{if} \; 0 < x_l.
\label{eq:leader-nash-foc-2}
\end{align}
It is straightforward to show that there is a unique solution is given by
\begin{align}
x_l
=
\left[\frac{1}{2}\left(\rv_x - (\nb-1)x\right)\right]_0^\infty.
\label{eq:leader-optimal-production-i}
\end{align}

\paragraph{Case (ii):} $0 \leq \rv_y + \qf  - (\nb-1)x < (\na+1)\ka$. We obtain
\begin{align*}
\psi_l\left(x_l;x\mathbf{1},f\mathbf{1}\right)
=
\begin{cases}
\beta\left(\rv_x - x_l - (\nb-1)x\right)x_l, & \mathrm{if} \; x_l \geq \rv_y + \qf - (\nb-1)x,
\\
\beta\left(\frac{1}{\na+1}\rv_x + \frac{\na}{\na+1}\left(\cm - \qf\right) - \frac{\nb-1}{\na+1}x - \frac{1}{\na+1}x_l \right)x_l, & \mathrm{otherwise}.
\end{cases}
\end{align*}
Therefore, $\psi_l\left(x_l;x\mathbf{1},f\mathbf{1}\right)$ is a piecewise smooth function in $x_l$ over $\mathbb{R}_+$. The first and second derivatives are given by
\begin{align*}
\frac{\partial}{\partial x_l}\psi_l\left(x_l;x\mathbf{1},f\mathbf{1}\right)
&=
\begin{cases}
\beta\left(\rv_x - (\nb-1)x - 2x_l\right), & \mathrm{if} \; x_l > \rv_y + \qf - (\nb-1)x,
\\
\frac{\beta}{\na+1}\left(\rv_x + \na\left(\cm - \qf\right) - (\nb-1)x - 2 x_l\right), & \mathrm{otherwise},
\end{cases}
\\
\frac{\partial^2}{\partial x_l^2}\psi_l\left(x_l;x\mathbf{1},f\mathbf{1}\right)
&=
\begin{cases}
-2\beta, & \mathrm{if} \; x_l > \rv_y + \qf - (\nb-1)x,
\\
-\frac{2}{\na+1}\beta, & \mathrm{otherwise},
\end{cases}
\\
& < 0.
\end{align*}
Moreover, we have
\begin{align*}
&\left.\frac{\partial^-}{\partial x_l}\psi_l\left(x_l;x\mathbf{1},f\mathbf{1}\right)\right|_{x_l = \rv_y + \qf - (\nb-1)x}
\\
&=
\beta\left(\frac{1}{\na+1}\rv_x + \frac{\na}{\na+1}\left(\cm + \qf\right) - \frac{\nb-1}{\na+1}x - \frac{2}{\na+1}\left(\rv_y + \qf - (\nb-1)x\right)\right)
\\
&=
\beta\left(-\rv_x + 2\left(\cm - \qf\right) + \frac{\na}{\na+1}\left(\rv_y + \qf\right) + \frac{\nb-1}{\na+1}x \right)
\\
&\geq
\beta\left(-\rv_x + 2\left(\cm - \qf\right) + \left(\nb-1\right)x \right)
\\
&=
\left.\frac{\partial^+}{\partial x_l}\psi_l\left(x_l;x\mathbf{1},f\mathbf{1}\right)\right|_{x_l = \rv_y + \qf - (\nb-1)x},
\end{align*}
where the inequality follows from~\eqref{eq:leader-nash-range-2}. Therefore, $\psi_l\left(x_l;x\mathbf{1},\qf\mathbf{1}\right)$ is concave in $x_l$ over $\mathbb{R}_+$. Therefore, $x_l$ is a solution to~\eqref{eq:leader-nash-aux} if and only if it satisfies the first order optimality conditions~\eqref{eq:leader-nash-foc-1}~--~\eqref{eq:leader-nash-foc-2}. It is straightforward to show that there is a unique solution given by
\begin{align}
x_l
=
\begin{cases}
0, & \mathrm{if} \; \eqref{eq:leader-nash-range-2a} \; \mathrm{holds},
\\
\frac{1}{2}\left(\rv_x + \na\left(\cm - \qf\right) - \left(\nb-1\right)x\right), & \mathrm{if} \; \eqref{eq:leader-nash-range-2b} \; \mathrm{holds},
\\
\rv_x - \cm + \qf - (\nb-1)x, & \mathrm{if} \; \eqref{eq:leader-nash-range-2c} \; \mathrm{holds},
\\
\frac{1}{2}\left(\rv_x - (\nb-1)x\right), & \mathrm{if} \; \eqref{eq:leader-nash-range-2d} \; \mathrm{holds},
\end{cases}
\label{eq:leader-optimal-production-ii}
\end{align}
where the cases~\eqref{eq:leader-nash-range-2a}~--~\eqref{eq:leader-nash-range-2d} are defined by:
\begin{subequations}
\begin{align}
\rv_x + \na\left(\cm - \qf\right)  \leq (\nb-1)x &,
\label{eq:leader-nash-range-2a}
\\
(\nb-1)x &< \min\left(\rv_x + \na(\cm - \qf), \rv_x - (\na+2)(\cm - \qf)\right),
\label{eq:leader-nash-range-2b}
\\
\rv_x - (\na+2)(\cm - \qf) \leq (\nb-1)x &\leq \rv_x - 2(\cm - \qf),
\label{eq:leader-nash-range-2c}
\\
\rv_x - 2(\cm - \qf) < (\nb-1)x &.
\label{eq:leader-nash-range-2d}
\end{align}
\end{subequations}

\paragraph{Case (iii):} $(\na+1)\ka \leq \rv_y + \qf - (\nb-1)x$. We obtain
\begin{align*}
&\psi_l\left(x_l;x\mathbf{1},f\mathbf{1}\right)
\\
&=
\begin{cases}
\beta\left(\rv_x - x_l - (\nb-1)x\right)x_l, & \mathrm{if} \; x_l \geq \rv_y + \qf - (\nb-1)x,
\\
\beta\left(\rv_x - x_l -(\nb-1)x - \na\ka\right)x_l, & \mathrm{if} \; x_l \leq \rv_y + \qf - (\nb-1)x - (\na+1)\ka, 
\\
\beta\left(\frac{1}{\na+1}\rv_x + \frac{\na}{\na+1}\left(\cm - \qf\right) - \frac{\nb-1}{\na+1}x - \frac{1}{\na+1}x_l \right)x_l, & \mathrm{otherwise}.
\end{cases}
\end{align*}
Therefore, $\phi_l\left(x_l;x\mathbf{1},\qf\mathbf{1}\right)$ is a piecewise smooth function in $x_l$ over $\mathbb{R}_+$. The first and second derivatives are given by
\begin{align*}
&\frac{\partial}{\partial x_l}\psi_l\left(x_l;x\mathbf{1},f\mathbf{1}\right)
\\
&=
\begin{cases}
\beta\left(\rv_x - (\nb-1)x - 2x_l\right), & \mathrm{if} \; x_l > \rv_y + \qf - (\nb-1)x,
\\
\beta\left(\rv_x - (\nb-1)x - \na\ka - 2x_l\right), & \mathrm{if} \; x_l < \rv_y + \qf - (\nb-1)x - (\na+1)\ka,
\\
\frac{\beta}{\na+1}\left(\rv_x + \na\left(\cm - \qf\right) - (\nb-1)x - 2 x_l\right), & \mathrm{otherwise},
\end{cases}
\\
&\frac{\partial^2}{\partial x_l^2}\psi_l\left(x_l;x\mathbf{1},f\mathbf{1}\right)
\\
&=
\begin{cases}
-2\beta, & \mathrm{if} \; x_l > \rv_y + \qf - (\nb-1)x,
\\
-2\beta, & \mathrm{if} \; x_l < \rv_y + \qf - (\nb-1)x - (\na+1\ka,
\\
-\frac{2}{\na+1}\beta, & \mathrm{otherwise},
\end{cases}
\\
& < 0.
\end{align*}
Moreover, we have
\begin{align*}
&\left.\frac{\partial^-}{\partial x_l}\psi_l\left(x_l;x\mathbf{1},f\mathbf{1}\right)\right|_{x_l = \rv_y + \qf - (\nb-1)x}
\\
&=
\beta\left(\frac{1}{\na+1}\rv_x + \frac{\na}{\na+1}\left(\cm - \qf\right) - \frac{\nb-1}{\na+1}x - \frac{2}{\na+1}\left(\rv_y + \qf - (\nb-1)x\right)\right)
\\
&=
\beta\left(-\rv_x + 2(\cm-\qf) + \frac{\na}{\na+1}(\rv_x -\cm + \qf)+ \frac{\nb-1}{\na+1}x \right)
\\
&>
\beta\left(-\rv_x + 2(\cm - \qf) + (\nb-1)x \right)
\\
&=
\left.\frac{\partial^+}{\partial x_l}\psi_l\left(x_l;x\mathbf{1},f\mathbf{1}\right)\right|_{x_l = \rv_y + \qf - (\nb-1)x}.
\end{align*}
Therefore, $\phi_l\left(x_l;x\mathbf{1},\qf\mathbf{1}\right)$ is concave in $x_l$ over $\left[\rv_y + \qf - (\nb-1)x - (\na+1)\ka,\infty\right)$. However, it is straightforward to check that $\phi_l\left(x_l;x\mathbf{1},\qf\mathbf{1}\right)$ has a non-concave kink at $x_l = \rv_y + \qf - (\nb-1)x - (\na+1)\ka$, and therefore $\phi_l\left(x_l;x\mathbf{1},\qf\mathbf{1}\right)$ is not concave in $x_l$ over $\mathbb{R}_+$. Therefore, solve the following sub-problems:
\begin{align}
\sup_{x_l\in\left[0,\rv_y + \qf - (\nb-1)x - (\na+1)\ka\right]} \psi_l\left(x_l;x\mathbf{1},\qf\mathbf{1}\right),
\label{eq:leader-nash-iii-prob1}
\end{align}
and
\begin{align}
\sup_{x_l\in\left[\rv_y + \qf - (\nb-1)x - (\na+1)\ka,\infty\right)} \psi_l\left(x_l;x\mathbf{1},\qf\mathbf{1}\right).
\label{eq:leader-nash-iii-prob2}
\end{align}
The solution of the sub-problem with the larger optimal value is the solution to~\eqref{eq:leader-nash-aux}. Using the first order optimality conditions, the unique solution to~\eqref{eq:leader-nash-iii-prob1} is given by
\begin{align*}
x_l
=
\left[\frac{1}{2}\left(\rv_x - (\nb-1)x - \na\ka\right)\right]_0^{\rv_y + \qf - (\nb-1)x - (\na+1)\ka} =: z_1,
\end{align*}
and that to~\eqref{eq:leader-nash-iii-prob2} is given by
\begin{align*}
x_l
=
\begin{cases}
\rv_y + \qf - (\nb-1)x - (\na+1)\ka, & \mathrm{if} \; \eqref{eq:leader-nash-range-3a} \; \mathrm{holds},
\\
\frac{1}{2}\left(\rv_x + \na(\cm - \qf) - (\nb-1)x\right), & \mathrm{if} \; \eqref{eq:leader-nash-range-3b} \; \mathrm{holds},
\\
\rv_y + \qf - (\nb-1)x, & \mathrm{if} \; \eqref{eq:leader-nash-range-3c} \; \mathrm{holds},
\\
\frac{1}{2}\left(\rv_x - (\nb-1)x\right), & \mathrm{if} \; \eqref{eq:leader-nash-range-3d} \; \mathrm{holds},
\end{cases}
=: z_2.
\end{align*}
where the cases~\eqref{eq:leader-nash-range-3a}~--~\eqref{eq:leader-nash-range-3d} are defined by:
\begin{subequations}
\begin{align}
(\nb-1)x  &\leq \rv_x - (\na+2)(\cm - \qf) - 2(\na+1)\ka,
\label{eq:leader-nash-range-3a}
\\
\rv_x - (\na+2)(\cm - \qf) - 2(\na+1)\ka < (\nb-1)x  &< \rv_x - (\na+2)(\cm - \qf),
\label{eq:leader-nash-range-3b}
\\
\rv_x - (\na+2)(\cm - \qf) \leq (\nb-1)x  &\leq \rv_x - 2(\cm - \qf),
\label{eq:leader-nash-range-3c}
\\
\rv_x - 2(\cm - \qf) < (\nb-1)x &.
\label{eq:leader-nash-range-3d}
\end{align}
\end{subequations}
Therefore, the solution(s) to~\eqref{eq:leader-nash-aux} are given by:
\begin{subequations}
\begin{align}
x_l = z_1, & \quad \mathrm{if} \; \psi_l\left(z_1;x\mathbf{1},\qf\mathbf{1}\right) > \psi_l\left(z_2;x\mathbf{1},\qf\mathbf{1}\right),
\\
x_l = z_2, & \quad \mathrm{if} \; \phi_l\left(z_2;x\mathbf{1},\qf\mathbf{1}\right) > \phi_l\left(z_1;x\mathbf{1},\qf\mathbf{1}\right),
\\
x_l = z_1 \; \mathrm{or} \; z_2, & \quad \mathrm{if} \; \psi_l\left(z_1;x\mathbf{1},\qf\mathbf{1}\right) = \psi_l\left(z_2;x\mathbf{1},\qf\mathbf{1}\right).
\end{align}
\label{eq:leader-optimal-production-iii}
\end{subequations}

\emph{Step 2:} We solve for the symmetric leader productions that satisfy the condition that every leader is producing at its payoff maximizing quantity. We divide the analyses into three cases depending on the value of $\rv_y + \qf - (\nb-1)x$. 

\paragraph{Case (i):} $\rv_y + \qf - (\nb-1)x < 0$. Note that $x_l$ is given by~\eqref{eq:leader-optimal-production-i}. Therefore, the symmetric leader reactions are given by:
\begin{align*}
& x = 0 \; \mathrm{and} \; \rv_x < (\nb-1)x \; \mathrm{and} \; \eqref{eq:leader-nash-range-1} \; \mathrm{holds}
\\
& \mathrm{or} \; x = \frac{1}{2}\left(\rv_x - (\nb-1)x\right) \; \mathrm{and} \; \rv_x \geq (\nb-1)x \; \mathrm{and} \; \eqref{eq:leader-nash-range-1} \; \mathrm{holds}
\\
\iff
& x = 0 \; \mathrm{and} \; \rv_x < 0 \; \mathrm{and} \; \eqref{eq:leader-nash-range-1} \; \mathrm{holds}
\\ 
& \mathrm{or} \; x = \frac{1}{\nb+1}\rv_x \; \mathrm{and} \; \rv_x \geq 0 \; \mathrm{and} \; \eqref{eq:leader-nash-range-1} \; \mathrm{holds}.
\end{align*}

\paragraph{Case (ii):} $0 \leq \rv_y + \qf - (\nb-1)x < (\na+1)\ka$. Note that $x_l$ is given by~\eqref{eq:leader-optimal-production-ii}. Therefore, the symmetric leader reactions are given by:
\begin{align*}
& x = 0 \; \mathrm{and} \; \eqref{eq:leader-nash-range-2a} \; \mathrm{and} \; \eqref{eq:leader-nash-range-2} \; \mathrm{holds}
\\
& \mathrm{or} \; x = \frac{1}{2}\left((\na+1)\rv_x - \na(\rv_y + \qf) - (\nb-1)x\right) \; \mathrm{and} \; \eqref{eq:leader-nash-range-2b}  \; \mathrm{and} \; \eqref{eq:leader-nash-range-2} \; \mathrm{holds}
\\
& \mathrm{or} \; x = \rv_y + \qf - (\nb-1)x \; \mathrm{and} \; \eqref{eq:leader-nash-range-2c} \; \mathrm{and} \; \eqref{eq:leader-nash-range-2} \; \mathrm{holds}
\\
& \mathrm{or} \; x = \frac{1}{2}\left(\rv_x - (\nb-1)x\right) \; \mathrm{and} \; \eqref{eq:leader-nash-range-2d}  \; \mathrm{and} \; \eqref{eq:leader-nash-range-2} \; \mathrm{holds}
\\
\iff & 
x = 0 \; \mathrm{and} \; \frac{\rv_x}{\na} \leq \qf - \cm \; \mathrm{and} \; \eqref{eq:leader-nash-range-2} \; \mathrm{holds}
\\
& \mathrm{or} \; x = \frac{1}{\nb+1}\left(\rv_x + \na(\cm-\qf)\right) \; \mathrm{and} \; -\frac{\rv_x}{\na\nb + \nb + 1} < \qf - \cm < \frac{\rv_x}{\na}  \; \mathrm{and} \; \eqref{eq:leader-nash-range-2} \; \mathrm{holds}
\\
& \mathrm{or} \; x = \frac{1}{\nb}\left(\rv_x - \left(\cm - \qf\right)\right) \; \mathrm{and} \; -\frac{\rv_x}{\nb+1} \leq \qf - \cm \leq -\frac{\rv_x}{\na\nb + \nb + 1} \; \mathrm{and} \; \eqref{eq:leader-nash-range-2} \; \mathrm{holds}
\\
& \mathrm{or} \; x = \frac{1}{\nb+1}\rv_x \; \mathrm{and} \; \qf - \cm < -\frac{\rv_x}{\nb+1} \; \mathrm{and} \; \eqref{eq:leader-nash-range-2} \; \mathrm{holds}.
\end{align*}

\paragraph{Case (iii):} $(\na+1)\ka \leq \rv_y + \qf - (\nb-1)x$. Note that $x_l$ is described by~\eqref{eq:leader-optimal-production-iii}. Therefore, the symmetric leader reactions are given by
\begin{align}
& x = z_1 \; \mathrm{and} \; \eqref{eq:leader-nash-range-3} \; \mathrm{holds} \; \mathrm{and} \; \psi_l\left(z_1;x\mathbf{1},\qf\mathbf{1}\right) > \psi_l\left(z_2;x\mathbf{1},\qf\mathbf{1}\right)
\label{eq:leader-nash-case-iiia}
\\
\mathrm{or} \; & x = z_2 \; \mathrm{and} \; \eqref{eq:leader-nash-range-3} \; \mathrm{holds} \; \mathrm{and} \; \psi_l\left(z_1;x\mathbf{1},\qf\mathbf{1}\right) < \psi_l\left(z_2;x\mathbf{1},\qf\mathbf{1}\right)
\label{eq:leader-nash-case-iiib}
\\
\mathrm{or} \; & x = z_1 \;\mathrm{or} \; z_2 \; \mathrm{and} \; \eqref{eq:leader-nash-range-3} \; \mathrm{holds} \; \mathrm{and} \; \psi_l\left(z_1;x\mathbf{1},\qf\mathbf{1}\right) = \psi_l\left(z_2;x\mathbf{1},\qf\mathbf{1}\right).
\label{eq:leader-nash-case-iiic}
\end{align}
We analyze the cases $x = z_1$ and $x = z_2$ separately. 

Suppose $x = z_1$ is a symmetric leader reaction. Since
\begin{align*}
&\left.\frac{\partial^-}{\partial x_l}\psi_l(x_l;x\mathbf{1},\qf\mathbf{1})\right|_{x_l = \rv_y + \qf - (\nb-1)x - (\na+1)\ka}
\\
&=
\beta\left(-\rv_x + 2(\cm - \qf) + (\nb-1)x + (\na+2)\ka\right)
\\
&\leq
\frac{\beta}{\na+1}\left(-\rv_x + (\na+2)(\cm - \qf) + (\nb-1)x + 2 (\na+1)\ka\right)
\\
&=
\left.\frac{\partial^+}{\partial x_l}\psi_l(x_l;x\mathbf{1},\qf\mathbf{1})\right|_{x_l = \rv_y + \qf - (\nb-1)x - (\na+1)\ka},
\end{align*}
we infer that $x < \rv_y + \qf - (\nb-1)x - (\na+1)\ka$. Therefore, we obtain
\begin{align*}
& x = z_1
\\
\iff & x = 0 \; \mathrm{and} \; \eqref{eq:leader-nash-range-3} \; \mathrm{holds} \; \mathrm{and} \; \psi_l\left(z_1;x\mathbf{1},\qf\mathbf{1}\right) \geq \psi_l\left(z_2;x\mathbf{1},\qf\mathbf{1}\right) \; \mathrm{and} \; \rv_x - (\nb-1)x - \na\ka \leq 0
\\
& \mathrm{or} \; x = \frac{1}{2}\left(\rv_x - (\nb-1)x - \na\ka\right) \; \mathrm{and} \; \eqref{eq:leader-nash-range-3} \; \mathrm{holds} \; \mathrm{and} \; \psi_l\left(z_1;x\mathbf{1},\qf\mathbf{1}\right) \geq \psi_l\left(z_2;x\mathbf{1},\qf\mathbf{1}\right) \\
& \quad \; \mathrm{and} \; 0 < \frac{1}{2}\left(\rv_x - (\nb-1)x - \na\ka\right) < \rv_y + \qf - (\nb-1)x - (\na+1)\ka
\\
\iff & x = 0 \; \mathrm{and} \; \eqref{eq:leader-nash-range-3} \; \mathrm{holds} \; \mathrm{and} \; \rv_x - \na\ka \leq 0
\\
& \mathrm{or} \; x = \frac{1}{\nb+1}\left(\rv_x - \na\ka\right) \; \mathrm{and} \; \eqref{eq:leader-nash-range-3} \; \mathrm{holds} \; \mathrm{and} \; \psi_l\left(z_1;x\mathbf{1},\qf\mathbf{1}\right) \geq \psi_l\left(z_2;x\mathbf{1},\qf\mathbf{1}\right) \\
& \quad \; \mathrm{and} \; \rv_x - \na\ka > 0 \; \mathrm{and} \; \qf - \cm > -\frac{1}{\nb+1}\left(\rv_x - \na\ka\right) + \ka
\end{align*}
The second equivalence follows from solving for $x$ in the equations, and the fact that in the case $x = 0$, the inequalities $\eqref{eq:leader-nash-range-3} \; \mathrm{and} \; \rv_x - \na\ka \leq 0 \implies \psi_l\left(z_1;x\mathbf{1},\qf\mathbf{1}\right) \geq \psi_l\left(z_2;x\mathbf{1},\qf\mathbf{1}\right)$. 

Suppose $x = z_2$. Then, using the same arguments in~\eqref{eq:leader-nash-case-iiia}, we infer that $x < \rv_y + \qf - (\nb-1)x - (\na+1)\ka$. Therefore, we obtain
\begin{align*}
& x = z_2
\\
\iff 
& x = \frac{1}{2}\left(\rv_x + \na(\cm - \qf) - (\nb-1)x\right) \; \mathrm{and} \; \eqref{eq:leader-nash-range-3} \; \mathrm{and} \; \eqref{eq:leader-nash-range-3b} \; \mathrm{holds} \\
& \quad \; \mathrm{and} \; \psi_l\left(z_1;x\mathbf{1},\qf\mathbf{1}\right) \leq \psi_l\left(z_2;x\mathbf{1},\qf\mathbf{1}\right)
\\
& \mathrm{or} \; x = \rv_x + (\qf - \cm) - (\nb-1)x \; \mathrm{and} \; \eqref{eq:leader-nash-range-3} \; \mathrm{and} \; \eqref{eq:leader-nash-range-3c} \; \mathrm{holds} \\
& \quad \; \mathrm{and} \; \psi_l\left(z_1;x\mathbf{1},\qf\mathbf{1}\right) \leq \psi_l\left(z_2;x\mathbf{1},\qf\mathbf{1}\right)
\\
& \mathrm{or} \; x = \frac{1}{2}\left(\rv_x - (\nb-1)x\right) \; \mathrm{and} \; \eqref{eq:leader-nash-range-3} \; \mathrm{and} \; \eqref{eq:leader-nash-range-3d} \; \mathrm{holds} \\
& \quad \; \mathrm{and} \; \psi_l\left(z_1;x\mathbf{1},\qf\mathbf{1}\right) \leq \psi_l\left(z_2;x\mathbf{1},\qf\mathbf{1}\right)
\\
\iff 
& x = \frac{1}{\nb+1}\left(\rv_x + \na(\cm - \qf)\right) \; \mathrm{and} \; \eqref{eq:leader-nash-range-3} \; \mathrm{holds} \; \mathrm{and} \; \psi_l\left(z_1;x\mathbf{1},\qf\mathbf{1}\right) \leq \psi_l\left(z_2;x\mathbf{1},\qf\mathbf{1}\right) \\
& \quad \; \mathrm{and} \; -\frac{1}{\na\nb+\nb+1}\rv_x < \qf - \cm < \frac{1}{\na\nb + \nb + 1}\left(-\rv_x + (\nb+1)(\na+1)\ka\right)
\\
& \mathrm{or} \; x = \frac{1}{\nb}\left(\rv_x + (\qf - \cm)\right) \; \mathrm{and} \; \eqref{eq:leader-nash-range-3} \; \mathrm{holds} \; \mathrm{and} \; \psi_l\left(z_1;x\mathbf{1},\qf\mathbf{1}\right) \leq \psi_l\left(z_2;x\mathbf{1},\qf\mathbf{1}\right) \\
& \quad \; \mathrm{and} \; -\frac{1}{\nb+1}\rv_x \leq \qf - \cm \leq - \frac{1}{\na\nb + \nb + 1}\rv_x
\\
& \mathrm{or} \; x = \frac{1}{\nb+1}\rv_x \; \mathrm{and} \; \eqref{eq:leader-nash-range-3} \; \mathrm{holds} \; \mathrm{and} \; \mathrm{and} \; \qf - \cm < -\frac{1}{\nb+1}\rv_x
\end{align*}
The second equivalence follows from solving for $x$ in the equations, and the fact that in the case $x = \frac{1}{\nb+1}\rv_x$, the inequalities $\eqref{eq:leader-nash-range-3} \; \mathrm{and} \; \qf-\cm < -\frac{1}{\nb+1}\rv_x \implies \psi_l\left(z_1;x\mathbf{1},\qf\mathbf{1}\right) \leq \psi_l\left(z_2;x\mathbf{1},\qf\mathbf{1}\right)$.

\emph{Step 3:} We explain how the solutions obtained in step 2 is equivalent to $\Qb$. Observe that step 2 obtains five cases for $x$:
\begin{align*}
x &= \frac{1}{\nb+1}\rv_x,
\\
x &= \frac{1}{\nb}\left(\rv_x - (\cm - \qf)\right),
\\
x &= \frac{1}{\nb+1}\left(\rv_x + \na(\cm - \qf)\right),
\\
x &= \frac{1}{\nb+1}\left(\rv_x - \na\ka\right),
\\
x &= 0.
\end{align*}
We analyze each case separately.

\paragraph{Case (i):} $x = \frac{1}{\nb+1}\rv_x$. This case is characterized by
\begin{align*}
& \rv_x \geq 0 \; \mathrm{and} \; \eqref{eq:leader-nash-range-1} \; \mathrm{holds}
\\
& \mathrm{or} \; \qf - \cm < - \frac{\rv_x}{\nb+1} \; \mathrm{and} \; \eqref{eq:leader-nash-range-2} \; \mathrm{holds}
\\
& \mathrm{or} \; \qf - \cm < -\frac{\rv_x}{\nb+1} \; \mathrm{and} \; \eqref{eq:leader-nash-range-3} \; \mathrm{holds}
\\ \iff
& \rv_x \geq 0 \; \mathrm{and} \; \qf - \cm < -\frac{\rv_x}{\nb+1}.
\end{align*}
The equivalence is due to the following facts. First, $\eqref{eq:leader-nash-range-2} \implies \rv_x \geq 0$ and $\eqref{eq:leader-nash-range-3} \implies \rv_x \geq 0$. Second, $\eqref{eq:leader-nash-range-1} \;\mathrm{and} \; \rv_x \geq 0 \implies \qf - \cm < -\frac{1}{\nb+1}\rv_x$. Third, $\eqref{eq:leader-nash-range-1},~\eqref{eq:leader-nash-range-2},~\eqref{eq:leader-nash-range-3} \implies \mathrm{True}$.

\paragraph{Case (ii):} $x = \frac{1}{\nb}\left(\rv_x - (\cm - \qf)\right)$. This case is characterized by
\begin{align*}
&-\frac{\rv_x}{\nb+1} \leq \qf - \cm \leq - \frac{\rv_x}{\na\nb+\nb+1} \; \mathrm{and} \; \eqref{eq:leader-nash-range-2} \; \mathrm{holds}
\\
&\mathrm{or} \; -\frac{\rv_x}{\nb+1} \leq \qf - \cm \leq -\frac{\rv_x}{\na\nb+\nb+1} 
\\ 
&\quad \; \mathrm{and} \; \psi_l\left(z_1;x\mathbf{1},\qf\mathbf{1}\right) \leq \psi_l\left(z_2;x\mathbf{1},\qf\mathbf{1}\right) \; \mathrm{and} \; \eqref{eq:leader-nash-range-3} \; \mathrm{holds}
\\ \iff
&-\frac{\rv_x}{\nb+1} \leq \qf - \cm \leq \min\left( - \frac{\rv_x}{\na\nb+\nb+1}, -\rv_x + \nb(\na+1)\ka\right)
\\
&\mathrm{or} \; -\frac{\rv_x}{\nb+1} \leq \qf - \cm \leq -\frac{\rv_x}{\na\nb+\nb+1} 
\\ & \quad \; \mathrm{and} \; \qf - \cm \leq \eta_4 \; \mathrm{and} \; \qf - \cm > - \rv_x + \nb(\na+1)\ka
\\ \iff
& -\frac{\rv_x}{\nb+1} \leq \qf - \cm \leq \min\left(-\frac{\rv_x}{\nb\na+\nb+1}, \eta_4\right).
\end{align*}
The first equivalence is due to the following facts. First, $\eqref{eq:leader-nash-range-2} \iff -\rv_x < \qf - \cm \leq -\rv_x + \nb(\na+1)\ka$. Second, $\psi_l\left(z_1;x\mathbf{1},\qf\mathbf{1}\right) \leq \psi_l\left(z_2;x\mathbf{1},\qf\mathbf{1}\right) \iff \qf - \cm \leq \eta_4$. Third, $\eqref{eq:leader-nash-range-3} \iff \qf - \cm > - \rv_x + \nb(\na+1)\ka$. The second equivalence is due to the fact that $-\frac{\rv_x}{\nb\na + \nb + 1} \leq -\rv_x + \nb(\na+1)\ka \implies \eta_4 > -\frac{\rv_x}{\nb\na+\nb+1}$.

\paragraph{Case (iii):} $x = \frac{1}{\nb+1}\left(\rv_x + \na(\cm-\qf)\right)$. This case is characterized by
\begin{align*}
& -\frac{\rv_x}{\na\nb+\nb+1} < \qf - \cm < \frac{\rv_x}{\na} \; \mathrm{and} \; \eqref{eq:leader-nash-range-2} \; \mathrm{holds}
\\
& \mathrm{or} \; -\frac{\rv_x}{\na\nb+\nb+1} < \qf - \cm < \frac{1}{\na\nb+\nb+1}\left(-\rv_x + (\nb+1)(\na+1)\ka\right) \\
& \quad \; \mathrm{and} \; \psi_l\left(z_1;x\mathbf{1},\qf\mathbf{1}\right) \leq \psi_l\left(z_2;x\mathbf{1},\qf\mathbf{1}\right) \; \mathrm{and} \; \eqref{eq:leader-nash-range-3} \; \mathrm{holds}
\\ \iff
& -\frac{\rv_x}{\na\nb+\nb+1} < \qf - \cm < \frac{\rv_x}{\na} \; \mathrm{and} \; \qf - \cm \leq \frac{-2\rv_x + (\nb+1)(\na+1)\ka}{\nb+1+\na(\nb-1)}
\\
& \mathrm{or} \; -\frac{\rv_x}{\na\nb+\nb+1} < \qf - \cm < \frac{1}{\na\nb+\nb+1}\left(-\rv_x + (\nb+1)(\na+1)\ka\right) \\
& \quad \; \mathrm{and} \; \qf - \cm \leq \eta_3 \; \mathrm{and} \; \qf - \cm > \frac{-2\rv_x + (\nb+1)(\na+1)\ka}{\nb+1+\na(\nb-1)}
\\ \iff
& \rv_x < \na\ka \; \mathrm{and} \; -\frac{\rv_x}{\na\nb+\nb+1} < \qf - \cm < \frac{\rv_x}{\na}
\\
& \mathrm{or} \; \rv_x \geq \na\ka \; \mathrm{and} \; -\frac{\rv_x}{\na\nb+\nb+1} < \qf - \cm \leq \eta_3.
\end{align*}
The first equivalence is due to the following facts. First, $\eqref{eq:leader-nash-range-2} \iff -\frac{2\rv_x}{\nb+1+\na(\nb-1)} < \qf - \cm \leq \frac{-2\rv_x + (\nb+1)(\na+1)\ka}{\nb+1+\na(\nb-1)}$. Second, $\eqref{eq:leader-nash-range-3} \iff \qf - \cm > \frac{-2\rv_x + (\nb+1)(\na+1)\ka}{\nb+1+\na(\nb-1)}$. Third, $\psi_l\left(z_1;x\mathbf{1},\qf\mathbf{1}\right) \leq \psi_l\left(z_2;x\mathbf{1},\qf\mathbf{1}\right) \iff \qf - \cm \leq \eta_3$. The second equivalence is due to the following facts. First, $\frac{\rv_x}{\na} \leq \frac{-2\rv_x + (\nb+1)(\na+1)\ka}{\nb+1+\na(\nb-1)} \iff \rv_x \leq \na\ka$. Second, $\rv_x \geq \na\ka \implies \eta_3 < \frac{1}{\na\nb+\nb+1}\left(-\rv_x + (\nb+1)(\na+1)\ka\right)$.

\paragraph{Case (iv):} $x = \frac{1}{\nb+1}\left(\rv_x - \na\ka\right)$. This case is characterized by
\begin{align*}
& \rv_x - \na\ka > 0 \; \mathrm{and} \; \qf - \cm > -\frac{1}{\nb+1}\left(\rv_x - \na\ka\right) + \ka \\
&\quad \; \mathrm{and} \; \psi_l\left(z_1;x\mathbf{1},\qf\mathbf{1}\right) \geq \psi_l\left(z_2;x\mathbf{1},\qf\mathbf{1}\right) \; \mathrm{and} \; \eqref{eq:leader-nash-range-3} \; \mathrm{holds}
\\ \iff
& \rv_x - \na\ka > 0 \; \mathrm{and} \; \qf - \cm > -\frac{1}{\nb+1}\left(\rv_x - \na\ka\right) + \ka \; \mathrm{and} \; \psi_l\left(z_1;x\mathbf{1},\qf\mathbf{1}\right) \geq \psi_l\left(z_2;x\mathbf{1},\qf\mathbf{1}\right)
\\ \iff
& \na\ka < \rv_x < \left(1 + \frac{(\nb+1)\sqrt{\na+1}}{(\sqrt{\na+1}-1)^2}\right)\na\ka \; \mathrm{and} \; \qf - \cm \geq \eta_1
\\
& \mathrm{or} \; \left(1 + \frac{(\nb+1)\sqrt{\na+1}}{(\sqrt{\na+1}-1)^2}\right)\na\ka \leq \rv_x \; \mathrm{and} \; \qf - \cm \geq \eta_2.
\end{align*}
The first equivalence is due to the fact that $\eqref{eq:leader-nash-range-3} \implies \qf - \cm > -\frac{1}{\nb+1}(\rv_x - \na\ka) + \ka$. The second equivalence is due to the fact that $\psi_l\left(z_1;x\mathbf{1},\qf\mathbf{1}\right) \geq \psi_l\left(z_2;x\mathbf{1},\qf\mathbf{1}\right) \iff \na\ka < \rv_x < \left(1 + \frac{(\nb+1)\sqrt{\na+1}}{(\sqrt{\na+1}-1)^2}\right)\na\ka \; \mathrm{and} \; \qf - \cm \geq \eta_1 \; \mathrm{or} \; \left(1 + \frac{(\nb+1)\sqrt{\na+1}}{(\sqrt{\na+1}-1)^2}\right)\na\ka \leq \rv_x \; \mathrm{and} \; \qf - \cm \geq \eta_2$.

\paragraph{Case (v):} $x = 0$. This case is characterized by
\begin{align*}
& \rv_x < 0 \; \mathrm{and} \; \eqref{eq:leader-nash-range-1}
\\ 
& \mathrm{or} \; \frac{\rv_x}{\na} \leq \qf - \cm \; \mathrm{and} \; \eqref{eq:leader-nash-range-2}
\\
& \mathrm{or} \; \rv_x - \na\ka \leq 0 \; \mathrm{and} \; \eqref{eq:leader-nash-range-3}
\\ \iff 
& \rv_x < 0 \; 
\\
&\mathrm{or} \; \rv_x \geq 0 \; \mathrm{and} \; \rv_x + \na(\cm - \qf) \leq 0 \; \mathrm{and} \; 0 \leq \rv_x + (\qf - \cm) < (\na+1)\ka
\\
&\mathrm{or} \; \rv_x \geq 0 \; \mathrm{and} \; \rv_x - \na\ka \leq 0 \; \mathrm{and} \; (\na+1)\ka \leq \rv_x + (\qf - \cm).
\end{align*}
The equivalence is due to the following facts. First, $\eqref{eq:leader-nash-range-2} \;\mathrm{and} \; \rv_x < 0 \implies \rv_x + \na(\cm - \qf) < 0$. Second, $\eqref{eq:leader-nash-range-3} \; \mathrm{and} \; \rv_x < 0 \implies \rv_x - \na\ka < 0$. 
\endproof

\subsection{Forward Market Equilibrium}
\label{app:market}

\begin{theorem}
    \label{prop:MarketEquilibrium}
Suppose $\rv_x > 0$. Let $Q \subseteq\mathbb{R}\times\mathbb{R}_+$ denote the set of all symmetric Nash equilibria, i.e., $(\qf,\qb) \in Q$ if $(\qf\mathbf{1},\qb\mathbf{1})$ is a Nash equilibrium of the forward market. Let:
\begin{align*}
    Q_1 &:= \left\{(\qf,\qb)\in\mathbb{R}\times\mathbb{R}_+ \left| 
    \begin{array}{l}
        \qb = \frac{1}{\nb+1}\rv_x \\
        \qf < \cm - \frac{1}{\nb+1}\rv_x
    \end{array} \right.\right\},    
    \\
    Q_2 &:= \left\{(\qf,\qb) \in\mathbb{R}\times\mathbb{R}_+ \left| 
    \begin{array}{l}
        \qb = \frac{1}{\nb}\left(\rv_x - (\cm - \qf)\right) \\ 
            \max\left(0, \cm - \frac{\rv_x}{\nb+1} \right)\leq \qf \leq \cm + \min\left(-\frac{\rv_x}{\nb\na+\nb+1}, \eta_4\right)
    \end{array} \right.\right\}, 
    \\
    Q_3 &:= \left\{(\qf,\qb) \in\mathbb{R}\times\mathbb{R}_+ \left| 
    \begin{array}{l}
        \qb = \frac{\na+1}{\na^2+\nb\na+\nb+1}\left(\rv_x + \na^2\cm\right)
        \\
        \qf = \frac{\na-1}{\na^2+\nb\na+\nb+1}\left(\rv_x - (\nb\na+\nb+1)\cm\right)
    \end{array} \right.\right\},
    \\
    Q_4 &:= \left\{(\qf,\qb)\in\mathbb{R}\times\mathbb{R}_+  \left| 
    \begin{array}{l}
        \qb = \frac{1}{\nb+1}\left(\rv_x - \na\ka\right)
        \\
        \qf \geq \cm + 
        \begin{cases}
            \eta_1, & \mathit{if} \; \na\ka < \rv_x \leq \left(1 + \frac{(\nb+1)\sqrt{\na+1}}{(\sqrt{\na+1}-1)^2}\right)\na\ka,
            \\
            \eta_2, & \mathit{if} \; \left(1 + \frac{(\nb+1)\sqrt{\na+1}}{(\sqrt{\na+1}-1)^2}\right)\na\ka < \rv_x.
        \end{cases} 
    \end{array}\right.\right\}.
\end{align*}
where $\eta_1,\eta_2,\eta_4$ are as defined in Proposition~\ref{prop:BaseloadReaction}. Then,
\begin{align*}
    Q=
    \left\{(\qf,\qb)\in\mathbb{R}\times\mathbb{R}_+\left|
    \begin{array}{l}
        (\qf,\qb)\in Q_1 \; \mathit{if} \; \rv_x \leq (\nb+1)\cm,
        \\
        \mathit{or} \; (\qf,\qb)\in Q_2 \; \mathit{if} \; \rv_x \leq \min\left( (\nb\na+\nb+1)\cm, \zeta_1 \right),
        \\
        \mathit{or} \; (\qf,\qb)\in Q_3 \; \mathit{if} \; (\nb\na+\nb+1)\cm < \rv_x \leq \zeta_2,
        \\
        \mathit{or} \; (\qf,\qb)\in Q_4 \; \mathit{if} \; (\nb+1)(\cm+\ka) + \na\ka \leq \rv_x.
    \end{array}
    \right.\right\},
\end{align*}
where
\begin{align}
    \zeta_1 &= \nb\na\ka + \left( \nb+1  \right)\cm + 2\nb\sqrt{\na\ka\cm},
    \\
    \zeta_2 &= \left(\nb\na+\nb+1\right)\cm + \frac{\na^2+\na\nb+\nb+1}{\na(\na+1)+2(1-\sqrt{\na+1})}\left( \na\ka-(\sqrt{\na+1}-1)^2\cm \right)
    \label{}
\end{align}
Moreover, for each $(\qf,\qb)\in Q$, 
\begin{align*}
    y_j(\qf\mathbf{1},\qb\mathbf{1}) = 0 & \iff (\qf,\qb)\in Q_1\cup Q_2,
    \\
    0 < y_j(\qf\mathbf{1},\qb\mathbf{1}) < k & \iff (\qf,\qb)\in Q_3,
    \\
    y_j(\qf\mathbf{1},\qb\mathbf{1}) = \ka & \iff (\qf,\qb) \in Q_4.
\end{align*}
\end{theorem}

    \proof{Proof.}
The symmetric equilibria are given by the intersection of the follower and leader reactions obtained in Propositions~\ref{prop:PeakerReaction} and~\ref{prop:BaseloadReaction}. We divide the analyses into three separate cases depending on the value of the follower productions $y_j(\qf\mathbf{1},\qb\mathbf{1})$. 

\emph{Case (i):} $y_j(\qf\mathbf{1},\qb\mathbf{1}) = 0$. Using Propositions~\ref{prop:PeakerReaction} and~\ref{prop:BaseloadReaction}, we infer that $(\qf,\qb)$ is a symmetric equilibrium with $y_j(\qf\mathbf{1},\qb\mathbf{1})=0$ if and only if
\begin{subequations}
    \label{eq:nash-y0}
\begin{align}
\qf &\leq -\left(\rv_x - \cm - \nb\qb\right),
\label{eq:nash-y0-a}
\\
0 & \geq \rv_x - \cm - \nb\qb,
\label{eq:nash-y0-b}
\\
x &= 
\begin{cases}
\frac{1}{\nb+1}\left[\rv_x\right]_0^\infty, & \mathrm{if} \; \qf - \cm < -\frac{\rv_x}{\nb+1},
\\
\frac{1}{\nb}\left(\rv_x - \cm + \qf\right), & \mathrm{if} \; -\frac{\rv_x}{\nb+1} \leq \qf - \cm \leq \min\left(-\frac{\rv_x}{\nb\na+\nb+1},\eta_4\right).
\end{cases}
\label{eq:nash-y0-c}
\end{align}
\end{subequations}

Suppose $x = \frac{1}{\nb+1}\left[\rv_x\right]_0^\infty$. Since $\rv_x > 0$, we infer that $x = \frac{1}{\nb+1}\rv_x$. Substituting into~\eqref{eq:nash-y0-a} and~\eqref{eq:nash-y0-b} yields
\begin{align*}
\eqref{eq:nash-y0-a} & \iff \qf < \cm - \frac{\rv_x}{\nb+1},
\\
\eqref{eq:nash-y0-b} & \iff \rv_x \leq (\nb+1)\cm.
\end{align*}
The above inequalities, together with~\eqref{eq:nash-y0-c}, imply that $(\qf,\qb)$ satisfies~\eqref{eq:nash-y0} with $x = \frac{1}{\nb+1}\rv_x$, if and only if  $(\qf,\qb)\in Q_1$ and $\rv_x \leq (\nb+1)\cm$. 

Suppose $x = \frac{1}{\nb}\left(\rv_x - \cm + \qf\right)$. Substituting into~\eqref{eq:nash-y0-a} and~\eqref{eq:nash-y0-b} yields
\begin{align*}
\eqref{eq:nash-y0-a} & \iff \qf \leq \qf \iff \mathrm{True},
\\
\eqref{eq:nash-y0-b} & \iff \qf \geq 0.
\end{align*}
Therefore, there exists $(\qf,\qb)$ that satisfies~\eqref{eq:nash-y0} with $\qb = \frac{1}{\nb}(\rv_x-\cm+\qf)$ if and only if
\begin{align*}
    & \left[\max\left( 0,\cm -\frac{\rv_x}{\nb+1}\right),\cm + \min\left( -\frac{\rv_x}{\nb\na+\nb+1},\eta_4 \right)  \right]\neq\varnothing
    \\
    \iff & 
    0 \leq \cm + \min\left( -\frac{\rv_x}{\nb\na+\nb+1},\eta_4 \right)
    \\
    \iff & 
    \rv_x \leq \left( \nb\na+\nb+1 \right)\cm \; \mathrm{and} \; 0 \leq \cm + \eta_4
    \\
    \iff & 
    \rv_x \leq \left( \nb\na+\nb+1 \right)\cm \; \mathrm{and} \; \rv_x \leq \zeta_1.
\end{align*}
Therefore, $(\qf,\qb)$ satisfies~\eqref{eq:nash-y0} with $x = \frac{1}{\nb}\left(\rv_x - \cm + \qf\right)$, if and only if $(\qf,\qb)\in Q_2$ and $\rv_x \leq \min\left( (\nb\na+\nb+1)\cm,\zeta_1 \right)$.

\emph{Case (ii):} $0 \leq y_j(\qf\mathbf{1},\qb\mathbf{1}) \leq \ka$. Using Propositions~\ref{prop:PeakerReaction} and~\ref{prop:BaseloadReaction}, we infer that $(\qf,\qb)$ is a symmetric equilibrium if and only if
\begin{subequations}
    \label{eq:nash-y}
\begin{align}
& \qf = \frac{\na-1}{\na^2+1}\left(\rv_x-\cm-\nb\qb\right),
\label{eq:nash-y-a}
\\
& 0 \leq \rv_x - \cm - \nb\qb \leq \xi_1,
\label{eq:nash-y-b}
\\
& \qb = \frac{1}{\nb+1}\left[\rv_x+\na(\cm-\qf)\right]_0^\infty,
\label{eq:nash-y-c}
\\
& -\frac{\rv_x}{\nb\na+\nb+1} < \qf - \cm \leq \max\left(\ka - (\rv_x - \na\ka),\eta_3\right).
\label{eq:nash-y-d}
\end{align}
\end{subequations}
We show that $x > 0$. Suppose otherwise. Substituting into~\eqref{eq:nash-y-a} implies that $\qf = \frac{\na-1}{\na^2+1}(\rv_x - \cm)$. Substituting further into~\eqref{eq:nash-y-c} yields
\begin{align*}
\rv_x + \na\left(\cm - \frac{\na-1}{\na^2+1}(\rv_x - \cm)\right) \leq 0
\iff
\rv_x + \na\cm < 0,
\end{align*}
which is a contradiction since $\rv_x > 0$, $\cm \geq 0$, and $\na \geq 2$. Therefore, we assume that $x > 0$. Solving~\eqref{eq:nash-y-a} and~\eqref{eq:nash-y-c} gives
\begin{align*}
\qf
&=
\frac{\na-1}{\na^2+\nb\na+\nb+1}\left(\rv_x -(\nb\na+\nb+1)\cm\right),
\\
\qb
&=
\frac{\na+1}{\na^2+\nb\na+\nb+1}\left(\rv_x + \na^2\cm\right).
\end{align*}
Substituting for $\qb$ yields
\begin{align*}
\eqref{eq:nash-y-b}
\iff
(\nb\na+\nb+1)\cm \leq \rv_x \leq (\nb\na+\nb+1)\cm + \frac{(\na^2+\nb\na+\nb+1)(\na-1)}{\na^2-2\sqrt{\na}+1}\ka.
\end{align*}
Substituting for $\qf$ yields
\begin{align*}
\eqref{eq:nash-y-d}
& \iff
(\nb\na+\nb+1)\cm < \rv_x \; \mathrm{and} \; \rv_x \leq 
\begin{cases}
\frac{\na(\nb+1)}{\nb+\na}\cm + \left(\na + \frac{\nb+1}{\nb+\na}\right)\ka, & \mathrm{if} \; \rv_x \leq \na\ka,
\\
\zeta_2, & \mathrm{if} \; \rv_x > \na\ka,
\end{cases}
\\
& \iff
(\nb\na+\nb+1)\cm < \rv_x \; \mathrm{and} \; \rv_x \leq 
\begin{cases}
\na\ka, & \mathrm{if} \; \rv_x \leq \na\ka,
\\
\zeta_2, & \mathrm{if} \; \rv_x > \na\ka,
\end{cases}
\\
& \iff
(\nb\na+\nb+1)\cm < \rv_x \leq \zeta_2.
\end{align*}
The first equivalence is due to the fact that $\ka-(\rv_x-\na\ka) \geq \eta_3 \iff \rv_x \leq \na\ka$. The second equivalence is due to the fact that $\cm \geq 0$ and $\ka > 0$. Next, using the fact that $\na \geq 2, \cm \geq 0, \ka > 0$, we obtain
\begin{align*}
    \zeta_2 < (\nb\na+\nb+1)\cm + \frac{(\na^2+\nb\na+\nb+1)(\na-1)}{\na^2-2\sqrt{\na}+1}\ka,
\end{align*}
from which it follows that $(\qf,\qb)$ satisfies~\eqref{eq:nash-y} if and only if $(\qf,\qb)\in Q_3$ and $(\nb\na+\nb+1)\cm \leq \rv_x \leq \zeta_2$.

\emph{Case (iii):} $y_j(\qf\mathbf{1},\qb\mathbf{1}) = \ka$. From Propositions~\ref{prop:PeakerReaction} and~\ref{prop:BaseloadReaction}, we infer that $(\qf,\qb)$ is a symmetric equilibrium if and only if
\begin{subequations}
    \label{eq:nash-yk}
\begin{align}
    & \qf \geq -(\rv_x - \cm - \nb\qb)+(\na+1)\ka,
    \label{eq:nash-yk-a}
    \\
    & (\na+1)\ka \leq \rv_x - \cm - \nb\qb,
    \label{eq:nash-yk-b}
    \\
    & x = \frac{1}{\nb+1}\left[\rv_x - \na\ka\right]_0^\infty,
    \label{eq:nash-yk-c}
    \\
    & \qf - \cm \geq 
    \begin{cases}
        \ka - (\rv_x - \na\ka), & \mathrm{if} \; \rv_x < \na\ka,
        \\
        \eta_2, & \mathrm{if} \; \rv_x \geq \left(1 + \frac{(\nb+1)\sqrt{\na+1}}{(\sqrt{\na+1}-1)^2}\right)\na\ka,
        \\
        \eta_1, & \mathrm{otherwise}.
    \end{cases}
    \label{eq:nash-yk-d}
\end{align}
\end{subequations}
We divide the analyses into three cases depending on the value of $\rv_x$.

Suppose $0 < \rv_x \leq \na\ka$. Then,~\eqref{eq:nash-yk-c} implies $x = 0$. However, substituting into~\eqref{eq:nash-yk-b} implies that $\rv_x - \na\ka \geq \ka + \cm > 0$ which is a contradiction. Therefore, there does not exist an equilibrium such that $0 < \rv_x \leq \na\ka$. 

Suppose $\na\ka < \rv_x \leq \left(1+\frac{(\nb+1)\sqrt{\na+1}}{(\sqrt{\na+1}-1)^2}\right)\na\ka$. Then,~\eqref{eq:nash-yk-c} implies $x = \frac{1}{\nb+1}\left(\rv_x - \na\ka\right)$. Substituting for $x$ yields
\begin{align*}
    \eqref{eq:nash-yk-a} & \iff \qf \geq \cm - \frac{1}{\nb+1}\left(\rv_x - \na\ka\right)+\ka,
    \\
    \eqref{eq:nash-yk-b} & \iff \cm - \frac{1}{\nb+1}\left(\rv_x - \na\ka\right) + \ka \leq 0.
\end{align*}
From~\eqref{eq:nash-yk-d}, we infer that $\qf \geq \cm + \eta_1$. Since $\na \geq 2 \implies \eta_1 \geq \ka - \frac{\rv_x - \na\ka}{\nb+1}$, it follows that the symmetric equilibria are characterized by
\begin{align}
    & x = \frac{1}{\nb+1}(\rv_x - \na\ka) \; \mathrm{and} \; 
    \qf \geq \cm + \eta_1 \ \mathrm{and} \; 
    \cm - \frac{1}{\nb+1}(\rv_x - \na\ka) + \ka \leq 0.
    \label{eq:nash-yk-sol1}
\end{align}

Suppose $\left(1+\frac{(\nb+1)\sqrt{\na+1}}{(\sqrt{\na+1}-1)^2}\right)\na\ka < \rv_x$. Then, we again have
\begin{align*}
    \eqref{eq:nash-yk-a} & \iff \qf \geq \cm - \frac{1}{\nb+1}\left(\rv_x - \na\ka\right)+\ka,
    \\
    \eqref{eq:nash-yk-b} & \iff \cm - \frac{1}{\nb+1}\left(\rv_x - \na\ka\right) + \ka \leq 0.
\end{align*}
From~\eqref{eq:nash-yk-d}, we infer that $\qf \geq \cm + \eta_2$. Since $\na \geq 2 \implies \eta_2 \geq \ka - \frac{\rv_x - \na\ka}{\nb+1}$, it follows that the symmetric equilibria are characterized by
\begin{align}
    & x = \frac{1}{\nb+1}(\rv_x - \na\ka) \; \mathrm{and} \; 
    \qf \geq \cm + \eta_2 \ \mathrm{and} \; 
    \cm - \frac{1}{\nb+1}(\rv_x - \na\ka) + \ka \leq 0.
    \label{eq:nash-yk-sol2}
\end{align}

By combining the characterizations in~\eqref{eq:nash-yk-sol1} and~\eqref{eq:nash-yk-sol2}, we infer that $(\qf,\qb)$ satisfies~\eqref{eq:nash-yk} if and only if $(\qf,\qb)\in Q_4$ and $(\nb+1)(\cm+\ka) + \na\ka\leq \rv_x$.
\endproof

\subsection{Stackelberg Equilibrium}

\begin{theorem}
    \label{prop:StackelbergEquilibrium}
    Suppose followers' forward positions $\mathbf{f} = \mathbf{0}$. Let $\Qb \subseteq\mathbb{R}_+$ denote the set of symmetric leader reactions, i.e., for each $x\in\Qb$ and $i\in\nb$,
    \begin{align*}
        \psi_i(\qb;\qb\mathbf{1},\mathbf{0})\geq\psi_i(\bar{\qb};\qb\mathbf{1},\mathbf{0}),\quad\forall\bar{x}\in\mathbb{R}_+.
    \end{align*}
    Let:
    \begin{align*}
        \qb_1 &= \frac{1}{\nb+1}\left[\rv_x\right]_0^\infty,
        \\
        \qb_2 &= \frac{1}{\nb}\left( \rv_x-\cm \right),
        \\
        \qb_3 &= \frac{1}{\nb+1}\left[ \rv_x+\na\cm \right]_0^\infty,
        \\
        \qb_4 &= \frac{1}{\nb+1}\left[ \rv_x-\na\ka \right]_0^\infty.
    \end{align*}
    Then,
    \begin{align*}
        \Qb = 
        \left\{ 
        \qb\in\mathbb{R}_+
        \left|
        \begin{array}[]{l}
            \qb = \qb_1 \; \mathit{if} \; \rv_x < (\nb+1)\cm,
            \\
            \mathit{or} \; \qb = \qb_2 \; \mathit{if} \; (\nb+1)\cm \leq \rv_x \leq \min\left( (\nb\na+\nb+1)\cm,\zeta_1 \right),
            \\
            \mathit{or} \; \qb = \qb_3 \; \mathit{if} \; (\nb\na+\nb+1)\cm < \rv_x \leq \zeta_2,
            \\
            \mathit{or} \; \qb = \qb_4 \; \mathit{if} \; \rv_x \geq
            \begin{cases}    
                \na\ka + \frac{\na(\nb+1)}{2(\sqrt{\na+1}-1)}(\cm+\ka), & \mathrm{if} \; (\sqrt{\na+1}-1)^2\cm < \na\ka,
                \\
                \na\ka + (\nb+1)\left(\cm + \sqrt{\na\ka\cm}\right), & \mathrm{otherwise}. 
            \end{cases}
        \end{array}
        \right.
        \right\}.
    \end{align*}
    where
    \begin{align*}
        \zeta_1 &:=
        \nb\na\ka + (\nb+1)\cm + 2\nb\sqrt{\na\ka\cm},
        \\
        \zeta_2 &:=
        (\nb\na+\nb+1)\cm + \frac{(\nb+1)\sqrt{\na+1}}{2(\sqrt{\na+1}-1)}\left(\na\ka-\left( \sqrt{\na+1}-1 \right)^2\cm \right).
    \end{align*}
    Moreover, for each $\qb\in\Qb$,
    \begin{align*}
        y_j(\mathbf{0},\qb\mathbf{1}) = 0 & \iff \qb = \qb_1 \; \mathit{or} \; \qb_2,
        \\
        0 < y_j(\mathbf{0},\qb\mathbf{1}) < \ka & \iff \qb = \qb_3,
        \\
        y_j(\mathbf{0},\qb\mathbf{1}) = \ka & \iff \qb = \qb_4.
    \end{align*}
\end{theorem}

    \proof{Proof.}
    The result is obtained by substituting $\qf = 0$ into Proposition~\ref{prop:BaseloadReaction} and simplifying the inequalities in $\Qb$. For the case of $\qb = \qb_1$, we have
    \begin{align*}
        \qf - \cm < -\frac{\rv_x}{\nb+1}
        \iff
        \rv_x < (\nb+1)\cm.
    \end{align*}
    For the case of $\qb = \qb_2$, we have
    \begin{align*}
        & -\frac{\rv_x}{\nb+1} \leq \qf - \cm \leq \min\left( -\frac{\rv_x}{\nb\na+\nb+1}, \eta_4 \right)
        \\ \iff
        & -\frac{\rv_x}{\nb+1} \leq -\cm \; \mathrm{and} \; -\cm \leq -\frac{\rv_x}{\nb\na+\nb+1} \; \mathrm{and} \; -\cm \leq \eta_4
        \\ \iff
        & (\nb+1)\cm \leq \rv_x \; \mathrm{and} \; \rv_x \leq (\nb\na+\nb+1)\cm \; \mathrm{and} \; \rv_x \leq \zeta_1
        \\ \iff
        & (\nb+1)\cm \leq \rv_x \leq \min\left( (\nb\na+\nb+1)\cm,\zeta_1 \right),
    \end{align*}
    where the second equivalence is due to the fact that $-\cm \leq \eta_4 \iff \rv_x \leq \zeta_1$. 
    For the case of $\qb = \qb_3$, we have
    \begin{align*}
        & -\frac{\rv_x}{\nb\na+\nb+1} < \qf - \cm \leq \max\left( \eta_3,\ka-(\rv_x-\na\ka) \right)
        \\ \iff
        & -\frac{\rv_x}{\nb\na+\nb+1} < -\cm \leq 
        \begin{cases}
            \ka-(\rv_x-\na\ka), & \mathrm{if} \; \rv_x \leq \na\ka,
            \\
            \eta_3, & \mathrm{if} \; \rv_x > \na\ka,
        \end{cases}
        \\ \iff
        & (\nb\na+\nb+1)\cm < \rv_x \leq
        \begin{cases}
            \na\ka, & \mathrm{if} \; \rv_x \leq \na\ka,
            \\
            \zeta_2, & \mathrm{if} \; \rv_x > \na\ka,
        \end{cases}
        \\ \iff
        & (\nb\na+\nb+1)\cm < \rv_x \leq \zeta_2.
    \end{align*}
    The first equivalence is due to the fact that $\ka - (\rv_x-\na\ka) \geq \eta_3 \iff \rv_x \leq \na\ka$. The second equivalence is due to the fact that $\cm \geq 0$ and $\ka > 0$.  
    For the case of $\qb = \qb_4$, we have
    \begin{align*}
        \qf - \cm \geq 
        \begin{cases}
            \ka - (\rv_x - \na\ka), & \mathrm{if} \; \rv_x < \na\ka,
            \\
            \eta_2, & \mathrm{if} \; \rv_x \geq \left( 1+\frac{(\nb+1)\sqrt{\na+1}}{(\sqrt{\na+1}-1)^2} \right)\na\ka,
            \\
            \eta_1, & \mathrm{otherwise}.
        \end{cases}
    \end{align*}
    Suppose $\rv_x < \na\ka$. Then, the above inequality implies that $-\cm \geq \ka-(\rv_x-\na\ka) \implies \rv_x \geq \cm+(\na+1)\ka > \na\ka$, which is a contradiction. Henceforth, we assume that $\rv_x \geq \na\ka$, and obtain
    \begin{align*}
        & -\cm \geq 
        \begin{cases}
            \eta_2, & \mathrm{if} \; \rv_x \geq \left( 1+\frac{(\nb+1)\sqrt{\na+1}}{(\sqrt{\na+1}-1)^2} \right)\na\ka,
            \\
            \eta_1, & \mathrm{if} \; \na\ka \leq \rv_x < \left( 1+\frac{(\nb+1)\sqrt{\na+1}}{(\sqrt{\na+1}-1)^2} \right)\na\ka,
        \end{cases}
        \\ \iff
        & \rv_x \geq 
        \begin{cases}    
            \na\ka + (\nb+1)\left( \cm + \sqrt{\na\ka\cm} \right), & \mathrm{if} \; \rv_x \geq \left( 1+\frac{(\nb+1)\sqrt{\na+1}}{(\sqrt{\na+1}-1)^2} \right)\na\ka,
            \\
            \na\ka + \frac{\na(\nb+1)}{2(\sqrt{\na+1}-1)}(\cm+\ka), & \mathrm{if} \; \na\ka \leq \rv_x < \left( 1+\frac{(\nb+1)\sqrt{\na+1}}{(\sqrt{\na+1}-1)^2} \right)\na\ka,
        \end{cases}
        \\ \iff
        & \rv_x \geq 
        \begin{cases}    
            \na\ka + \frac{\na(\nb+1)}{2(\sqrt{\na+1}-1)}(\cm+\ka), & \mathrm{if} \; (\sqrt{\na+1}-1)^2\cm < \na\ka,
            \\
            \na\ka + (\nb+1)\left( \cm+\sqrt{\na\ka\cm} \right), & \mathrm{otherwise}. 
        \end{cases}
    \end{align*}
    The last equivalence is due to the fact that $\frac{\na(\nb+1)}{2(\sqrt{\na+1}-1)}(\cm+\ka) < \frac{(\nb+1)\sqrt{\na+1}}{(\sqrt{\na+1}-1)^2}\na\ka \iff (\sqrt{\na+1}-1)^2\cm < \na\ka$.
    \endproof

\section{Proofs of Structural Results}
\label{sec:StructuralProofs}

\subsection{Proof of Lemma~\ref{lem:struct-follower}}
    From Proposition~\ref{prop:PeakerReaction}, note that $0 < \qa_j(\qf\mathbf{1},\qb\mathbf{1}) < \ka$ if and only if $\qf = \frac{\na-1}{\na^2+1}\left( \rv_x-\cm-\nb\qb \right)$. Substituting into the follower productions from Proposition~\ref{prop:SpotEquilibrium2} gives
    \begin{align*}
        \qa_j(\qf\mathbf{1},\qb\mathbf{1})
        =
        \frac{\na}{\na^2+1}\xi,
    \end{align*}
    which is strictly increasing in $\xi$. Since $\xi \leq \frac{(\na^2+1)(\na-1)}{\na^2-2\sqrt{\na}+1}\ka$, we obtain
    \begin{align*}
        \bar{\qa}
        &=\left( 1-\frac{\na-2\sqrt{\na}+1}{\na^2-2\sqrt{\na}+1} \right)\ka
        \\
        &\geq
        \left( 1-\frac{\na+1}{\na^2-2\sqrt{\na}} \right)\ka
        \\
        &\geq
        \left( 1-\frac{1}{\na}\frac{\na+1}{\na}\frac{\na^2}{\na^2-2\sqrt{\na}} \right)\ka,
    \end{align*}
    which gives the first claim. 
    
    Next, from Proposition~\ref{prop:PeakerReaction}, $\ubar{\xi}$ and $\bar{\xi}$ are given by
    \begin{align*}
        \ubar{\xi}
        &=
        (\na+1)\ka,
        \\
        \bar{\xi}
        &=
        \frac{(\na^2+1)(\na-1)}{\na^2-2\sqrt{\na}+1}\ka,
    \end{align*}
    from which we obtain
    \begin{align*}
        \frac{\ubar{\xi}-\bar{\xi}}{\ubar{\xi}}
        &=
        \frac{2(\na^2-\na\sqrt{\na}-\sqrt{\na}+1)}{(\na^2-2\sqrt{\na}+1)(\na+1)}
        \\
        &\leq
        \frac{2(\na^2+1)}{\na(\na^2-2\sqrt{\na})}
        \\
        &=
        \frac{2}{\na}\frac{\na^2+1}{\na^2}\frac{\na^2}{\na^2-2\sqrt{\na}},
    \end{align*}
    which gives the rest of the second claim.

\subsection{Proof of Lemma~\ref{lem:struct-leader}}
    From Proposition~\ref{prop:BaseloadReaction}, $\ubar{\qf}$ and $\bar{\qf}$ are given by
    \begin{align*}
        \ubar{\qf}
        &=
        \cm - \frac{\rv_x}{\nb+1},
        \\
        \bar{\qf}
        &=
        \cm + \min\left( -\frac{\rv_x}{\nb\na+\nb+1},\eta_4 \right).
    \end{align*}
    The first claim follows from Proposition~\ref{prop:BaseloadReaction}. Next,
    \begin{align*}
        \bar{\qf}-\ubar{\qf}
        &=
        \frac{\rv_x}{\nb+1} + \min\left( -\frac{\rv_x}{\nb\na+\nb+1}, \eta_4 \right)
        \\
        & \leq
        \frac{\rv_x}{\nb+1} - \frac{\rv_x}{\nb\na+\nb+1}
        \\
        &=
        \frac{\nb\na\rv_x}{(\nb+1)(\nb\na+\nb+1)}
        \\
        &\leq
        \frac{\nb\na\rv_x}{\nb^2(\na+1)}
        \\
        &=
        \frac{\rv_x}{\nb}\frac{\na}{\na+1},
    \end{align*}
    which gives the second claim.

\subsection{Proof of Lemma~\ref{lem:struct-leader-2}}
    From Proposition~\ref{prop:BaseloadReaction}, note that $0 < \qa_j(\qf\mathbf{1},\qb\mathbf{1}) < \ka \implies \qb=\qb_3$. Substituting into the follower productions from Proposition~\ref{prop:SpotEquilibrium2} gives
    \begin{align*}
        \qa_j(\qf\mathbf{1},\qb\mathbf{1})
        =
        \left[ \frac{1}{\na+1}\left( \rv_x+(\qf-\cm)-\frac{\nb}{\nb+1}(\rv_x+\na(\cm-\qf)) \right) \right]_0^\ka,
    \end{align*}
    which is strictly increasing in $\qf$. Note that $\qb = \qb_3$ is a reaction if and only if
    \begin{align*}
        -\frac{\rv_x}{\nb\na+\nb+1} < \qf-\cm \leq \max\left( \eta_3,\ka-(\rv_x-\na\ka) \right) \iff -\frac{\rv_x}{\nb\na+\nb+1} < \qf-\cm \leq \eta_3,
    \end{align*}
    where we used the fact that $\rv_x > \na\ka \implies \eta_3 \geq \ka-(\rv_x-\na\ka)$. Since 
    \begin{align*}
    \rv_x \leq \na\ka\left( 1+\frac{(\nb+1)\sqrt{\na+1}}{(\sqrt{\na+1}-1)^2}+\frac{(\nb-1)\sqrt{\na+1}}{\sqrt{\na+1}-1} \right) \implies -\frac{\rv_x}{\na\nb+\nb+1} \leq \eta_3,
    \end{align*}
    we infer the case for $\bar{\qa} = 0$. Otherwise, substituting for $\eta_3$ gives 
    \begin{align*}
        \bar{\qa}
        &=
        \qa_j\left( (\cm+\eta_3)\mathbf{1},\qb_3\mathbf{1} \right)
        \\
        &=
        \ka + \frac{\rv_x-\na\ka}{\na+1}\frac{1}{\nb+1}\left( \frac{2(\nb+1)(\na+1)-(\nb+1)(\na+2)\sqrt{\na+1}}{\na(2+(\nb-1)\sqrt{\na+1}} \right)
        \\
        &\geq
        \ka - \frac{\rv_x-\na\ka}{\na} \left( \frac{\na+2}{(\na+1)(\nb-1)} \right),
    \end{align*}
    from which we obtain the first claim. From Proposition~\ref{prop:BaseloadReaction}, we infer that $\bar{\qf}=\eta_3$ and $\ubar{\qf}=\eta_1$ when $\rv_x \leq \na\ka\left( 1+\frac{(\nb+1)\sqrt{\na+1}}{(\sqrt{\na+1}-1)^2} \right)$. Therefore, we obtain
    \begin{align*}
        \bar{\qf} - \ubar{\qf}
        &= \frac{\rv_x-\na\ka}{\na}\left( 2(\sqrt{\na+1}-1 \right)\left( \frac{2+(\nb-1)\sqrt{\na+1}-(\nb+1)}{(\nb+1)(2+(\nb-1)\sqrt{\na+1}} \right)
        \\
        &\leq 2 \left( \frac{\rv_x-\na\ka}{\na} \right) \left( \frac{\sqrt{\na+1}-1}{\nb-1} \right)
        \\
        &\leq \frac{\rv_x-\na\ka}{\nb\sqrt{\na}} 2 \left( \frac{\sqrt{\na+1}}{\sqrt{\na}}\frac{\nb}{\nb-1} \right),
    \end{align*}
    which gives the second claim.

\subsection{Proof of Lemma~\ref{lem:struct-market}}
    From Theorem~\ref{prop:MarketEquilibrium}, note that $0 < \qa_j(\qf\mathbf{1},\qb\mathbf{1}) < \ka \implies (\qf,\qb)\in Q_3$. Substituting into the follower productions from Proposition~\ref{prop:SpotEquilibrium2} gives
    \begin{align*}
        \qa_j(\qf\mathbf{1},\qb\mathbf{1})
        =
        \frac{\na}{\na^2+\na\nb+\nb+1} \rv_x,
    \end{align*}
    which is strictly increasing in $\rv_x$. Since $\rv_x\leq\zeta_2$, it follows that
    \begin{align*}
        \bar{\qa}
        &=
        \frac{\na}{\na^2+\na\nb+\nb+1}\zeta_2
        \\
        &=
        \left( 1-\frac{\na+2-2\sqrt{\na+1}}{\na^2+\na+2-2\sqrt{\na+1}} \right)\ka
        \\
        &\geq
        \left( 1-\frac{\na}{\na^2+\na} \right)\ka
        \\
        &=
        \left( 1-\frac{1}{\na}\frac{\na}{\na+1} \right)\ka,
    \end{align*}
    from which we obtain the first claim. 
    
    Next, from Theorem~\ref{prop:MarketEquilibrium}, we infer that $\bar{\rv}_x = \zeta_2$ and $\ubar{\rv} = (\nb+\na+1)\ka$. It is easy to show that $\zeta_2 < (\nb+\na+1)\ka \iff \nb < \na\sqrt{\na+1}-1$. Moreover, 
    \begin{align*}
        (\nb+\na+1)\ka - \zeta_2
        &=
        \frac{2}{\na^2+(\sqrt{\na+1}-1)^2}\left( (\na^2+\na+\nb+1)-(\nb+\na+1)\sqrt{\na+1} \right)\ka.
    \end{align*}
    Therefore, if $\ubar{\rv}_x \leq \bar{\rv}_x$, then
    \begin{align*}
        \frac{\ubar{\rv}_x-\bar{\rv}_x}{\ubar{\rv}_x}
        &=
        \frac{2}{\na^2+(\sqrt{\na+1}-1)^2}\left( 1+\frac{\na^2}{\nb+\na+1}-\sqrt{\na+1} \right)
        \\
        &\leq
        \frac{2\na}{\na^2+(\sqrt{\na+1}-1)^2}
        \\
        &\leq
        \frac{2\na}{\na^2},
    \end{align*}
    from which we obtain the first part of the second claim. If $\ubar{\rv}_x \geq \bar{\rv}_x$, then
    \begin{align*}
        \frac{\bar{\rv}_x-\ubar{\rv}_x}{\ubar{\rv}_x}
        &=
        \frac{2}{\na^2+(\sqrt{\na+1}-1)^2}\left( \sqrt{\na+1}-1-\frac{\na^2}{\nb+\na+1} \right)
        \\
        &\leq
        \frac{2}{\na^2+(\sqrt{\na+1}-1)^2}\left( \sqrt{\na+1} \right)
        \\
        &\leq
        \frac{2}{\na^2}\sqrt{\na+2}
        \\
        &\leq
        \frac{2}{\na\sqrt{\na}}\sqrt{\frac{\na+2}{\na}},
    \end{align*}
    from which we obtain the rest of the second claim.

\subsection{Proof of Lemma~\ref{lem:struct-market-3}}
    From Theorem~\ref{prop:MarketEquilibrium}, note that $0< \qa_j(\qf\mathbf{1},\qb\mathbf{1})<\ka\iff(\qf,\qb)\in Q_3$. Substituting into the follower productions from Proposition~\ref{prop:SpotEquilibrium2} gives
    \begin{align*}
        \qa_j(\qf\mathbf{1},\qb\mathbf{1})
        =
        \frac{\na}{\na^2+\nb\na+\nb+1}\left( \rv_x-(\nb\na+\nb+1)\cm \right),
    \end{align*}
    which is strictly increasing in $\rv_x$. Since $\rv_x \leq \zeta_2$, it follows that
    \begin{align*}
        \bar{\qa}
        &=
        \frac{\na}{\na^2+\nb\na+\nb+1}\left( \zeta_2-(\nb\na+\nb+1)\cm \right)
        \\
        &=
        \frac{\na^2}{\na^2+(\sqrt{\na+1}-1)^2}\left( \ka-\frac{(\sqrt{\na+1}-1)^2}{\na}\cm \right)
        \\
        &\geq
        \left( 1-\frac{\na+2-2\sqrt{\na+1}}{\na^2} \right)\left( \ka-\frac{(\sqrt{\na+1}-1)^2}{\na}\cm \right)
        \\
        &\geq
        \left( 1-\frac{1}{\na} \right)\left( \ka-\frac{(\sqrt{\na+1}-1)^2}{\na}\cm \right),
    \end{align*}
    from which we obtain the first claim.

    Next, from Theorem~\ref{prop:MarketEquilibrium}, we infer that if $(\sqrt{\na+1}-1)^2\cm<\na\ka$, then $\bar{\rv}_x = \zeta_2$ and $\ubar{\rv}_x = (\nb+\na+1)\ka$, and it is straightforward to show that $\zeta_2<(\nb+\na+1)\ka$ if and only if the first case in~\eqref{eq:struct-market-3-cond} holds. Otherwise, then $\bar{\rv}_x = \zeta_1$ and $\ubar{\rv}_x = (\nb+\na+1)\ka$, and it is straightforward to show that $\zeta_1<(\nb+\na+1)\ka$ if and only if the second case in~\eqref{eq:struct-market-3-cond} holds.

\subsection{Proof of Lemma~\ref{lem:struct-stack}}
    From Theorem~\ref{prop:StackelbergEquilibrium}, note that $0 < \qa_j(\mathbf{0},\qb\mathbf{1}) < \ka \iff \qb=\qb_3$. Substituting into the follower productions from Proposition~\ref{prop:SpotEquilibrium2} gives
    \begin{align*}
        \qa_j(\mathbf{0},\qb\mathbf{1})
        =\frac{1}{(\na+1)(\nb+1)}\left( \rv_x - (\nb\na+\nb+1)\cm \right),
    \end{align*}
    which is strictly increasing in $\rv_x$. Since $\rv_x \leq \zeta_2$, it follows that
    \begin{align*}
        \bar{\qa}
        &=
        \frac{1}{(\na+1)(\nb+1)}\left( \zeta_2-(\nb\na+\nb+1)\cm \right)
        \\
        &=\left( 1+\frac{1}{\sqrt{\na+1}} \right)\frac{\ka}{2},
    \end{align*}
    from which we obtain the first claim.

    Next, from Theorem~\ref{prop:StackelbergEquilibrium}, we infer that $\bar{\rv}_x = \zeta_2$ and $\ubar{\rv}_x = \na\ka+\frac{\na(\nb+1)}{2(\sqrt{\na+1}-1)}(\cm+\ka)$. It is straightforward to show that $\bar{\rv}_x \geq \ubar{\rv}_x$.

\subsection{Proof of Lemma~\ref{lem:struct-stack-3}}
    From Theorem~\ref{prop:StackelbergEquilibrium}, note that $0<\qa_j(\mathbf{0},\qb\mathbf{1})<\ka\iff\qb=\qb_3$. Substituting into the follower productions from Proposition~\ref{prop:SpotEquilibrium2} gives
    \begin{align*}
        \qa_j(\mathbf{0},\qb\mathbf{1})
        =
        \frac{1}{(\na+1)(\nb+1)}\left( \rv_x-(\nb\na+\nb+1)\cm \right),
    \end{align*}
    which is strictly increasing in $\rv_x$. Since $\rv_x \leq \zeta_2$, it follows that
    \begin{align*}
        \bar{\qa}
        &= \frac{1}{(\na+1)(\nb+1)}\left( \zeta_2-(\nb\na+\nb+1)\cm \right)
        \\
        &= \left( 1+\frac{1}{\sqrt{\na+1}} \right)\frac{1}{2}\left( \ka-\frac{(\sqrt{\na+1}-1)^2}{\na}\cm \right),
    \end{align*}
    from which we obtain the first claim.

    Next, from Theorem~\ref{prop:StackelbergEquilibrium}, we infer that, if $(\sqrt{\na+1}-1)^2\cm < \na\ka$, then $\bar{\rv}_x = \zeta_2$ and $\ubar{\rv}_x = \na\ka+\frac{\na(\nb+1)}{2(\sqrt{\na+1}-1)}(\cm+\ka)$, and it is straightforward to show that $\bar{\rv}_x \geq \ubar{\rv}_x$. Otherwise, then $\bar{\rv}_x = \zeta_1$ and $\ubar{\rv}_x = \na\ka + (\nb+1)(\cm+\sqrt{\na\ka\cm})$, and it is straightforward to show that $\bar{\rv}_x \geq \ubar{\rv}_x$. 

\subsection{Proof of Lemma~\ref{lem:struct-compare}}
    The proof proceeds in three steps. In step 1, we compute an equilibria with the smallest (resp. largest) market production in the forward (resp. Stackelberg) market. In step 2, we compute an equilibria with the smallest (resp. largest) social welfare in the forward (resp. Stackelberg) market. In step 3, we show that the worst case ratios of productions and efficiencies are both strictly increasing in $\rv_x$. The bounds in the lemma are obtained by evaluating those ratios at $\rv_x = \bar{\rv}_x$.

    \emph{Step 1:} We compute an equilibria with the smallest (resp. largest) market production in the forward (resp. Stackelberg) market. First, we tackle the forward market. Substituting $\cm = 0$ into Theorem~\ref{prop:MarketEquilibrium}, we infer that $(\qf,\qb)\in Q$ if and only if $(\qf,\qb)\in Q_3$ or $(\qf,\qb)\in Q_4$. By substituting into Theorem~\ref{prop:SpotEquilibrium2}, and using the fact that $\qa_j(\qf\mathbf{1},\qb\mathbf{1})=0$ for all $(\qf,\qb)\in Q_4$, we obtain the following market productions:
    \begin{align*}
        \nb\qb + \na\qa_j(\qf\mathbf{1},\qb\mathbf{1})
        &=
        \begin{cases}
            \frac{1}{\na^2+\nb\na+\nb+1}\left( \na^2+\nb\na+\nb \right)\rv_x, & \mathrm{if} \; (\qf,\qb)\in Q_3,
            \\
            \frac{1}{\nb+1}\left( \nb\rv_x+\na\ka \right), & \mathrm{if} \; (\qf,\qb)\in Q_4.
        \end{cases}
    \end{align*}
    Note that
    \begin{align*}
        &\frac{1}{\nb+1}\left[ \nb\rv_x+\na\ka \right]
        \\
        &=\frac{1}{\na^2+\nb\na+\nb+1}\left[ (\na^2+\nb\na+\nb)\rv_x + \frac{-\na^2-\nb\na}{\nb+1}\rv_x + \frac{\na(\na^2+\nb\na+\nb+1)}{\nb+1}\ka\right]
        \\
        &\leq\frac{1}{\na^2+\nb\na+\nb+1}\left[ (\na^2+\nb\na+\nb)\rv_x + \frac{\na(\nb+1)(1-\na-\nb)}{\nb+1}\ka \right]
        \\
        &\leq\frac{1}{\na^2+\nb\na+\nb+1}\left( \na^2+\nb\na+\nb \right)\rv_x,
    \end{align*}
    where the first inequality is due to the fact that $\rv_x \geq \ubar{\rv}_x$ and the second inequality is due to the fact that $\nb\geq1$, $\na\geq2$, and $\ka > 0$. Therefore, we infer that the smallest equilibrium production in the forward market is given by
    \begin{align*}
        \qa_F
        &= \ka,
        \\
        \qb_F
        &= \frac{1}{\nb+1}(\rv_x-\na\ka).
    \end{align*}

    Next, we tackle the Stackelberg market. Substituting $\cm = 0$ into Theorem~\ref{prop:StackelbergEquilibrium}, we infer that $(0,\qb_s)\in \Qb(0)$ if and only if $\qb_s = \qb_3$ or $\qb_s = \qb_4$. Suppse
    \begin{align}
        \rv_x < \na\ka + \frac{\na(\nb+1)}{2(\sqrt{\na+1}-1)}\ka.
        \label{eq:struct-stack-4-pf-1}
    \end{align}
    Then, from Theorem~\ref{prop:StackelbergEquilibrium}, we conclude that $\qb_s = \qb_3$ is the only Stackelberg equilibrium, and hence it is also the equilibrium with the largest market production. Suppose, instead, that~\eqref{eq:struct-stack-4-pf-1} does not hold. By substituting into Proposition~\ref{prop:SpotEquilibrium2}, and using the fact that $\qa_j(\mathbf{0},\qb_4\mathbf{1}) = 0$, we obtain the following market productions:
    \begin{align*}
        \nb\qb_s + \na\qa_j(\mathbf{0},\qb_s\mathbf{1})
        =
        \begin{cases}
            \frac{\nb\na+\nb+\na}{(\nb+1)(\na+1)}\rv_x, & \mathrm{if} \; \qb_s = \qb_3,
            \\
            \frac{1}{\nb+1}\left( \nb\rv_x+\na\ka \right), & \mathrm{if} \; \qb_s = \qb_4.
        \end{cases}
    \end{align*}
    Note that
    \begin{align*}
        &\frac{\nb\na+\nb+\na}{(\nb+1)(\na+1)}\rv_x
        \\
        &=\frac{1}{\nb+1}\left[ \nb\rv_x + \frac{\na}{\na+1}\rv_x \right]
        \\
        &\geq\frac{1}{\nb+1}\left[ \nb\rv_x + \frac{\na}{\na+1}\left( \na+\frac{\na(\nb+1)}{2(\sqrt{\na+1}-1)} \right)\ka \right]
        \\
        &\geq\frac{1}{\nb+1}\left[ \nb\rv_x + \frac{\na}{\na+1}\left( \na+1 \right)\ka \right]
        \\
        &=\frac{1}{\na+1}\left[ \nb\rv_x + \na\ka \right],
    \end{align*}
    where the first inequality is due to the fact that~\eqref{eq:struct-stack-4-pf-1} does not hold and the second inequality is due to the fact that $\nb \geq 1$ and $\na \geq 2$. Therefore, we infer that the largest equilibrium production in the Stackelberg market is given by
    \begin{align*}
        \qa_S &= \frac{1}{\na+1}\left( \rv_x - \nb\qb_s \right),
        \\
        \qb_S &= \frac{1}{\nb+1}\rv_x.
    \end{align*}

    \emph{Step 2:} We compute the equilibria with the smallest (resp. largest) social welfare in the forward (resp. Stackelberg) market. Substituting the demand function into the social welfare gives
    \begin{align*}
        \mathsf{SW}(\qa,\qb)
        &=\beta\left( \rv_x\left( \nb\qb+\na\qa \right)-\cm\na\qa-\frac{1}{2}\left( \nb\qb+\na\qa \right)^2 \right)
        \\
        &=
        \beta\left( \rv_x\left( \nb\qb+\na\qa \right)-\frac{1}{2}\left( \nb\qb+\na\qa \right)^2 \right),
    \end{align*}
    where the second equality is obtained by substituting $\cm = 0$. Given any two equilibrium productions $(\qa,\qb)$ and $(\qa',\qb')$, we have
    \begin{align*}
        & \mathsf{SW}(\qa,\qb) \geq \mathsf{SW}(\qa',\qb')
        \\ \iff &
        \rv_x(\nb\qb+\na\qa) - \frac{1}{2}(\nb\qb+\na\qa)^2 \geq \rv_x(\nb\qb'+\na\qa') - \frac{1}{2}(\nb\qb'+\na\qa')^2
        \\ \iff &
        \frac{1}{2}\left( (\nb\qb+\na\qa)-(\nb\qb'+\na\qa') \right)\left( \rv_x - (\nb\qb+\na\qa)+ \rv_x - (\nb\qb'+\na\qa') \right) \geq 0
        \\ \iff &
        \frac{1}{2}\left( (\nb\qb+\na\qa)-(\nb\qb'+\na\qa') \right)\left( \frac{1}{\beta}\left( P(\nb\qb+\na\qa)-\cb \right) + \frac{1}{\beta}\left( P(\nb\qb'+\na\qa')-\cb \right) \right)\geq 0
        \\ \iff & 
        \nb\qb + \na\qa \geq \nb\qb' + \na\qa',
    \end{align*}
    where the last equivalence follows from the fact that, since $(\qa,\qb)$ and $(\qa',\qb')$ are equilibrium productions, the profit margins $P(\nb\qb+\na\qa)-\cb > 0$ and $P(\nb\qb'+\na\qa')-\cb > 0$. Therefore, the equilibria with the smallest (resp. largest) social welfare in the forward (resp. Stackelberg) market are those with the smallest (resp. largest) market productions, which were obtained in step 1.

    \emph{Step 3:} We show that the worst-case ratios of productions and social welfares are strictly increasing in $\rv_x$. From step 1, the ratio of productions is bounded from above by
    \begin{align*}
        r_P := \frac{\nb\qb_S+\na\qa_S}{\nb\qb_F+\na\qa_F}.
    \end{align*}
    Taking derivatives gives
    \begin{align*}
        \frac{\partial r_P}{\partial \rv_x}
        &=\frac{(\nb\qb_F+\na\qa_F)\left( \nb\frac{\partial\qb_S}{\partial\rv_x}+\na\frac{\partial\qa_S}{\partial\rv_x} \right) - (\nb\qb_S+\na\qa_S)\left( \nb\frac{\partial\qb_F}{\partial\rv_x} + \na\frac{\partial\qa_F}{\partial\rv_x} \right)}{(\nb\qb_F+\na\qa_F)^2}
        \\
        &=\frac{\frac{\na(\nb\na+\nb+\na)\ka}{(\nb+1)^2(\na+1)}}{(\nb\qb_F+\na\qa_F)^2}
        \\
        &>0.
    \end{align*}
    Next, the ratio of social welfares is bounded from above by
    \begin{align*}
        r_W := \frac{\mathsf{SW}(\qa_S,\qb_S)}{\mathsf{SW}(\qa_F,\qb_F)}.
    \end{align*}
    Taking derivatives gives
    \begin{align*}
        \frac{\partial r_W}{\partial \rv_x}
        &=\frac{\mathsf{SW}(\qa_F,\qb_F)\frac{\partial\mathsf{SW}(\qa_S,\qb_S)}{\partial\rv_x}-\mathsf{SW}(\qa_S,\qb_S)\frac{\partial\mathsf{SW}(\qa_F,\qb_F)}{\partial\rv_x}}{\mathsf{SW}(\qa_F,\qb_F)^2}
        \\
        &=
        \frac{\frac{\beta}{2(\nb+1)^4(\na+1)^2} 
            (\nb\na+\nb+\na)(\nb\na+\nb+\na+2)(\rv_x-\na\ka)\na\ka\rv_x        
        }{\mathsf{SW}(\qa_F,\qb_F)^2}
        \\
        &>0,
    \end{align*}
    where the inequality is due to $\rv_x \geq \ubar{\rv}_x > \na\ka$. Therefore, $r_P$ and $r_W$ are both strictly increasing in $\rv_x$ over $\left[ \ubar{\rv}_x,\bar{\rv}_x \right]$. By substituting $\rv_x = \bar{\rv}_x$ into $r_P$ and $r_W$, we obtain the desired result.  



\end{document}